\documentclass[preprint]{elsarticle}
\usepackage{lineno,hyperref}
\usepackage{latexsym}
\usepackage{graphicx}
\usepackage{multicol,multirow}
\usepackage{amsmath,amssymb,amsfonts}
\usepackage{mathrsfs}
\usepackage{amsthm}
\usepackage{subcaption}
\usepackage{upgreek}
\usepackage{xcolor}

\newtheorem*{remark}{Remark}
\DeclareMathOperator*{\argmin}{arg\,min}

\usepackage{xcolor}

\renewcommand{\textcolor}[2]{#2} 
\journal{Computer Methods in Applied Mechanics and Engineering}

\begin{document}

\begin{frontmatter}

\title{Deep autoencoders for physics-constrained data-driven nonlinear materials modeling}

\author[1]{Xiaolong He\fnref{first_author}}
\author[2]{Qizhi He\fnref{first_author}\corref{corresponding_author}}
\fntext[first_author]{These authors contributed equally to this work.}

\author[1]{Jiun-Shyan Chen\fnref{first_author}\corref{corresponding_author}}
\cortext[corresponding_author]{Corresponding author}
\ead{qizhi.he@pnnl.gov (Qizhi He); js-chen@ucsd.edu}

\address[1]{Department of Structural Engineering, University of California, San Diego, La Jolla, CA, 92093, USA}
\address[2]{Physical and Computational Sciences Directorate, Pacific Northwest National Laboratory, Richland, WA, 99354, USA}






\begin{abstract}
Physics-constrained data-driven computing is an emerging computational paradigm that allows simulation of complex materials directly based on material database and bypass the classical constitutive model construction. 
However, it remains difficult to deal with high-dimensional applications and extrapolative generalization. 
This paper introduces deep learning techniques under the data-driven framework to address these fundamental issues in nonlinear materials modeling. 
To this end, an autoencoder neural network architecture is introduced to learn the underlying low-dimensional representation (embedding) of the given material database. The offline trained autoencoder and the discovered embedding space are then incorporated in the online data-driven computation such that the search of optimal material state from database can be performed on a low-dimensional space, aiming to enhance the robustness and predictability with projected material data.
To ensure numerical stability and representative constitutive manifold, a convexity-preserving interpolation scheme tailored to the proposed autoencoder-based data-driven solver is proposed for constructing the material state.
In this study, the applicability of the proposed approach is demonstrated by modeling nonlinear biological tissues.
A parametric study on data noise, data size and sparsity, training initialization, and model architectures, is also conducted to examine the robustness and convergence property of the proposed approach.
\end{abstract}
\end{frontmatter}


\section{Introduction}\label{sec1}
Constitutive modeling is traditionally based on constitutive or material laws to describe the explicit relationship among strain, stress, and state variables based on experimental observations, physical hypothesis, and mathematical simplifications.
However, the phenomenological modeling process inevitably introduces errors due to limited data and mathematical assumptions in model parameter calibration.
Moreover, constitutive laws rely on pre-defined functions and often lack generality to capture full aspects of material behaviors \cite{HeJBM2020}.

With advancements in computing power and proliferation of digital data \cite{goodfellow2016deep}, machine learning (ML) based data-driven approaches, such as artificial neural networks (NNs), have emerged as a promising alternative for constitutive modeling due to their ability in extracting features and complex patterns in data \cite{bock2019review}. 
This type of approaches usually requires parameterization of machine learning models with given constitutive data \textit{a priori}, and thus, they are classified as model-based approaches.
For example, NNs have been applied to modeling a variety of materials, including concrete materials \cite{ghaboussi1991knowledge}, hyper-elastic materials \cite{Shen2005}, viscoplastic material of steel \cite{furukawa1998implicit}, and homogenized properties of composite structures \cite{lefik2009artificial}.
Accordingly, the NN-based constitutive models were integrated into finite element codes to predict path- or rate-dependent materials behaviors \cite{lefik2003artificial,hashash2004numerical,jung2006neural,Stoffel2019}.
Recently, deep neural networks (DNNs) with special mechanistic architectures, such as recurrent neural networks (RNNs) haven been applied to path-dependent materials \cite{wang2018multiscale,mozaffar2019deep,heider2020so}. In \cite{heider2020so}, for example, graph neural networks (GNNs) based representation of topology information of crystal microstructures (e.g., Euler angles) were introduced for learning anisotropic elasto-plastic responses.
Besides, the deep material network approach was proposed in \cite{liu2019deep,Liu2019a} to simulate multiscale heterogeneous materials, in which the NNs were constructed based on hierarchical mechanistic building blocks derived from analytical two-phase elasticity model and homogenization theory.

Another strategy in data-driven \textcolor{blue}{materials modeling} is to bypass the constitutive modeling step by formulating an optimization problem to search for the physically admissible state that satisfies equilibrium and compatibility and minimizes the distance to a material dataset \cite{kirchdoerfer2016data,ibanez2018manifold,Conti2018}.
In this data-driven approach, the search of material data at each integration point from the material dataset is determined via a distance-minimization function and is called the distance-minimizing data-driven (DMDD) computing.
This data-driven computing paradigm has been extended to dynamics \cite{kirchdoerfer2018data}, problems with geometrical nonlinearity \cite{nguyen2018data,HeJBM2020}, inelasticity \cite{eggersmann2019model}, anisotropy \cite{he2020physics}, material identification and constitutive manifold construction \cite{Ibanez2017,Leygue2018,Ayensa-Jimenez2019,Stainier2019}.
A variational framework for data-driven computing was proposed in \cite{he2019physics,Nguyen2020} to allow versatile in the employment of special approximation functions and numerical methods.

To better handle noise induced by outliers and intrinsic randomness in the experimental datasets, data-driven computing integrated with statistical models or machine learning techniques were developed, including incorporating maximum entropy estimation with a cluster analysis \cite{kirchdoerfer2017data}, regularization based on data sampling statistics \cite{ayensa2018new}, and locally convex reconstruction inspired from manifold learning for nonlinear dimensionality reduction \cite{he2019physics}.
In the local convexity data-driven (LCDD) computing proposed in \cite{he2019physics}, the convexity condition is imposed on the reconstructed material graph to avoid convergence issues that usually arise in standard data fitting approaches.
Recently, Eggersmann et al. \cite{Eggersmann} introduced tensor voting, an instance-based machine learning technique, to the entropy-based data-driven scheme \cite{kirchdoerfer2017data} to construct locally linear tangent spaces for achieving higher-order convergence in data-driven solvers.

The above mentioned (distance-minimizing) data-driven computing approaches are distinct from the NN-based constitutive modeling approaches.
The latter constructs the surrogate models of constitutive laws independent to the solution procedure of boundary value problems, whereas the former is "model free" and it incorporates the material data search into the solution procedure of boundary value problems. Thus, the former approach was also called model-free data-driven computing \cite{kirchdoerfer2017data,eggersmann2019model}.
The model-free data-driven approach circumvents the need of using material tangent during solution iteration processes, which offers another unique feature in computational mechanics.
From the perspective of machine learning algorithms, the NN approach is considered as supervised learning, and it requires pre-defined input-output functions, e.g. strain-stress laws. It has been shown that selecting the response function by NNs within the constitutive framework is non-trivial \cite{Settgast2019a,wang2018multiscale,heider2020so}.
On the other hand, the data-driven computing approaches with unsupervised learning such as clustering and manifold learning are capable of discovering the underlying data structure for the constitutive manifold
\cite{he2019physics,ibanez2018manifold}.
However, this type of data-driven approaches could encounter difficulties in high dimensional applications when material data sampling involves multidimensional and history-dependent state variables. As demonstrated in \cite{kirchdoerfer2016data,Conti2018,he2019physics}, higher data dimension results in lower convergence rate with respect to data size, and thus demands effective dimensionality reduction for improved effectiveness of data-driven computing. 
Further, the model-free data-driven schemes rely on the direct search of nearest neighbors from the material dataset, leading to limited extrapolative generalization when the distribution of data points becomes sparse.

In this work, we aim to develop a novel data-driven computing approach to overcome the curse of dimensionality and the lack of generalization in classical model-free data-driven computing approaches \cite{kirchdoerfer2016data,ibanez2018manifold,he2019physics}.
To this end, we propose to introduce a novel autoencoders based deep neural network architecture \cite{demers1993non,hinton2006reducing} under the LCDD framework \cite{he2019physics}
to achieve two major objectives: dimensionality reduction and generalization of physically meaningful constitutive manifold.
It should be emphasized that there have been various nonlinear dimensionality reduction (i.e. manifold learning) techniques developed for complex high-dimensional data \cite{scholkopf1998nonlinear,tenenbaum2000global,roweis2000nonlinear,belkin2003laplacian,van2009dimensionality}, but these methods remain challenging in out-of-sample extension (OOSE), that is to project new unseen data onto the learned low-dimensional representation space, due to the lack of explicit mapping function. While nonparametric OOSE \cite{van2009dimensionality} is an option to construct the mapping, the computational complexity increases substantially with the size of dataset.
On the other hand, owing to its deep learning architecture, an autoencoder is capable of naturally defining the mapping functions between high- and low-dimensional representation and capturing highly varying nonlinear manifold with good generalization capability \cite{van2009dimensionality,goodfellow2016deep}, as to be discussed in Section \ref{sec.AEDD_local}.

By integrating autoencoders and the discovered low-dimensional embedding into the data-driven solver, the proposed framework is referred to as auto-embedding data-driven (AEDD) computing, which can also be considered as a hybrid of the NN-based constitutive modeling and the classical model-free data-driven computing.
In this approach, the autoencoders are first trained in an offline stage to extract a representative low-dimensional manifold (embedding) of the given material data. Autoencoders also provide effective noise filtering through the data compression processes.
The trained autoencoders are then incorporated in the data-driven solver during the online computation with customized convexity-preserving reconstruction.
Hence, all operations related to distance measure, including the search of the closest material points in the dataset are performed in the learned embedding space, circumventing the difficulties resulting from high dimensionality and data noise.
To ensure numerical stability and representative constitutive manifold parameterized by the trained autoencoder networks, an efficient convexity-preserving interpolation is proposed to locally approximate the optimal material data to a given physically admissible state.
Specifically, in this work, we present two different solvers to perform locally convex reconstruction, and demonstrate the one directly providing interpolation approximation without using decoders outperforms the one that fully uses the encoder-decoder network structure.
Furthermore, it is shown that the proposed method is computationally tractable, since the additional autoencoder training is conducted offline and the online data-driven computation mainly involve lower-dimensional variables in the embedding space.

The remainder of this paper is organized as follows. The background of physics-constrained data-driven framework is first introduced in Section \ref{sec2}, including the data-driven equations in variational form and the material local solvers.
In Section \ref{sec.AEDD_local}, the basic theory of deep neural networks and autoencoders are presented. Section \ref{sec.AEDD} introduces the convexity-preserving interpolation for data-driven materials modeling with improved performance.
Finally, in Section \ref{sec.results}, the effectiveness of the proposed AEDD framework are examined, and a parametric study is conducted to investigate the effects of autoencoder architecture, data noise, data size and sparsity, and neural network initialization on AEDD's performance.
The proposed method is also applied to biological tissue modeling to demonstrate the enhanced effectiveness and generalization capability of AEDD over the other data-driven schemes.
Concluding remarks and discussions are summarized in Section \ref{sec.conclusion}.
\section{Formulations of physics-constrained data-driven nonlinear modeling}\label{sec2}

This section provides the basic equations of the physics-constrained data-driven computational framework for nonlinear solids \cite{nguyen2018data,he2019physics,HeJBM2020}, followed by a review of two material data-driven local solvers and the associated computational approaches.

\subsection{Governing equations of nonlinear mechanics}\label{sec.goveqn}
The equations governing the deformation of a solid in a domain $\Omega^X$ bounded by a Neumann boundary $\Gamma_t^X$ and a Dirichlet boundary $\Gamma_u^X$ in the undeformed configuration are given as
\begin{equation}\label{eq.2.3.1}
    \begin{cases}
    DIV(\mathbf{F}(\mathbf{u}) \cdot \mathbf{S}) + \mathbf{b} = \mathbf{0}, & \text{in} \hspace{0.1cm} \Omega^X,\\
    \mathbf{E} =\mathbf{E}(\mathbf{u}) = (\mathbf{F}^T\mathbf{F}-\mathbf{I})/2, & \text{in} \hspace{0.1cm} \Omega^X,\\
    (\mathbf{F}(\mathbf{u}) \cdot \mathbf{S}) \cdot \mathbf{N} = \mathbf{t}, & \text{on} \hspace{0.1cm} \Gamma_t^X,\\
    \mathbf{u} = \mathbf{g}, & \text{on} \hspace{0.1cm} \Gamma_u^X,
    \end{cases}
\end{equation}
where $\mathbf{u}$ is the displacement vector, $\mathbf{E}$ is the Green Lagrangian strain tensor, $\mathbf{S}$ is the second Piola-Kirchhoff (2nd-PK) stress tensor, and $DIV$ denotes the divergence operator. Without loss of generality, the governing equations (\ref{eq.2.3.1}) are defined in the reference (undeformed) configuration \cite{Belytschko:2000tz}, which is denoted by the superscript "$X$". In Eq. (\ref{eq.2.3.1}), $\mathbf{F}$ is the deformation gradient related to $\mathbf{u}$, defined as $\mathbf{F}(\mathbf{u}) = \partial (\mathbf{X}+\mathbf{u}) / \partial \mathbf{X}$, where $\mathbf{X}$ is the material coordinate, and $\mathbf{b}$, $\mathbf{N}$, $\mathbf{t}$, and $\mathbf{g}$ are the body force, the surface normal on $\Gamma_t^X$, the traction on $\Gamma_t^X$, and the prescribed displacement on $\Gamma_u^X$, respectively.

The first equation in (\ref{eq.2.3.1}) is the equilibrium. The second equation in (\ref{eq.2.3.1}) is the compatibility. The third and forth equations in (\ref{eq.2.3.1}) are the Neumann and Dirichlet boundary conditions, respectively. To obtain the solutions to the boundary value problem in Eq. (\ref{eq.2.3.1}), material laws that describe the relation between stress and strain are required, e.g.,
\begin{equation}\label{eq.material_law}
\mathbf{S} = \mathbf{f} (\mathbf{E}), \quad \text{in} \hspace{0.1cm} \Omega^X.
\end{equation}
The material law is typically constructed with a pre-defined function $\mathbf{f}$ based on experimental observation, mechanics principles, and mathematical simplification with model parameters calibrated from limited material data \cite{ghaboussi1991knowledge,Sussman2009} or by computational homogenization approaches such as $\text{FE}^2$ \cite{Feyel1999,Feyel2003}, which inevitably introduce \textcolor{blue}{materials modeling} empiricism and errors \cite{Latorre2014}.
Moreover, the consistent tangent stiffness associated with the material law is often required in nonlinear computation \cite{Belytschko:2000tz}.

For complex material systems, phenomenological material models are difficult to construct.
The physics-constrained data-driven computing framework \cite{kirchdoerfer2016data, ibanez2018manifold, he2019physics} offers an alternative which directly utilizes material data and bypasses the
need of phenomenological model construction.

\subsection{Data-driven modeling of nonlinear elasticity}\label{sec.pcdd}
In this framework, the material behavior is described by means of
strain and stress tensors $(\hat{\mathbf{E}},\hat{\mathbf{S}})$ given by the material genome database. A material database $\mathbb{E}=\{(\hat{\mathbf{E}}_I,\hat{\mathbf{S}}_I)\}_{I=1}^M \in \mathcal{E}$ is defined to store the material data,
where $M$ is the number of material data points, and the hat symbol "\textsuperscript{$\wedge$}" is used to denote material data. Here, $\mathcal{E}$ denotes the admissible set of material database, which will be further discussed in Section \ref{sec.local_solver}.
To search for the most suitable (closest) strain-stress pairs $(\hat{\mathbf{E}}^*,\hat{\mathbf{S}}^*)$ for a given state $(\mathbf{E},\mathbf{S})$, an energy-like distance function extended from \cite{kirchdoerfer2016data} is defined:
\begin{equation}\label{eq.2.3.6}
\mathcal{F}(\mathbf{E},\mathbf{S};\hat{\mathbf{E}}^*,\hat{\mathbf{S}}^*) = 
    \underset{(\hat{\mathbf{E}},\hat{\mathbf{S}}) \in \mathcal{E}}{\mathrm{min}} \hspace{0.1cm}
     \int_{\Omega^X} \left( d_E^2(\mathbf{E},\hat{\mathbf{E}}) + d_S^2(\mathbf{S},\hat{\mathbf{S}}) \right) d\Omega,
\end{equation}
with
\begin{align}
    d_E^2(\mathbf{E},\hat{\mathbf{E}}) & = \frac{1}{2}(\mathbf{E}-\hat{\mathbf{E}}):\hat{\mathbb{C}}:(\mathbf{E}-\hat{\mathbf{E}}), \label{eq.2.3.7}\\
    d_S^2(\mathbf{S},\hat{\mathbf{S}}) & = \frac{1}{2}(\mathbf{S}-\hat{\mathbf{S}}):\hat{\mathbb{C}}^{-1}:(\mathbf{S}-\hat{\mathbf{S}}),\label{eq.2.3.8}
\end{align}
where $\hat{\mathbb{C}}$ is a predefined symmetric and positive-definite tensor used to properly regulate the distances between $(\mathbf{E}, \mathbf{S})$ and $(\hat{\mathbf{E}},\hat{\mathbf{S}})$. Usually, the selection of the coefficient matrix $\hat{\mathbb{C}}$ depends on the a priori knowledge of the given dataset. One widely adopted normalization scheme in machine learning is based on the variance of the data, but the selection is not unique. For example, the Mahalanobis distance of data is employed to compute the coefficient matrices \cite{ayensa2018new}. Recently, He et al. \cite{HeJBM2020} proposed to use the ratio of the standard deviations of the associated components of the stress–strain data to construct a diagonal coefficient matrix. Henceforth, the strain-stress pair $(\hat{\mathbf{E}},\hat{\mathbf{S}}) \in \mathcal{E}$ extracted from the material database is called the \textit{material data (state)}, whereas $(\mathbf{E},\mathbf{S})$ is called the \textit{physical state} if it satisfies the physically admissible set $\mathcal{C}$ given by the equilibrium and compatibility equations in Eq. (\ref{eq.2.3.1}), denoted as $(\mathbf{E},\mathbf{S}) \in \mathcal{C}$.

\textcolor{blue}{As a result, the data-driven modeling problem can be formulated as:
\begin{equation}\label{eq.dd_problem}
\underset{(\mathbf{E},\mathbf{S}) \in \mathcal{C}}{\mathrm{min}} \mathcal{F}(\mathbf{E},\mathbf{S};\hat{\mathbf{E}}^*,\hat{\mathbf{S}}^*) 
= \underset{(\mathbf{E},\mathbf{S}) \in \mathcal{C}}{\mathrm{min}} \underset{(\hat{\mathbf{E}},\hat{\mathbf{S}}) \in \mathcal{E}}{\mathrm{min}} \hspace{0.1cm}
     \int_{\Omega^X} \left( d_E^2(\mathbf{E},\hat{\mathbf{E}}) + d_S^2(\mathbf{S},\hat{\mathbf{S}}) \right) d\Omega.
\end{equation}
This data-driven problem is solved by fixed-point iterations, where the minimization of $\mathcal{F}$ with respect to $(\mathbf{E},\mathbf{S})$ and $(\hat{\mathbf{E}},\hat{\mathbf{S}})$ are performed iteratively until the intersection of two sets, $\mathcal{C}$ and $\mathcal{E}$, is found within a prescribed tolerance. We denote the minimization corresponding to the material data as the \textit{local step}, i.e. Eq. (\ref{eq.2.3.6}), while the one associated with the physical state as the \textit{global step}, which will be discussed as follows.}

Given the optimal material data $(\hat{\mathbf{E}}^*,\hat{\mathbf{S}}^*)$, the global step of the data-driven problem (\ref{eq.dd_problem}) is expressed as the following constrained minimization problem \cite{HeJBM2020}:
\begin{equation}\label{eq.2.3.11}
    \begin{aligned}
        \underset{\mathbf{u}, \mathbf{S}}{\mathrm{min}}
        \hspace{0.1cm} \mathcal{F}(\mathbf{E(\mathbf{u})},\mathbf{S};\hat{\mathbf{E}}^*,\hat{\mathbf{S}}^*) = 
        \underset{\mathbf{u}, \mathbf{S}}{\mathrm{min}} \hspace{0.1cm} & \int_{\Omega^X} \left(d_E^2(\mathbf{E}(\mathbf{u}),\hat{\mathbf{E}}^*) + d_S^2(\mathbf{S},\hat{\mathbf{S}}^*) \right) d\Omega \\
        \text{subject to:} \hspace{0.1cm}
        & DIV(\mathbf{F}(\mathbf{u}) \cdot \mathbf{S}) + \mathbf{b} = \mathbf{0} \hspace{0.2cm}\text{in}\hspace{0.2cm} \Omega^X, \\
        & (\mathbf{F}(\mathbf{u}) \cdot \mathbf{S}) \cdot \mathbf{N} = \mathbf{t} \hspace{0.2cm}\text{on}\hspace{0.2cm} \Gamma_t^X.
    \end{aligned}
\end{equation}
With the Lagrange multipliers $\boldsymbol{\lambda}$ and $\boldsymbol{\eta}$, Eq. (\ref{eq.2.3.11}) is transformed to the minimization of the following functional:
\begin{equation}\label{eq.Lagrange1a}
    \begin{split}
        & \mathcal{F}(\mathbf{E(\mathbf{u})},\mathbf{S};\hat{\mathbf{E}}^*,\hat{\mathbf{S}}^*) + \\
        & \int_{\Omega^X}
        \boldsymbol{\lambda} \cdot [DIV(\mathbf{F}(\mathbf{u}) \cdot \mathbf{S}) + \mathbf{b}] d\Omega +
        \int_{\Gamma_t^X}
        \boldsymbol{\eta} \cdot
        [(\mathbf{F}(\mathbf{u}) \cdot \mathbf{S}) \cdot \mathbf{N} - \mathbf{t}] d\Gamma.
    \end{split}
\end{equation}
The Euler-Lagrange equation of Eq. (\ref{eq.Lagrange1a}) indicates that $\boldsymbol{\eta} = - \boldsymbol{\lambda}$ on $\Gamma_t^X$ and $\boldsymbol{\lambda} = \boldsymbol{0}$ on $\Gamma_u^X$ \cite{Felippa1994}. Consequently, we have
\begin{equation}\label{eq.Lagrange1b}
    \begin{split}
        & \mathcal{F}(\mathbf{E(\mathbf{u})},\mathbf{S};\hat{\mathbf{E}}^*,\hat{\mathbf{S}}^*) + \\
        &
        \int_{\Omega^X} \boldsymbol{\lambda} \cdot
        [DIV(\mathbf{F}(\mathbf{u}) \cdot \mathbf{S}) + \mathbf{b}] d\Omega - 
        \int_{\Gamma_t^X} \boldsymbol{\lambda} \cdot [
        (\mathbf{F}(\mathbf{u}) \cdot \mathbf{S}) \cdot \mathbf{N} - \mathbf{t}] d\Gamma . 
    \end{split}
\end{equation}
By means of integration by parts and the divergence theorem, Eq. (\ref{eq.Lagrange1b}) is reformulated as
\textcolor{blue}{
\begin{equation}\label{eq.Lagrange2}
    \begin{split}
        & \mathcal{F}(\mathbf{E(\mathbf{u})},\mathbf{S};\hat{\mathbf{E}}^*,\hat{\mathbf{S}}^*) - \\
        & \int_{\Omega^X}
        [\nabla \boldsymbol{\lambda} : (\mathbf{F}(\mathbf{u}) \cdot \mathbf{S}) - \boldsymbol{\lambda} \cdot \mathbf{b}] d\Omega +
        \int_{\Gamma_t^X}
        \boldsymbol{\lambda} \cdot \mathbf{t} d\Gamma.
    \end{split}
\end{equation}
}
The stationary conditions of Eq. (\ref{eq.Lagrange2}) read:
\begin{subequations}\label{eq.variation}
    \begin{align}
    \delta \mathbf{u}: \hspace{0.2cm} & \int_{\Omega^X}
        \delta \mathbf{E}(\mathbf{u}) : \hat{\mathbb{C}}: (\mathbf{E}(\mathbf{u})-\hat{\mathbf{E}}^*) d\Omega  =
        \int_{\Omega^X}
        \delta \mathbf{F}^T(\mathbf{u}) \cdot \nabla \boldsymbol{\lambda} : \mathbf{S} d\Omega, \\
    \delta \mathbf{S}: \hspace{0.2cm} & \int_{\Omega^X}
        \delta \mathbf{S}:(\hat{\mathbb{C}}^{-1}:\mathbf{S} -
        \mathbf{F}^T(\mathbf{u}) \cdot \nabla \boldsymbol{\lambda}) d\Omega =
        \int_{\Omega^X}
        \delta \mathbf{S}:\hat{\mathbb{C}}^{-1}: \hat{\mathbf{S}}^* d\Omega, \\
    \delta \mathbf{\boldsymbol{\lambda}}: \hspace{0.2cm} &
        \int_{\Omega^X} \delta \nabla \boldsymbol{\lambda}:
        (\mathbf{F}(\mathbf{u}) \cdot \mathbf{S}) d\Omega =
        \int_{\Omega^X} \delta \boldsymbol{\lambda} \cdot
        \mathbf{b} d\Omega +  \int_{\Gamma_t^X}
        \delta \boldsymbol{\lambda} \cdot \mathbf{t}
        d\Gamma.
    \end{align}
\end{subequations}
As Eq. (\ref{eq.variation}b) provides correction between the physical stress $\mathbf{S}$ and the material stress data $\mathbf{S}^*$, a collocation approach is considered in Eq. (\ref{eq.variation}b) to yield:
\begin{equation}\label{eq.stress_updatre}
    \mathbf{S} =
        \hat{\mathbb{C}}:(\mathbf{F}^T(\mathbf{u}) \cdot \nabla \boldsymbol{\lambda}) + \hat{\mathbf{S}}^*,
\end{equation}
which represents the stress solution update. Substituting Eq. (\ref{eq.stress_updatre}) into Eqs. (\ref{eq.variation}a) and (\ref{eq.variation}c) yields:
\begin{subequations}\label{eq.global}
\begin{align}
 \begin{split}
    \int_{\Omega^X} \bigl[ \delta \mathbf{E}(\mathbf{u}) : \hat{\mathbb{C}}: (\mathbf{E}(\mathbf{u})-\hat{\mathbf{E}}^*) - & (\delta \mathbf{F}^T(\mathbf{u}) \cdot \nabla \boldsymbol{\lambda}) : \hat{\mathbb{C}}: (\mathbf{F}^T(\mathbf{u}) \cdot \nabla \boldsymbol{\lambda}) \bigl] d\Omega \\
    & = \int_{\Omega^X} (\delta \mathbf{F}^T(\mathbf{u}) \cdot \nabla \boldsymbol{\lambda}) : \hat{\mathbf{S}}^* d\Omega,
  \end{split}
  \\
  \begin{split}
    \int_{\Omega^X} (\mathbf{F}^T(\mathbf{u}) \cdot \delta \nabla     \boldsymbol{\lambda}):[ \hat{\mathbb{C}}:(\mathbf{F}^T(\mathbf{u})
    & \cdot \nabla \boldsymbol{\lambda}) + \hat{\mathbf{S}}^*] d\Omega \\
    & = \int_{\Omega^X} \delta \boldsymbol{\lambda} \cdot
        \mathbf{b} d\Omega +  \int_{\Gamma_t^X} \delta \boldsymbol{\lambda} \cdot \mathbf{t} d\Gamma.
  \end{split}
\end{align}
\end{subequations}
The solutions $\mathbf{u}$ and $\boldsymbol{\lambda}$ are solved from Eqs. (\ref{eq.global}) by means of the Newton-Raphson method \cite{Belytschko:2000tz}, and the physical state stress $\mathbf{S}$ is subsequently obtained from (\ref{eq.stress_updatre}).
As such, Eqs. (\ref{eq.stress_updatre})-(\ref{eq.global}) are the computational procedures to solve Eq. (\ref{eq.2.3.11}).
Moreover, in this boundary value problem, Eq. (\ref{eq.global}), the optimal material data $(\hat{\mathbf{E}}^*,\hat{\mathbf{S}}^*)$ not only provides the underlying material information learned from material database, but also serves as the material data-based connection to relate the strain in compatibility and the stress in equilibrium equations in Eqs. (\ref{eq.global}a) and (\ref{eq.global}b), respectively.

\textcolor{blue}{
In this study, the reproducing kernel particle method (RKPM) \cite{liu1995reproducing, chen1996reproducing} is employed to discretize the unknown fields $\boldsymbol{u}$ and $\boldsymbol{\lambda}$ in Eqs. (\ref{eq.global}) due to its capabilities of nodal approximation of state variables and enhanced smoothness that are particularly effective for data-driven computing. The formulation of the reproducing kernel approximation is given in \ref{appendix:RKPM}. Moreover, for effectiveness in data-driven computing, a stabilized nodal integration scheme (SCNI \cite{chen2002non}, see \ref{appendix:SCNI}) is used to integrate the weak formulations in Eq. (\ref{eq.global}) for reducing the number of integration points where the stress and strain material data need to be searched \cite{he2019physics}.
}

Combining the ease of introducing arbitrary order of continuity in the RK approximation as well as the employment of the SCNI scheme to integrate the weak equations in Eq. (\ref{eq.global}), it allows the material data search and variable evaluation to be performed only at the nodal points, avoiding the necessity of computing field and state variables separately at nodal points and Gauss points, respectively, for enhanced efficiency and accuracy of data-driven computing. \textcolor{blue}{In addition, under the SCNI framework, the nodal physical stress $\mathbf{S}$ is directly associated with the material stress data without introducing additional interpolation errors, see \cite{he2019physics,HeJBM2020} for more details.}
The Green Lagrangian strain tensor $\mathbf{E}$ is computed from the RK approximated displacement $\mathbf{u}$ and is evaluated at the nodal points. Note that stress update equation in Eq. (\ref{eq.stress_updatre}) is also carried out nodally.

The physical state $(\mathbf{E},\mathbf{S})$ evaluated at the integration points $\mathbf{x}_\alpha$ are denoted as $\{(\mathbf{E}_\alpha,\mathbf{S}_\alpha)\}_{\alpha = 1}^N$, where $N$ is the number of integration points. Note that due to the employment of nodal integration \cite{chen2002non,he2019physics}, the integration points share the same set of points as the nodal points. For simplicity of notation, we denote $\mathbf{z}_\alpha = (\mathbf{E}_\alpha,\mathbf{S}_\alpha)$ as the physics state and $\hat{\mathbf{z}}_\alpha = (\hat{\mathbf{E}}_\alpha,\hat{\mathbf{S}}_\alpha)$ the material data associated with the integration point $\alpha$ in the following discussions.

\textcolor{blue}{In data-driven computing, the local step (\ref{eq.2.3.6}) and the global step (\ref{eq.2.3.11}) are solved iteratively to search for the optimal material data from the admissible material set $\mathcal{E}$ that is closest to the physical state satisfying the physical constraints given in Eq. (\ref{eq.2.3.1}).
The convergence properties of this fixed-point iteration solver have been investigated in \cite{kirchdoerfer2016data,Conti2018}. It should be noted that the selection of the optimal material data by the material data-driven local solver is crucial to the effectiveness of data-driven computing \cite{kirchdoerfer2017data,he2019physics}, and further discussion will be presented in the following Sections \ref{sec.local_solver} and \ref{sec.AEDD}.}

\subsection{Material data-driven local solver} \label{sec.local_solver}
\begin{figure}[!ht]
\centering
    \begin{subfigure}{0.515\textwidth}
        \centering
        \includegraphics[width=1\linewidth]{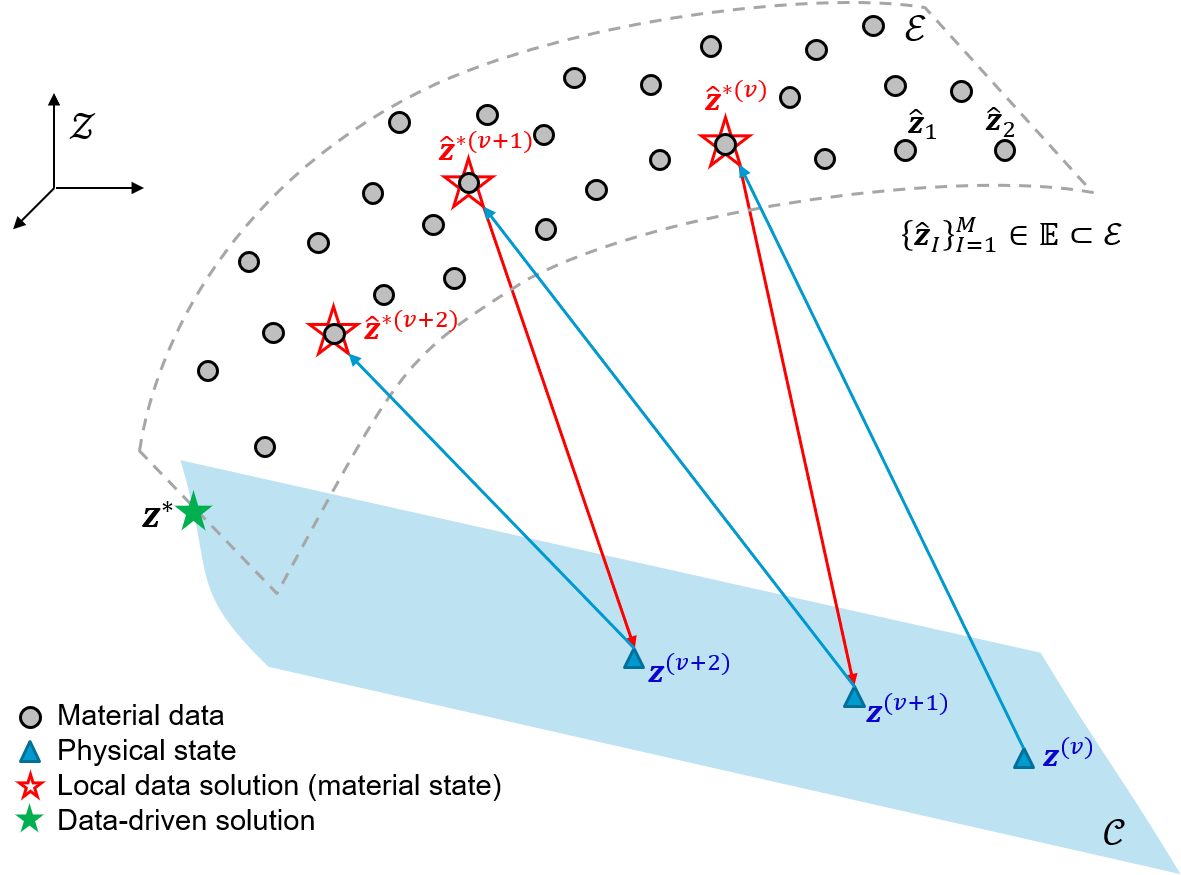}
        {\caption{DMDD}}
    \end{subfigure}
    \begin{subfigure}{0.475\textwidth}
        \centering
        \includegraphics[width=1\linewidth]{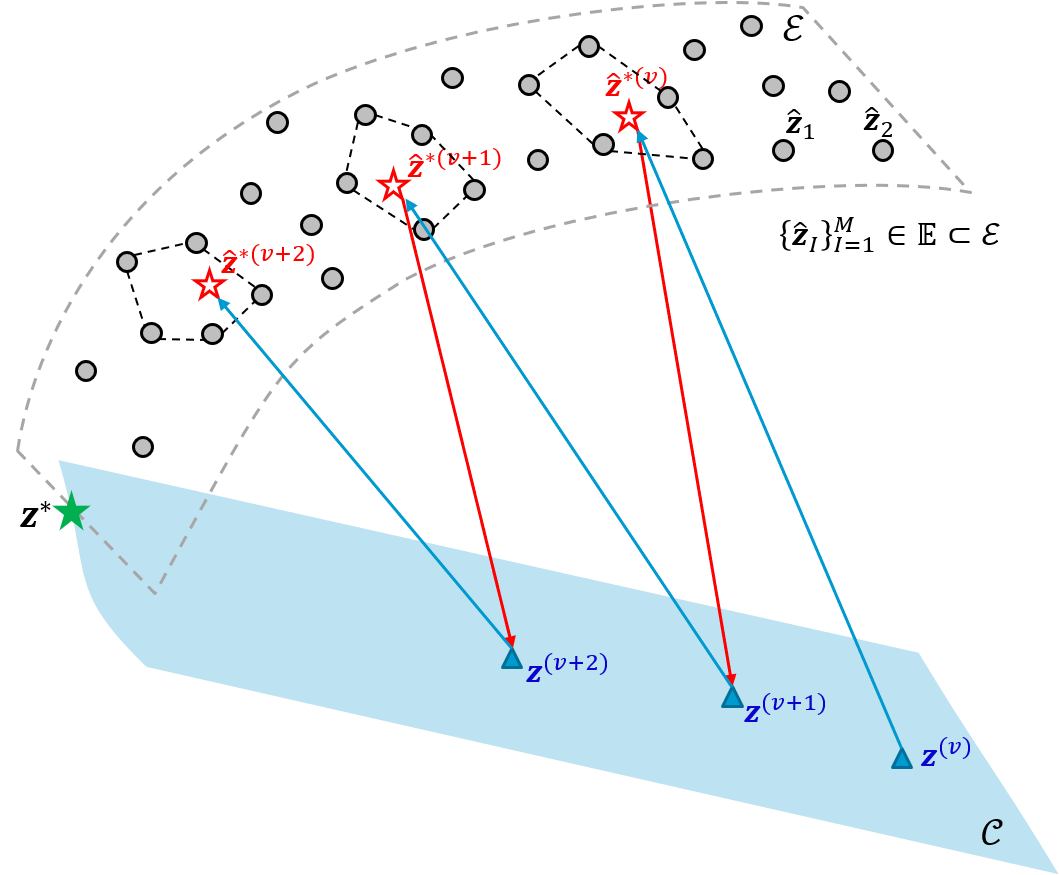}
        {\caption{LCDD}}
    \end{subfigure}
\caption{Geometric schematics of the (a) DMDD \cite{kirchdoerfer2016data} and (b) LCDD \cite{he2019physics} solvers. The data-driven solution $\mathbf{z}^*$ is given by the intersection of the admissible set of physical states $\mathcal{C}$ and the material admissible set $\mathcal{E}$.}\label{fig.DMDDvsLCDD}
\end{figure}

The effectiveness of data-driven computational paradigm discussed in Section \ref{sec.pcdd} relies heavily on the search of optimal material data $\hat{\mathbf{z}}_{\alpha}^* = (\hat{\mathbf{E}}_\alpha^*,\hat{\mathbf{S}}_\alpha^*)$, $\alpha=1,...,N$, from the material dataset $\mathbb{E}$. When adopting the distance-minimizing data-driven (DMDD) approach proposed in \cite{kirchdoerfer2016data}, the local material solver in Eq. (\ref{eq.2.3.6}) used to find the optimal material data can be defined as
\begin{equation}\label{eq.dmdd_loc1}
    (\hat{\mathbf{E}}_{\alpha}^*,\hat{\mathbf{S}}_{\alpha}^*) = \underset{(\hat{\mathbf{E}}_{\alpha},\hat{\mathbf{S}}_{\alpha}) \in \mathbb{E}}{\mathrm{arg \hspace{0.05cm} min}} \hspace{0.1cm} d_E^2(\mathbf{E}_{\alpha},\hat{\mathbf{E}}_{\alpha}) + d_S^2(\mathbf{S}_{\alpha},\hat{\mathbf{S}}_{\alpha}), \hspace{0.5cm} \alpha=1,...,N,
\end{equation}
where the distance functions $d_E$ and $d_S$ are referred to Eqs. (\ref{eq.2.3.7}) and (\ref{eq.2.3.8}), $\alpha$ denotes the indices of integration points, and $N$ is the total number of integration points.

For enhanced data-driven computing, especially with sparse noisy data, an alternative approach called local convexity data-driven (LCDD) computing was proposed in \cite{he2019physics} by introducing the underlying structure of material data via manifold learning. In this approach, the local solver is expressed as:
\begin{equation}\label{eq.lcdd_loc}
    (\hat{\mathbf{E}}_{\alpha}^*,\hat{\mathbf{S}}_{\alpha}^*) = \underset{(\hat{\mathbf{E}}_{\alpha},\hat{\mathbf{S}}_{\alpha}) \in \mathcal{E}_{\alpha}^{l}}{\mathrm{arg \hspace{0.05cm} min}} \hspace{0.1cm} d_E^2(\mathbf{E}_{\alpha},\hat{\mathbf{E}}_{\alpha}) + d_S^2(\mathbf{S}_{\alpha},\hat{\mathbf{S}}_{\alpha}), \hspace{0.5cm} \alpha=1,...,N,
\end{equation}
where $\mathcal{E}_{\alpha}^{l}:=\mathcal{E}^{l}(\mathbf{z}_\alpha)$ denotes a local convex subset formed by $k$ material data points closest to the given physical state $\mathbf{z}_\alpha = (\mathbf{E}_{\alpha},\mathbf{S}_{\alpha})$, as illustrated by the polygons in Fig. \ref{fig.DMDDvsLCDD}(b). The local minimization problem (\ref{eq.lcdd_loc}) is solved by means of a non-negative least-square algorithm with penalty relaxation, see details in \cite{he2019physics}. This approach introduces a local embedding reconstruction of datsa which is more robust in dealing noisy data and outliers.

Fig. \ref{fig.DMDDvsLCDD} shows the comparison of the DMDD and LCDD solvers, where $(v)$ is the iteration index and one iteration consists of solving one global (physical) step, i.e. Eqs. (\ref{eq.stress_updatre})--(\ref{eq.global}) and one material data-driven local step, e.g. Eq. (\ref{eq.dmdd_loc1}) or (\ref{eq.lcdd_loc}), as noted in Section \ref{sec.pcdd}.
The local step of the DMDD solver (\ref{eq.dmdd_loc1}) searches for the material data closest to the given physical state directly from the material dataset $\mathbb{E}$, see Fig. \ref{fig.DMDDvsLCDD}(a).
It has been shown that this heuristic solver suffers from noisy dataset and requires enormous data to guarantee satisfactory accuracy \cite{kirchdoerfer2017data,he2019physics}. 
On the other hand, the LCDD solver (\ref{eq.lcdd_loc}) searches for the optimal material data based on the locally constructed convex space $\mathcal{E}_{\alpha}^{l}$ informed by the neighboring data, as shown in Fig. \ref{fig.DMDDvsLCDD}(b).
The key idea behind the construction of $\mathcal{E}_{\alpha}^{l}$ is to provide a smooth, bounded and lower dimensional admissible space for optimal material data search in Eq. (\ref{eq.lcdd_loc}), and to preserve the convexity of the constructed local material manifold for enhanced robustness and stability in data-driven iterations.

While the material data-driven local solver in Eq. (\ref{eq.lcdd_loc}) locates the optimal data from the defined feasible set (constructed by a set of local neighboring points), the final solution of the associated data-driven modeling problem Eq. (\ref{eq.dd_problem}) is not
guaranteed to be globally optimal.
This is consistent to other existing data-driven approaches \cite{kirchdoerfer2016data,Conti2018,he2019physics} where the optimality fundamentally depends on the characteristic of the material dataset.
However, as demonstrated in references \cite{kirchdoerfer2016data,Conti2018,he2019physics}, if the material dataset is well posed, the proposed data-driven solver can converge optimally as the density of data points increases.

\begin{remark}
It should be noticed that in both Eqs. (\ref{eq.dmdd_loc1}) and (\ref{eq.lcdd_loc}) the nearest points are sought based on the metric functions $d_E$ and $d_S$.
Thus, it suffers from the notorious “dimensionality curse” when data-driven modeling attempts to scale up to high-dimensional material data.
Although the innate manifold learning in LCDD allows noise and dimensionality reduction, the proper definition of the metric functions in high-dimensional phase space remains challenging \cite{goodfellow2016deep}.
Besides, as the nearest neighbors are searched locally from the existing data points of the material dataset, it leads to limited extrapolative generalization to be demonstrated in Section \ref{sec3.2}.
Furthermore, the data search and the locally convex reconstruction through a constrained minimization solver at every local step during data-driven computation could result in high computational cost especially for the large and high dimensional material dataset.
\end{remark}

To address the issue of the curse of dimensionality, we propose to use autoencoders in the data-driven local solver for deep manifold learning of material data, allowing effective discovering of the underlying representation of stress-strain material data. To the best of the authors' knowledge, this is the first attempt to apply deep manifold learning in physics-constrained data-driven computing. In the following exposition, we demonstrate how autoencoder based deep learning enhances accuracy, robustness, and generalization ability of data-driven computing.

%

\section{Autoencoders \textcolor{blue}{for low-dimensional nonlinear representation of material data}}\label{sec.AEDD_local}
\textcolor{blue}{For effective data search in the local step, Eq. (\ref{eq.2.3.6}) or Eqs. (\ref{eq.dmdd_loc1}) and (\ref{eq.lcdd_loc}), in} this section, we first review the basic concepts of deep neural networks and autoencoders that are used for deep manifold learning. \textcolor{blue}{We then} present \textcolor{blue}{the employment of} autoencoders to construct low-dimensional nonlinear representation (embedding) of material data. \textcolor{blue}{Thereafter, optimum data search on the low-dimensional data manifold using a locally convex projection method is presented in Section \ref{sec.AEDD}.}

\subsection{Background: Autoencoders}\label{sec.autoencoder}

As the core of the deep learning \cite{goodfellow2016deep}, deep neural networks (DNNs) or often called multilayer perceptrons (MLPs), are used to represent a complex model relating data inputs, $\mathbf{x} \in \mathbb{R}^{d_{in}}$ and data outputs $\mathbf{y} \in \mathbb{R}^{d_{out}}$. A typical DNN is composed of an input layer, an output layer, and $L$ hidden layers. Each hidden layer transforms the outputs of the previous layer through two operators, i.e., an affine mapping followed by a nonlinear activation function $\sigma(\cdot)$, and outputs the results to the next layer, which can be written as:
\begin{equation}\label{eq.2.1.1}
    \mathbf{x}^{(l)} = \sigma(\mathbf{W}^{(l)} \mathbf{x}^{(l-1)} + \mathbf{b}^{(l)}), \hspace{0.5cm} l = 1,...,L,
\end{equation}
where $\mathbf{x}^{(l)} \in \mathbb{R}^{n_l}$ is the outputs of layer $l$ with $n_l$ neurons,
and $\mathbf{W}^{(l)} \in \mathbb{R}^{n_l \times n_{l-1}}$ and $\mathbf{b}^{(l)} \in \mathbb{R}^{n_l}$ are the weight matrix for linear mapping and the bias vector of layer $l$, respectively, where $n_0 = d_{in}$ is the input dimension.
Some of the commonly used activation functions include logistic sigmoid, rectified linear unit (ReLu), and leaky ReLu. In this study, a hyperbolic tangent function is used as the activation function for hidden layers, $\sigma(\cdot) = tanh(\cdot)$.
Note that the setup of the output layer depends on the type of machine learning tasks, e.g., classification, regression. For regression tasks, which is the application of this study, a linear function is used in the output layer where the last hidden layer information is mapped to the output vector $\Tilde{\mathbf{y}}$, expressed as: $\Tilde{\mathbf{y}} = \mathbf{W}^{(L+1)} \mathbf{x}^{(L)} + \mathbf{b}^{(L+1)}$, where $\Tilde{\mathbf{y}}$ denotes the DNN approximation of the output $\mathbf{y}$.
We denote $\boldsymbol{\theta}$ as the collection of all trainable weight and bias coefficients, $\boldsymbol{\theta} = \{\mathbf{W}^{(l)}, \mathbf{b}^{(l)}\}_{l=1}^{L+1}$.

\begin{figure}[ht!]
    \centering
    \includegraphics[width=0.7\textwidth]{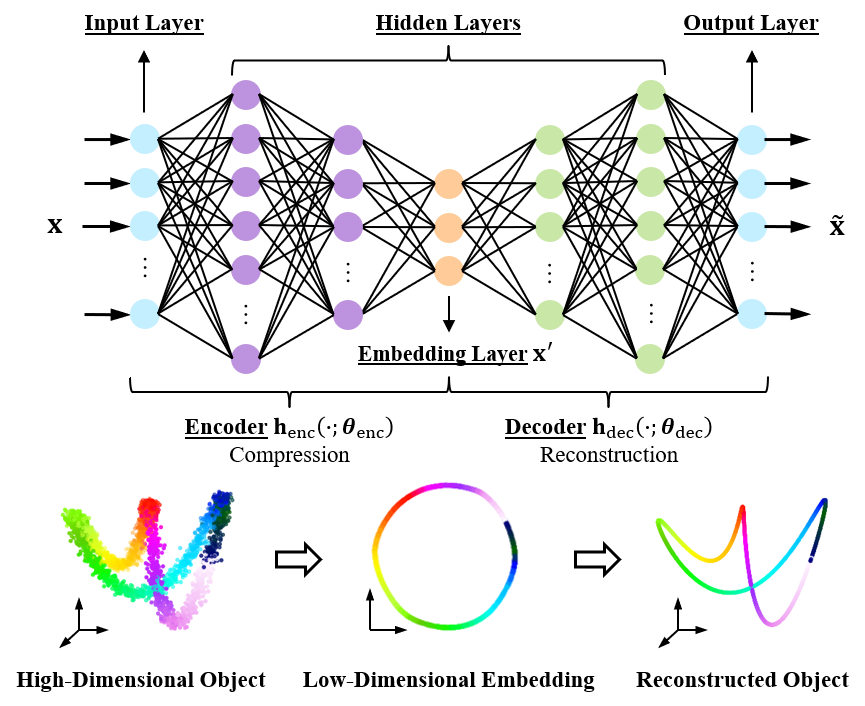}
    \caption{Schematic of an autoencoder consisting of an encoder and a decoder, where the dimension of the embedding layer is smaller than the input dimension. For a high-dimensional input object, the encoder learns a compressed low-dimensional embedding, on which the decoder optimally reconstructs the input object.}
    \label{fig2.1.autoencoder}
\end{figure}

Autoencoders \cite{demers1993non,hinton2006reducing} are an unsupervised learning technique in which special architectures of DNNs are leveraged for dimensionality reduction or representation learning. Specially, an autoencoder aims to optimally copy its input to output with the most representative features by introducing a low-dimensional embedding layer (or called \textbf{a code}).
As shown in Fig. \ref{fig2.1.autoencoder}, an autoencoder consists of two parts, an encoder function $\mathbf{h}_{\text{enc}}(\cdot;\boldsymbol{\theta}_{\text{enc}}) : \mathbb{R}^{d} \rightarrow \mathbb{R}^{p}$ and a decoder function $\mathbf{h}_{\text{dec}}(\cdot;\boldsymbol{\theta}_{\text{dec}}): \mathbb{R}^{p} \rightarrow \mathbb{R}^{d}$, such that the autoencoder is
\begin{subequations} \label{eq.auto}
    \begin{align}
        \Tilde{\mathbf{x}} = \mathbf{h} (\mathbf{x}; \boldsymbol{\theta}_{\text{enc}}, \boldsymbol{\theta}_{\text{dec}}) 
        & := (\mathbf{h}_{\text{dec}} \circ \mathbf{h}_{\text{enc}})(\mathbf{x}) \\
        & := \mathbf{h}_{\text{dec}}(\mathbf{h}_{\text{enc}}(\mathbf{x};\boldsymbol{\theta}_{\text{enc}});\boldsymbol{\theta}_{\text{dec}}),
    \end{align}
\end{subequations}
where $p < d$ is the embedding dimension, $\boldsymbol{\theta}_{\text{enc}}$ and $\boldsymbol{\theta}_{\text{dec}}$ are the DNN coefficients of encoder and deconder parts, respectively, and $\Tilde{\mathbf{x}}$ is the output of the autoencoder, a reconstruction of the original input $\mathbf{x}$.
With the latent dimension $p$ much less than the input dimension $d$, the encoder $\mathbf{h}_{\text{enc}}$ is trained to learn the compressed representation of  $\mathbf{x}$, denoted as the embedding $\mathbf{x}' \in \mathbb{R}^p$,
whereas the decoder $\mathbf{h}_{\text{dec}}$ reconstructs the input data by mapping the embedding representation back to the high-dimensional space.


It is important to note that similar to any other dimensionality reduction techniques~\cite{Lee2007}, the employment of autoencoders is based on the \textit{manifold hypothesis}, which presumes that the given high-dimensional input data, e.g., the material dataset $\mathbb{E}$, lies on a low-dimensional manifold $\mathcal{E}'$ that is embedded in a higher-dimensional vector space, as shown by the schematic figures at the bottom of Fig. \ref{fig2.1.autoencoder}.


\subsection{Nonlinear material embedding}\label{sec.auto_mat}
%
In this study, autoencoders are used to discover the intrinsic low-dimensional material embedding of the given material dataset $\mathbb{E}=\{\hat{\mathbf{z}}_I\}_{I=1}^M$, where $\hat{\mathbf{z}}_I=(\hat{\mathbf{E}}_I,\hat{\mathbf{S}}_I)$ and $M$ is the number of material data points.
Given the autoencoder architecture $\mathbf{h} (\cdot; \boldsymbol{\theta}_{\text{enc}}, \boldsymbol{\theta}_{\text{dec}})$ in Eq. (\ref{eq.auto}), the parameters $\boldsymbol{\theta}_{\text{enc}}^*$ and $\boldsymbol{\theta}_{\text{dec}}^*$ are computed by minimizing the following loss function:
\begin{equation}\label{eq.loss_auto}
    (\boldsymbol{\theta}_{\text{enc}}^*, \boldsymbol{\theta}_{\text{dec}}^*) = \underset{\boldsymbol{\theta}_{\text{enc}}, \boldsymbol{\theta}_{\text{dec}}} \argmin \frac{1}{M} \sum^M_{I=1} ||\mathbf{h} (\hat{\mathbf{z}}_I; \boldsymbol{\theta}_{\text{enc}}, \boldsymbol{\theta}_{\text{dec}}) - \hat{\mathbf{z}}_I||^2 + \beta \sum^{L+1}_{l=1} ||\mathbf{W}^{(l)}||_F^2,
\end{equation}
where $\beta$ is a regularization parameter, and $||\cdot||_F$ denotes the Frobenius norm.
Here, the loss function consists of the reconstruction error over all training data and a $L_2$-norm based weight regularization term used to prevent over-fitting issues \cite{goodfellow2016deep,mishne2019diffusion}.

The training procedures of autoencoders in terms of the loss function in Eq. (\ref{eq.loss_auto}) are performed offline. Thus, training on a large material dataset does not result in additional overhead on the online data-driven computation.
The details of the training algorithms for autoencoders are given in Section \ref{sec.auto_train}.

\begin{remark}
It is well known that for an autoencoder with a single hidden layer and linear activation function, the weights trained by the mean-squared-error cost function learn to span the same principal subspace as principal components analysis (PCA)~\cite{Jolliffe2002}. 
Autoencoders based on neural networks with nonlinear transform functions can be thought of as a generalizaiton of PCA, capable of learning nonlinear relationships.
\end{remark}

Given the trained autoencoder $\mathbf{h} (\cdot; \boldsymbol{\theta}_{\text{enc}}^*, \boldsymbol{\theta}_{\text{dec}}^*)$, we can define a low-dimensional embedding space, $\mathcal{E}'=\{ \mathbf{z}' \in \mathbb{R}^p \:| \: \mathbf{z}'= \mathbf{h}_{\text{enc}}(\mathbf{z};\boldsymbol{\theta}_{\text{enc}}^*), \forall \mathbf{z} \in \mathcal{Z} \}$, in which the material state is described by a lower-dimensional coordinate system $\mathbf{z}'$.
Here, the prime symbol $(\cdot)'$ is used to denote the quantities defined in the embedding space, and $\mathcal{Z}$ denotes the high-dimensional phase space where the material states $\hat{\mathbf{z}}$ and the physical states $\mathbf{z}$ are defined. 
For example, the embedding set of the given material data is
\begin{equation}\label{eq.mat_emb}
    \quad \mathbb{E}^{'} = \{\hat{\mathbf{z}}'_I\}_{I=1}^M \subset \mathcal{E}^{'},
\end{equation}
where $\hat{\mathbf{z}}'_I = \mathbf{h}_{\text{enc}}(\hat{\mathbf{z}}_I;\boldsymbol{\theta}_{\text{enc}}^*)$ for $\hat{\mathbf{z}}_I \in \mathbb{E}$.


Considering the data-driven application on learning the underlying structure of material data, autoencoders provide the following advantages:
\begin{enumerate}
    \item [1)] Deep neural network architecture enables autoencoders to capture highly complex nonlinear manifold with exponentially less data points than nonparametric methods based on nearest neighbor graph \cite{bengio2009learning,van2009dimensionality,goodfellow2016deep}.
    \item [2)] Autoencoders provide explicit mapping functions, i.e. $\mathbf{h}_{\text{enc}}$ and $\mathbf{h}_{\text{dec}}$, between the high- and low-dimensional representation so that the trained encoders allow efficient evaluation of the embedding of new input data.
	\item [3)] Through information compression by encoders, unwanted information of material data, such as noise and outliers, can be filtered while preserving its essential low-dimensional manifold structure \cite{goodfellow2016deep}.
\end{enumerate}
Compared to data-driven methods based on conventional manifold learning techniques~\cite{ibanez2018manifold,he2019physics,Eggersmann}, the explicit nonlinear mapping functions learned by autoencoders are particularly attractive to data-driven computing because not only can they encode the essential global structure of the given material data for enhanced generalization ability, they also greatly  reduce online computational cost by using the pretrained autoencoders.
Furthermore, as we can see in next section, due to the availability of low-dimensional embedding $\mathcal{E}'$, we can introduce a convexity-preserving interpolation scheme to effectively search for the optimal material data associated with the given physical state.


\subsection{Autoencoder architectures and training algorithms}\label{sec.auto_train}
As the architectures of encoder and decoder in an autoencoder are symmetric, we only use the encoder architecture to denote the autoencoder architecture. For example, the encoder architecture in Fig. \ref{fig2.1.autoencoder} is $4-6-4-3$, where the first and last values denote the numbers of artificial neurons in the input layer and embedding layer, respectively, and the other values denote the neuron numbers of the hidden layers in sequence. As such, the decoder architecture in this case is $3-4-6-4$.

The offline training on the given material datasets is performed by using the open-source Pytorch library \cite{paszke2017automatic}, and the optimal parameters $\boldsymbol{\theta}_{\text{enc}}^*$ and $\boldsymbol{\theta}_{\text{dec}}^*$ of autoencoders are obtained by minimizing the loss function (Eq. (\ref{eq.loss_auto})).
The regularization parameter $\beta$ is set as $10^{-5}$. A hyperbolic tangent function is adopted as the activation function for all layers of autoencoders, except for the embedding layer and the output layer, where a linear function is employed instead.
To eliminate the need of manually tuning the learning rate for training, an adaptive gradient algorithm, Adagrad \cite{duchi2011adaptive}, is employed,
where the initial learning rate is set to be $0.1$ and the number of training epochs is set to be 2000. The training datasets are standardized such that they have zero mean and unit variance to accelerate the training process.
It should be noted that the training of autoencoders could get trapped in local minima and this can be overcome by pretraining the network using Restricted Boltzmann Machines or by denoising autoencoders \cite{vincent2008extracting,van2009dimensionality,larochelle2009exploring}.

\section{Auto-embedding data-driven (AEDD) solver}\label{sec.AEDD}
We now develop the AEDD solver based on autoencoders to search for the optimal material data in the solution process of the local step (e.g., Eqs. (\ref{eq.dmdd_loc1}) or (\ref{eq.lcdd_loc})). We begin with introducing a simple interpolation scheme to preserve local convexity in the material data search, which is essential in enhancing the local solver performance, followed by presenting two AEDD solvers with the employment of convexity-preserving reconstruction.

\subsection{Convexity-preserving interpolation}\label{sec.interpolation}
With deep manifold learning by autoencoders, we are able to extract the underlying low-dimensional global manifold of material datasets, and to enhance the generalization capability of the material local solver.
During data-driven computing, material and physical states are projected onto the constructed material embedding space $\mathcal{E}'$, and a convexity-preserving local data reconstruction is introduced for enhanced stability and convergence in the local data search on the embedding space.

In this approach, because the material embedding points in low-dimensional space are explicitly given by the offline trained autoencoders, interpolation schemes on the embedding space is straightforward without suffering the high-dimensionality issues.
A convexity-preserving, partition-of-unity interpolation method is therefore introduced into the material data-driven solver, using Shepard function \cite{shepard1968two} or inverse distance weighting.
Shepard interpolation has been widely used in data fitting and function approximation with positivity constraint \cite{babuvska1997partition,wendland2004scattered,he2014topology}.

Here, the Shepard functions are applied to reconstruct the material embedding of a given physical state, $\mathbf{z}' = \mathbf{h}_{\text{enc}}(\mathbf{z};\boldsymbol{\theta}_{\text{enc}}^*)$, by its material embedding neighbors, expressed as
\begin{equation}\label{eq.shep_interp}
    \mathbf{z}'_{recon} = \mathcal{I} \left( \{\Psi_I(\mathbf{z}');\hat{\mathbf{z}}'_I\}_{I \in \mathcal{N}_k(\mathbf{z}')} \right) = \sum_{I \in \mathcal{N}_k(\mathbf{z}')} \Psi_I(\mathbf{z}') \hat{\mathbf{z}}'_I,
\end{equation}
where $\mathbf{z}'_{recon}$ is the reconstruction of $\mathbf{z}'$, $\hat{\mathbf{z}}'_I$ is the material data embedding in $\mathbb{E}'$ defined in Eq. (\ref{eq.mat_emb}), $\mathcal{N}_k(\mathbf{z}')$ is the index set of the $k$ nearest neighbor points of $\mathbf{z}'$ selected from $\mathbb{E}'$,
and the shape functions are 
\begin{equation}\label{eq.shep_shap}
    \Psi_I(\mathbf{z}') = \frac{\phi(\mathbf{z}'-\hat{\mathbf{z}}'_I)}{\sum_{J=1} \phi(\mathbf{z}'-\hat{\mathbf{z}}'_J)}.
\end{equation}
In Eqs. (\ref{eq.shep_interp}) and (\ref{eq.shep_shap}), $\phi$ is a positive kernel function representing the weight on the data set $\{\hat{\mathbf{z}}'_I\}_{I \in \mathcal{N}_k(\mathbf{z}')}$, and
$\mathcal{I}$ denotes the interpolation operator that constructs shape functions with respect to $\mathbf{z}'$ and its neighbors.
Note that these functions form a partition of unity, i.e., $\sum_{I \in \mathcal{N}_k(\mathbf{z}')} \Psi_I(\mathbf{z}')=1$ for transformation objectivity. Furthermore, they are convexity-preserving when the kernel function $\phi$ is a positive function.
Here, an inverse distance function is used as the kernel function
\begin{equation}\label{eq.shep_kern}
    \phi(\mathbf{z}'-\hat{\mathbf{z}}'_I) = \frac{1}{|| \mathbf{z}'- \hat{\mathbf{z}}'_I ||^2}.
\end{equation}
It is worth noting that the interpolation functions defined in Eqs. (\ref{eq.shep_shap}) and (\ref{eq.shep_kern}) are equivalent to the RK approximation function in (\ref{eq.2.2.2}) with zero-order basis.

\begin{figure}[!ht]
    \centering
    \includegraphics[width=0.5\linewidth]{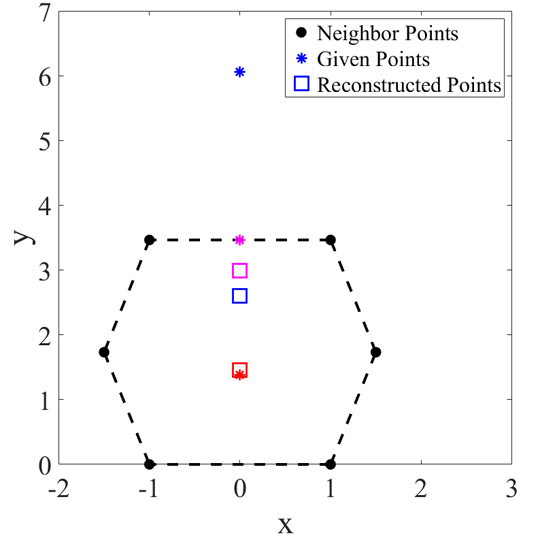}
    {\caption{Demonstration of the convexity-preserving reconstruction by Shepard interpolation in Eq. (\ref{eq.shep_interp}), where the asterisk and square denote the given and the reconstructed points, respectively, black dots represent the nearest neighbors of the given points, and the black dash line depicts a locally convex hull formed by the nearest neighbors.}\label{fig.shep}}
\end{figure}

Fig. \ref{fig.shep} demonstrates the locally convex reconstruction by the proposed interpolation in Eq. (\ref{eq.shep_interp}). 
For example, the given blue asterisk is mapped to the blue-square point by using the Shepard interpolation.
It can be seen that the three given points (inside (red), on-edge (pink), and outside (blue)) are all mapped to locations within the convex hull, showing the desired convexity-preserving capability.
The interpolation is simple and efficient as the interpolation functions in Eq. (\ref{eq.shep_shap}) can be constructed easily in a low-dimensional embedding space.

\subsection{Auto-embedding data-driven (AEDD) solver in data-driven computing}\label{AEDDsolvers}
The physics-constrained data-driven computing described in Section \ref{sec2} is conducted in the high-dimensional phase space $\mathcal{Z}$ (or called \textit{data space}), where the physical state $\mathbf{z}_\alpha \in \mathcal{C}$, the material data $\hat{\mathbf{z}}_\alpha \in \mathcal{E}$ and the material dataset $\mathbb{E}$ are defined in $\mathcal{Z}$.
We use the subscript "$\alpha$" to denote the quantities at integration points with the employment of numerical discretization, see Section \ref{sec.pcdd}.
To enhance solution accuracy and generalization capability of data-driven computing, deep manifold learning enabled by autoencoders is introduced into the material data-driven local solver.

Recall that autoencoders introduced in Section \ref{sec.auto_mat} are trained offline and the trained encoder $\mathbf{h}_{\text{enc}}$ and decoder $\mathbf{h}_{\text{dec}}$ functions are employed directly in the online data-driven computation.
As such, the encoder maps an arbitrary point from the data space to the embedding space, i.e. $\mathbf{z}'_\alpha = \mathbf{h}_{\text{enc}}(\mathbf{z}_\alpha)$, whereas the decoder performs the reverse mapping, i.e. $\Tilde{\mathbf{z}}_\alpha = \mathbf{h}_{\text{dec}}(\mathbf{z}'_\alpha)$. With the autoencoders and the proposed convexity-preserving data reconstruction in the embedding space introduced in Section \ref{sec.interpolation}, we propose the following two AEDD approaches for the material data-driven local solver. The objective is to find the optimal material data for a given physical state $\mathbf{z}_\alpha$ computed in Section \ref{sec.pcdd}.

\subsubsection{AEDD local solver: Solver I}\label{sec.AEDD_locI}
Let $\mathbf{z}_\alpha$ be the physical state obtained from the global step with physical constraints, Eq. (\ref{eq.2.3.11}), or the corresponding variational equations, Eqs. (\ref{eq.stress_updatre})--(\ref{eq.global}). We first introduce a local solver that uses decoders for reverse mapping from the embedding space to the data space, denoted as Solver I. In this approach, the local problem defined in Eq. (\ref{eq.2.3.6}) is reformulated by three steps, as described below:
\begin{subequations}\label{eq.AEDD_1}
    \begin{align}
        & \textit{Step 1}: \hspace{1cm} \mathbf{z}'_{\alpha} = \mathbf{h}_{\text{enc}} (\mathbf{z}_{\alpha}), \label{eq.2.3.13}\\
        & \textit{Step 2}: \hspace{1cm} \hat{\mathbf{z}}'^{*}_{\alpha} = \mathcal{I} \left( \{\Psi_I(\mathbf{z}'_{\alpha});\hat{\mathbf{z}}'_I\}_{I \in \mathcal{N}_k(\mathbf{z}'_{\alpha})} \right) \label{eq.2.3.14}\\
        & \textit{Step 3}: \hspace{1cm} \hat{\mathbf{z}}^{*}_{\alpha} = \mathbf{h}_{\text{dec}}(\hat{\mathbf{z}}'^{*}_{\alpha}), \label{eq.2.3.15}
    \end{align}
\end{subequations}
for $\alpha = 1,...,N$, where $\hat{\mathbf{z}}'_I \in \mathbb{E}^{'}$ (see Eq. (\ref{eq.mat_emb})), and $\mathcal{I}$ is the convexity-preserving interpolation operator defined in Eq. (\ref{eq.shep_interp}).

\begin{figure}[!ht]
\centering
    \begin{subfigure}{0.8\textwidth}
        \centering
        \includegraphics[width=1\linewidth]{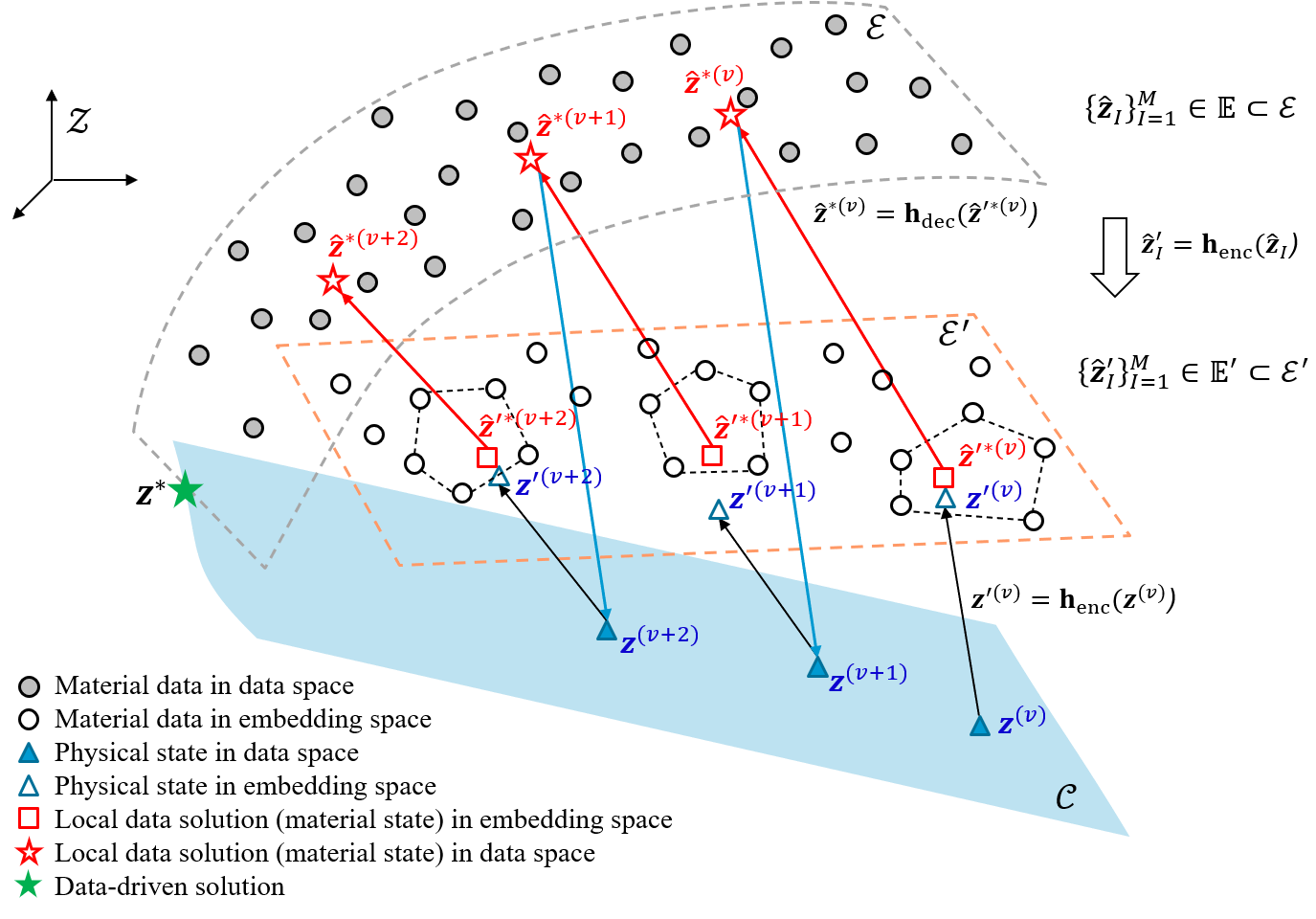}
        \caption{}
    \end{subfigure}
    \begin{subfigure}{0.8\textwidth}
        \centering
        \includegraphics[width=0.95\linewidth]{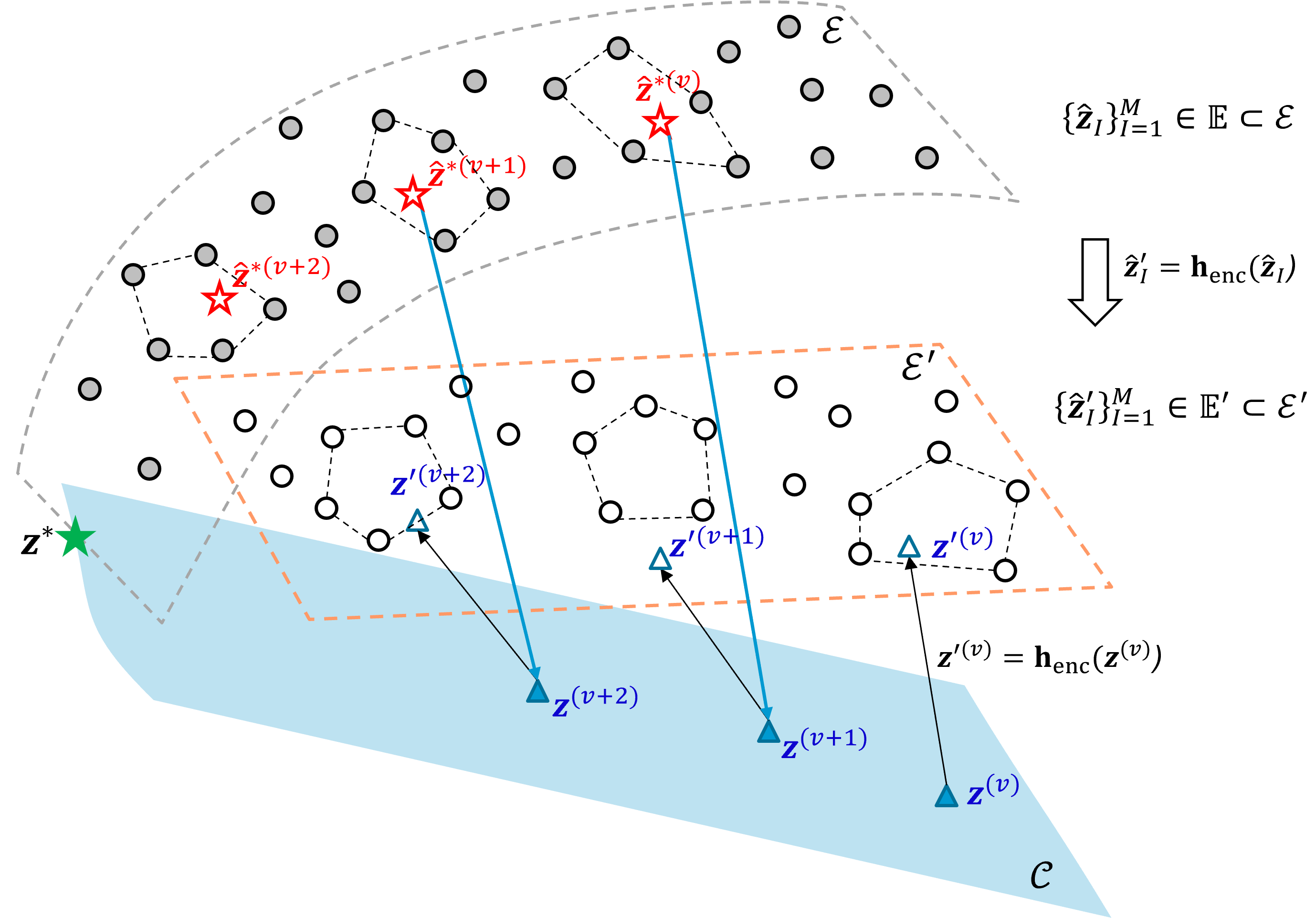}
        \caption{}
    \end{subfigure}
\caption{Geometric schematic of the proposed auto-embedding data-driven computational framework: (a) Solver I; (b) Solver II, corresponding to two different ways to reconstruct the optimal material data solution in high-dimensional data space. The material data points (the gray-filled circles), $\hat{\mathbf{z}}_I$, in the phase space are related to the material embedding points (the white-filled circles) $\hat{\mathbf{z}}_I'$ via the encoder function. The low-dimensional embedding manifold is represented by the orange dash line.}\label{fig2.3.1}
\end{figure}

The schematic of data-driven computing with Solver I is illustrated in Fig. \ref{fig2.3.1}(a), where the integration point index $\alpha$ is dropped for brevity.
For example, at the $v$-th global-local iteration, after the physical state $\mathbf{z}^{(v)}$ (the blue-filled triangle) is obtained from the global physical step (Eq. (\ref{eq.2.3.11})), \textit{Step 1} of the local solver (Eq. (\ref{eq.2.3.13})) maps the sought physical state from the data space to the embedding space by the encoder,  $\mathbf{z}'^{(v)} = \mathbf{h}_{\text{enc}}(\mathbf{z}^{(v)})$, depicted by the white-filled triangle in Fig. \ref{fig2.3.1}(a).
In \textit{Step 2}, $k$ nearest neighbors of $\mathbf{z}'^{(v)}$ based on Euclidean distance are sought in the embedding space and the optimal material embedding solution $\hat{\mathbf{z}}'^{*(v)}$ (the red square) is reconstructed by using the proposed convexity-preserving interpolation (Eqs. (\ref{eq.shep_interp})-(\ref{eq.shep_kern})). Lastly, in \textit{Step 3}, the optimal material embedding state $\hat{\mathbf{z}}'^{*(v)}$ is transformed from the embedding space to the data space by the decoder, $\hat{\mathbf{z}}^{*(v)} = \mathbf{h}_{\text{dec}}(\hat{\mathbf{z}}'^{*(v)})$ (the red star in Fig. \ref{fig2.3.1}(a)). Subsequently, this resultant material data solution $\hat{\mathbf{z}}^{*(v)}$ from the local solver in Eq. (\ref{eq.AEDD_1}) is used in the next physical solution update $\mathbf{z}^{(v+1)}$.
These processes complete one global-local iteration. The iterations proceed until the distance between the physical and material states is within a tolerance, yielding the data-driven solution denoted by the green star in Fig. \ref{fig2.3.1}(a), which ideally is the intersection between the physical manifold and material manifold in the data space.

Here, the nearest neighbors searching and locally convex reconstruction of the optimal material state are processed in the filtered (noiseless) low-dimensional embedding space, resulting in the enhanced robustness against noise and accuracy of the local data-driven solution. 

\subsubsection{AEDD local solver: Solver II}\label{sec.AEDD_locII}
Although autoencoders aim to transform the input material data to output data with maximally preserving essential features, see Fig. \ref{fig2.1.autoencoder}, the decoder functions $\mathbf{h}_{\text{dec}}$ do not exactly reproduce the given material data in the data space due to the information compression and errors inevitably introduced by training processes \cite{goodfellow2016deep}.
During \textit{Step 3} of Solver I (Eq. (\ref{eq.2.3.15})), the material embedding solution $\hat{\mathbf{z}}'^{*}_{\alpha}$ in Eq. (\ref{eq.AEDD_1}b) projecting back to the data space by decoders could involve data reconstruction errors. That is, the performance of AEDD Solver I is subject to the quality of the trained decoder functions.

To enhance the robustness and stability of data-driven computing, we propose the second AEDD local solver (Solver II) that circumvents the use of decoders and, instead, uses the interpolation scheme in Eq. (\ref{eq.shep_interp}) to perform locally convex reconstruction directly on material dataset.
The procedures of this solver are expressed as
\begin{subequations}\label{eq.AEDD_2}
    \begin{align}
        & \textit{Step 1}: \hspace{1cm} \mathbf{z}'_{\alpha} = \mathbf{h}_{\text{enc}} (\mathbf{z}_{\alpha}), \\
        & \textit{Step 2}: \hspace{1cm} \hat{\mathbf{z}}^{*}_{\alpha} = \mathcal{I} \left( \{\Psi_I(\mathbf{z}'_{\alpha});\hat{\mathbf{z}}_I\}_{I \in \mathcal{N}_k(\mathbf{z}'_{\alpha})} \right),
    \end{align}
\end{subequations}
for $\alpha = 1,...,N$, where $\hat{\mathbf{z}}_I \in \mathbb{E}$ are the material data given in the original data space. The key ingredient of this approach is that the modified locally convex reconstruction in \textit{Step 2} involves interpolation functions constructed in embedding space but interpolating material data that are in data space.
It can be viewed as a blending interpolation approach compared to that in Solve I.
The effectiveness of Solver I and II will be compared and discussed in Section \ref{sec3.2.1}.

In both Solver I and II, the interpolation functions $\Psi_I(\mathbf{z}'_{\alpha})$ are evaluated on the embedding space related to the embedded physical state $\mathbf{z}'_{\alpha}$ and its $k$ nearest neighbors in the material embedding data $\{\hat{\mathbf{z}}'_I\}_{I \in \mathcal{N}_k(\mathbf{z}'_{\alpha})} \subset \mathbb{E}'$. In Solver II, however, these functions are weighted on the un-projected material data $\{\hat{\mathbf{z}}_I\}_{I \in \mathcal{N}_k(\mathbf{z}'_{\alpha})} \subset \mathbb{E}$ corresponding to the $k$ selected neighbors.
Because the locally convex reconstruction in Eq. (\ref{eq.AEDD_2}b) gives the optimal material data solution immediately in the data space, the decoder is avoided in Solver II.

Fig. \ref{fig2.3.1}(b) shows a schematic of the proposed data-driven computing based on Solver II. 
Taking the $v$-th global-local iteration as an example, the physical state $\mathbf{z}^{(v)}$ obtained from the global physical step (Eq. (\ref{eq.2.3.11})) is mapped to the embedding space $\mathbf{z}'^{(v)}$ by the encoder. 
In \textit{Step 2} of Solver II, the same $k$ nearest neighbors search is performed on the embedding space $\mathcal{E}'$, while their corresponding material data in the original dataset are used in the data reconstruction via Eqs. (\ref{eq.shep_interp})--(\ref{eq.shep_kern}).
As shown in Fig. \ref{fig2.3.1}(b), the locally convex reconstruction of the material embedding state $\mathbf{z}'^{(v)}$ can be directly performed with
the material data $\{\hat{\mathbf{z}}_I\}_{I \in \mathcal{N}_k(\mathbf{z}'_{\alpha})}$ in the data space by using data indices, yielding the optimal material solution $\hat{\mathbf{z}}^{*(v)}$.

%

It is worth emphasizing that in comparison with the LCDD approach \cite{he2019physics}, the key difference in the proposed solver is that the neighbor search and data reconstruction are performed on the embedding space $\mathcal{E}'$, a lower-dimensional space constructed by the pre-trained encoder function. Thus, the proposed AEDD with Solver II can be considered as a enhanced generalization of LCDD for high-dimensional material data.
We use this approach as the default AEDD method, unless stated otherwise.





\section{Numerical results}\label{sec.results}
In this section, the proposed AEDD approach is first tested on a cantilever beam \textcolor{blue}{using synthetic material data generated by constitutive laws}. In this example, the effects of several factors on autoencdoers and the resulting AEDD data-driven solutions are investigated, including the size, sparsity, and the noise level of material datasets, neural network initialization during autoencoder training, and autoencoder architectures, aiming to validate the robustness and reliability.
In the second subsection, AEDD is applied to modeling biological tissues using experimental data measured from heart valve tissues to demonstrate the enhanced generalization capability.


For simplicity, we consider homogeneous material in the following numerical examples, and thus the same material dataset, e.g., $\mathbb{E}=\{\hat{\mathbf{z}}_I\}_{I=1}^M=\{({\hat{\mathbf{E}}}_{I},{\hat{\mathbf{S}}}_{I})\}_{I=1}^M$, with $M$ data points, is used for all integration points.

\subsection{Cantilever beam: Verification of the AEDD method}\label{sec3.1}
\begin{figure}[!ht]
    \centering
    \includegraphics[width=0.7\linewidth]{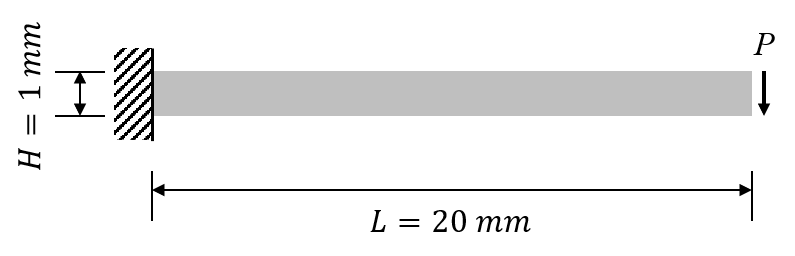}
    \caption{Schematic of a cantilever beam model subjected to a tip shear load, where $P = 10 EI/L^2$, and $I=H^3/12$.}\label{fig.bean.schematic}
\end{figure}

To verify the proposed AEDD framework (with Solver II in Section \ref{sec.AEDD_locII} by default), a cantilever beam subjected to a tip shear load is analyzed, as shown in Fig. \ref{fig.bean.schematic}. The Saint Venant-Kirchhoff phenomenological model with Young's modulus $E=4.8 \times 10^3 N/mm^2$ and Poisson's ratio $\nu=0$ is used to generate material datasets for training autoencoders.
The problem domain is discretized with $41\times5$ randomly distributed nodes.
The data-driven analysis is performed with 10 equal loading steps under a plane-strain condition.
Following the same setting in \cite{he2019physics}, the
weight matrix $\hat{\mathbb{C}}$ used in the distance metric (Eqs. (\ref{eq.2.3.7})-(\ref{eq.2.3.8})) and the physical solver (Eq. (\ref{eq.global})-(\ref{eq.stress_updatre})) is defined as
\begin{equation}\label{eq.3.1.1}
    \hat{\mathbb{C}} = \frac{E}{1-\nu^2}
    \begin{bmatrix}
    1 & 0 & 0 \\
    0 & 1 & 0 \\
    0 & 0 & (1-\nu)/2
    \end{bmatrix}
\end{equation}

\subsubsection{Preparation of material datasets}\label{sec3.1.1}
\begin{figure}[!ht]
\centering
    \begin{subfigure}{0.49\textwidth}
        \centering
        \includegraphics[width=1\linewidth]{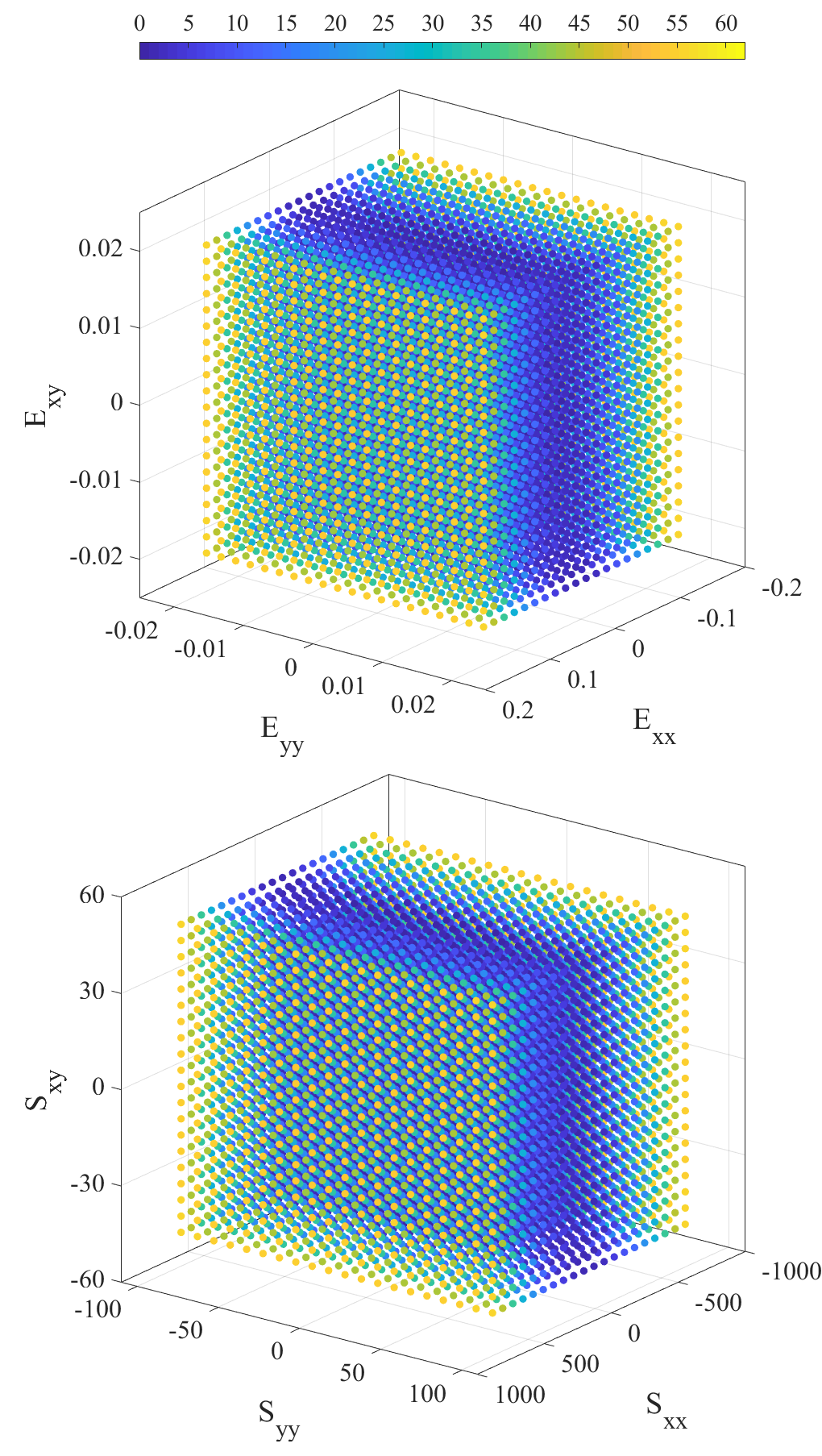}
        \caption{Noiseless}
    \end{subfigure}
    \begin{subfigure}{0.49\textwidth}
        \centering
        \includegraphics[width=1\linewidth]{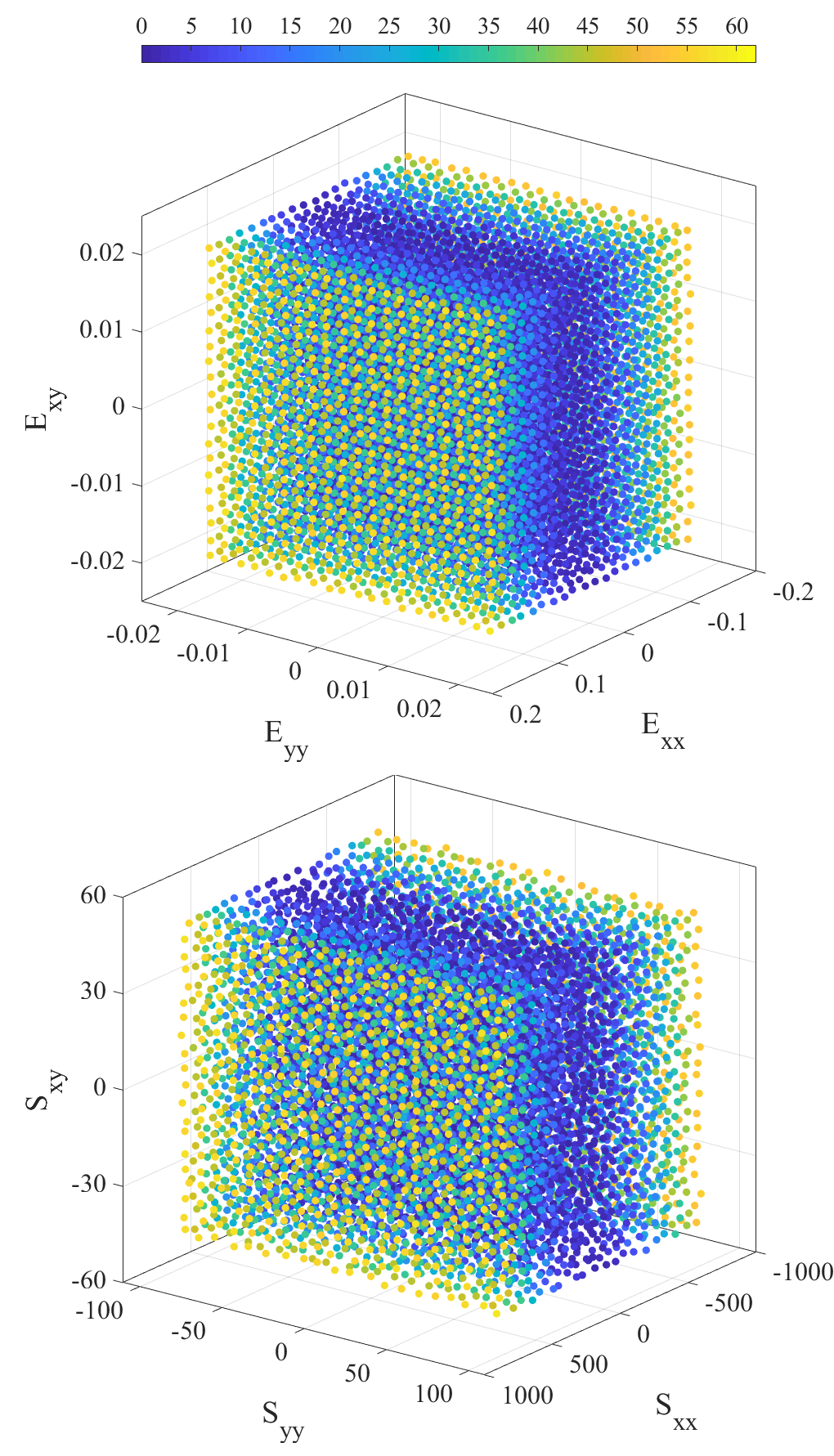}
        \caption{Noisy}
    \end{subfigure}
\caption{Material dataset with a size of $M=20^3$: (a) Noiseless; (b) Noisy; Top: strain components; Bottom: stress components}
\label{fig.beam.noisyData1}
\end{figure}

\begin{figure}[!ht]
\centering
    \begin{subfigure}{0.32\textwidth}
        \centering
        \includegraphics[width=1\linewidth]{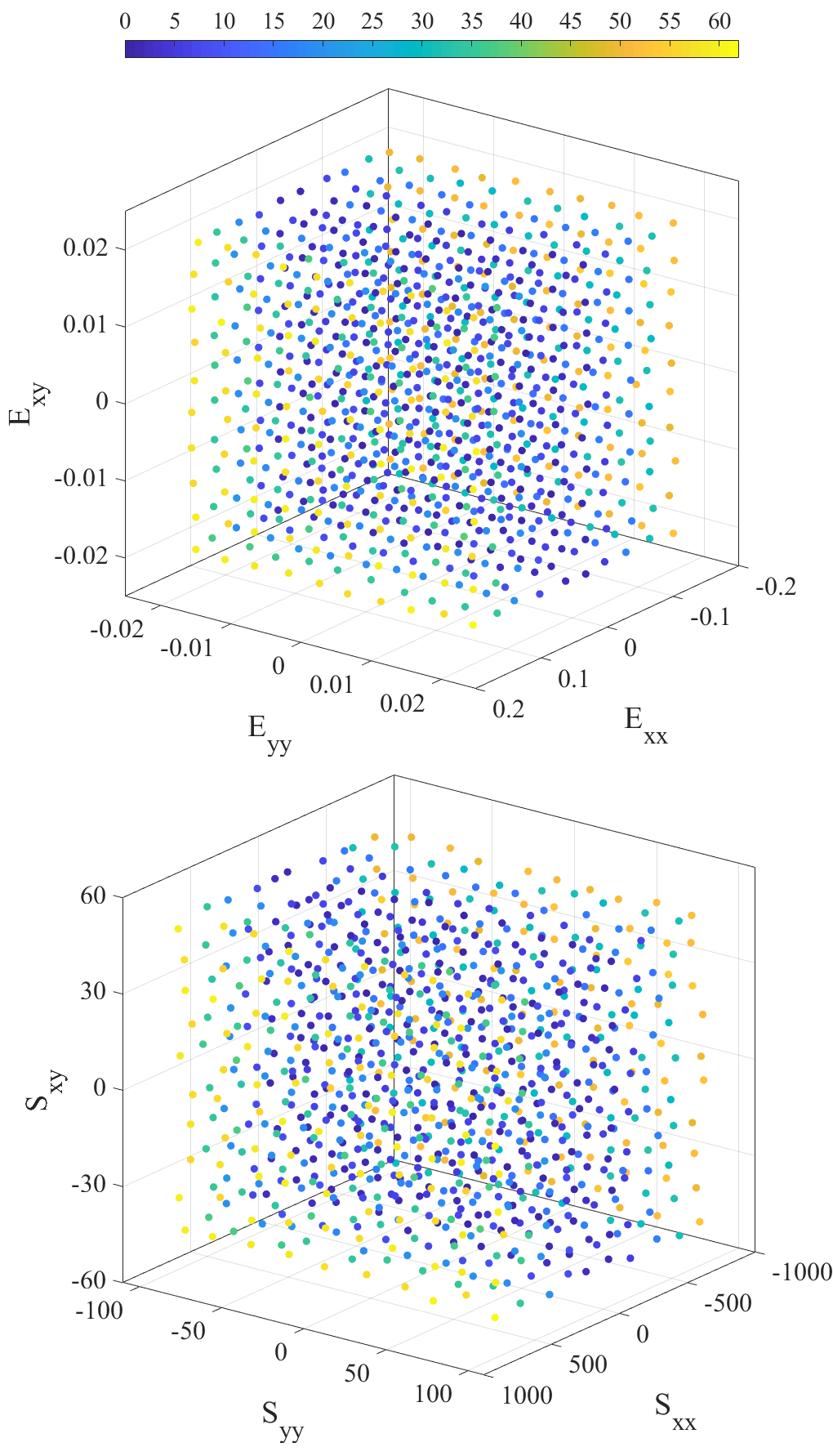}
        \caption{}
    \end{subfigure}
    \begin{subfigure}{0.32\textwidth}
        \centering
        \includegraphics[width=1\linewidth]{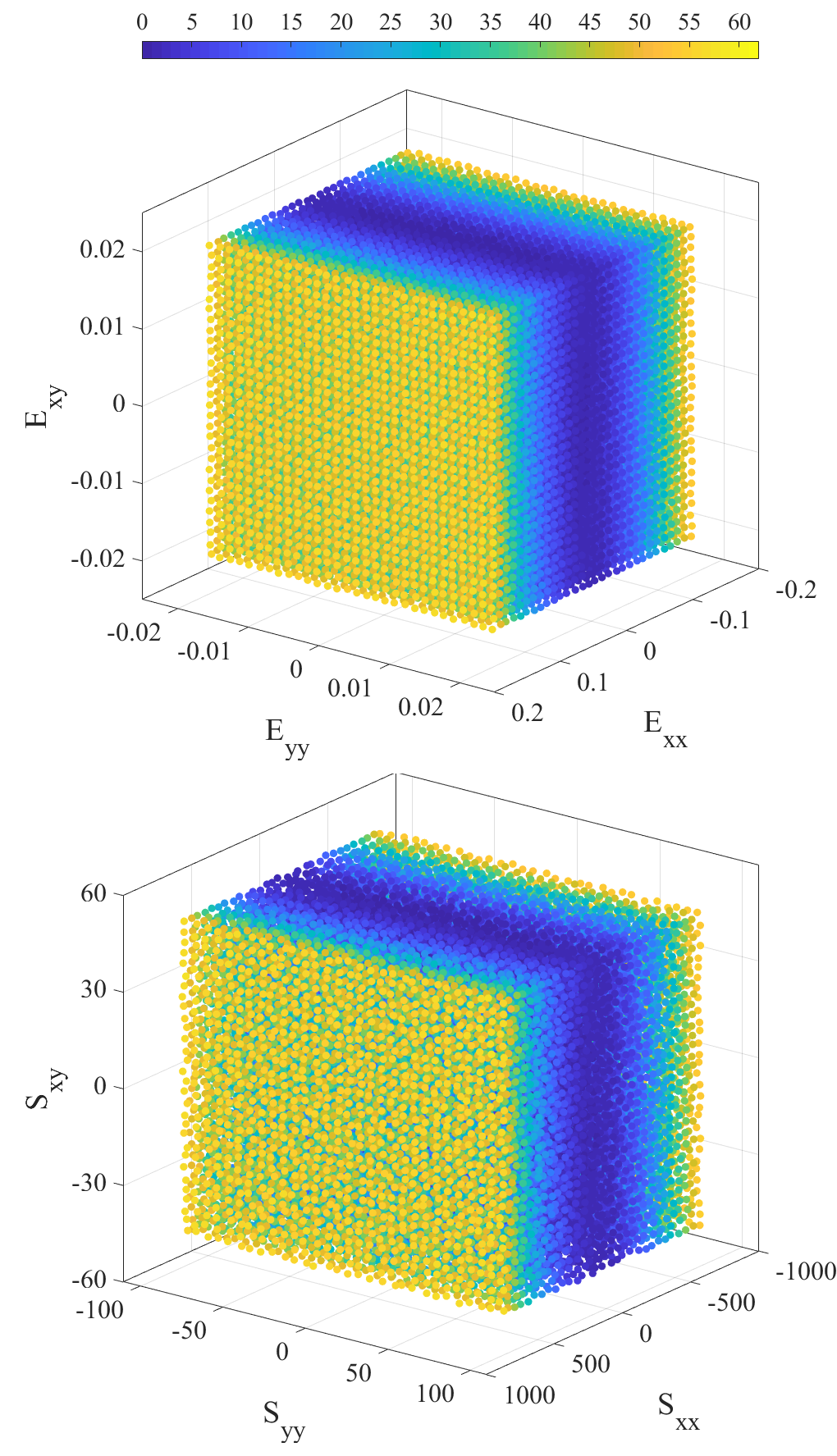}
        \caption{}
    \end{subfigure}
    \begin{subfigure}{0.32\textwidth}
        \centering
        \includegraphics[width=1\linewidth]{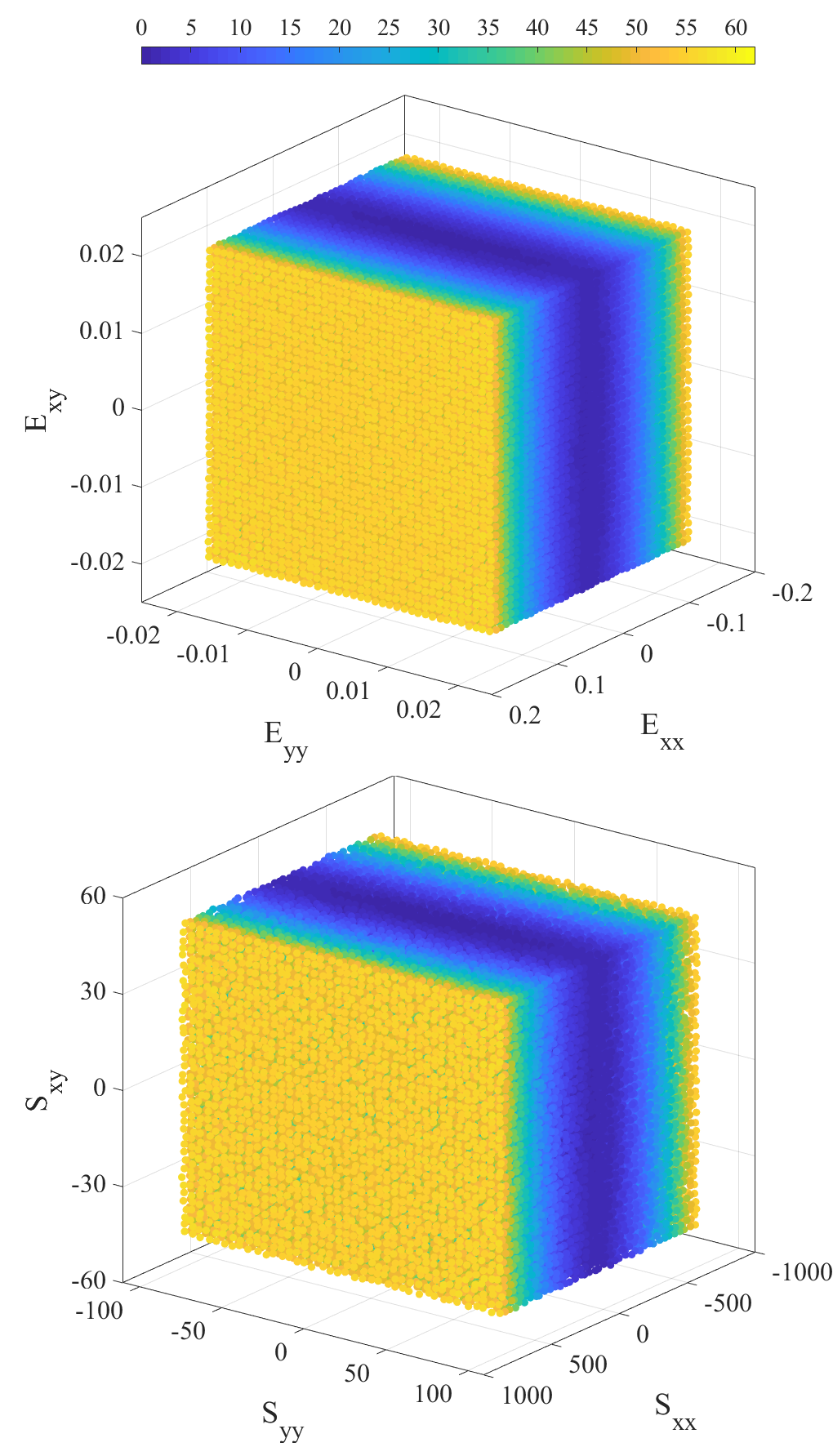}
        \caption{}
    \end{subfigure}
\caption{Noisy material datasets with: (a) $M=10^3$; (b) $M=30^3$; (c) $M=40^3$; Top: strain components; Bottom: stress components}
\label{fig.beam.noisyData2}
\end{figure}

To assess the robustness and convergence property of AEDD against noise presented in the given material datasets,
four manufactured noisy material datasets approximating the Saint Venant-Kirchhoff phenomenological model with different data sizes, i.e. $M=10^3$, $20^3$, $30^3$, and $40^3$, are considered. The generation procedure of these noisy datasets is described below.
%
First, an noiseless dataset,
$ \bar{\mathbb{E}} =\{\bar{\mathbf{z}}_{I}\}_{I=1}^M$
is generated, where each Green-Lagrangian strain component is uniformly distributed within the range $[-0.02, 0.02]$ and the 2nd-PK stress components are obtained by using the elastic tensor in Eq. (\ref{eq.3.1.1}) that relates strain to stress. The example with $M=20^3$ is shown in Fig. \ref{fig.beam.noisyData1}(a), where the strain and stress components are displayed separately for visualization.
Following \cite{kirchdoerfer2017data,he2019physics}, Gaussian perturbations scaled by a factor dependent on the size of datasets, $0.4\bar{\mathbf{z}}_{max}/\sqrt[3]{M}$, are added to each component of the noiseless dataset $ \bar{\mathbb{E}}$
to obtain the associated noisy datasets
$\mathbb{E} =\{\hat{\mathbf{z}}_{I}\}_{I=1}^M$, where $\bar{\mathbf{z}}_{max}$
is a vector of the maximum values for each component among the noiseless dataset.
The noisy dataset corresponding to $M=20^3$ is shown in Fig. \ref{fig.beam.noisyData1}(b).
Fig. \ref{fig.beam.noisyData2} shows the other three noisy material datasets.

\subsubsection{Effects of autoencoder architecture and initialization}\label{sec3.1.2}
In order to assess the effects of initialization during training and network architectures on autoencoders' accuracy and robustness,
five random initializations and four architectures of autoencoders are considered.
For the given noisy material datasets associated with this plane-strain cantilever beam problem, it is observed that autoencoders with an embedding layer of the dimension $p=1$ or $p=2$ could not capture a meaningful embedding representation.
This is consistent to the observation in \cite{he2019physics} where the number of neighbor points to construct the locally convex embedding is suggested to be larger than the number of intrinsic dimensionality, which is 2 of the employed linear elastic database.
Hence, it requires the embedding dimension to be greater than 2.
As described in Section \ref{sec.auto_train}, the encoder architecture is used to represent the autoencoder architecture. Four encoder architectures, $6-4-3$, $6-5-4$, $6-5-4-3$, and $6-10-8-5$, are considered in the following tests.

\begin{figure}[!ht]
\centering
    \begin{subfigure}{0.24\textwidth}
        \centering
        \includegraphics[width=1\linewidth]{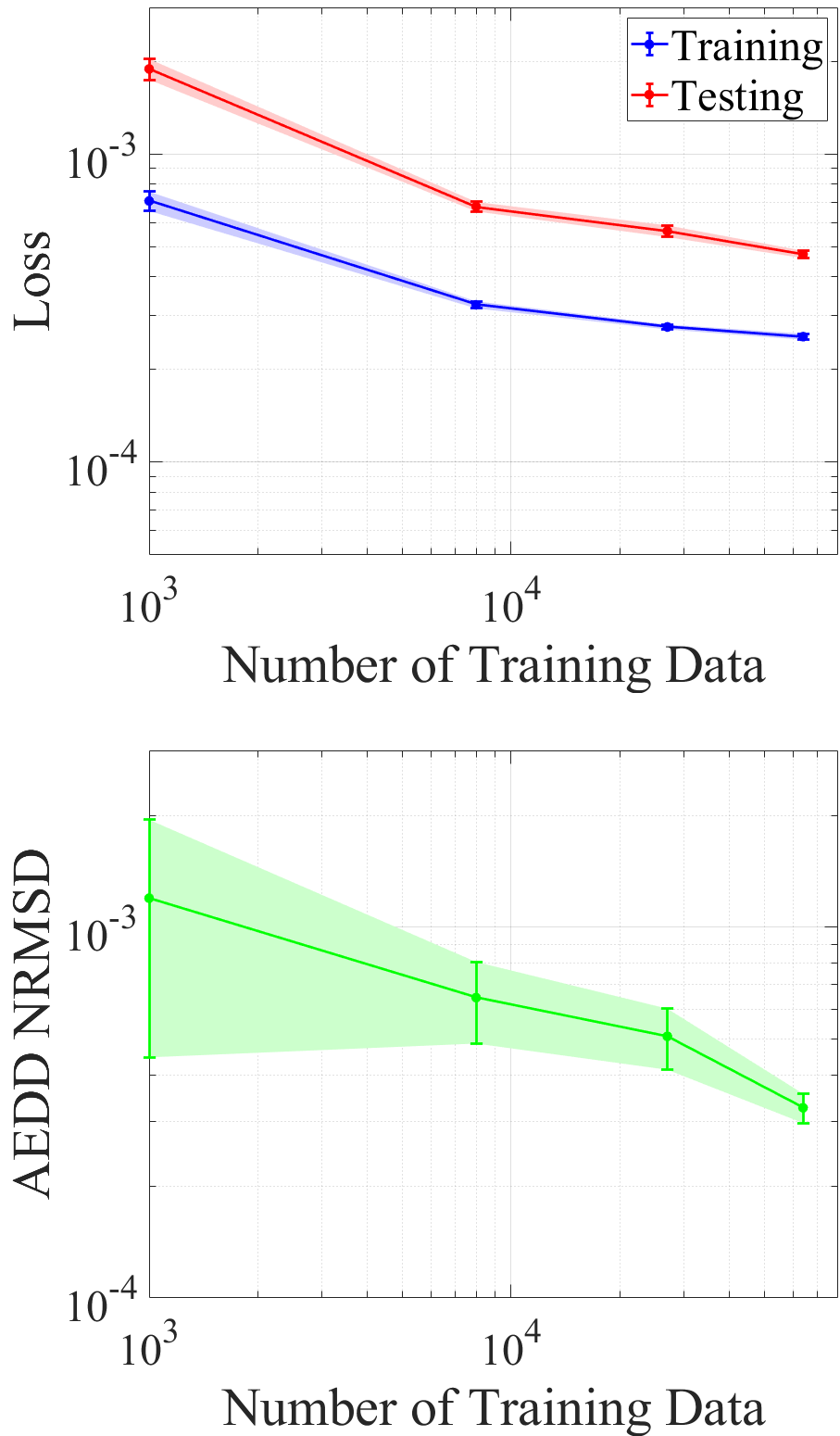}
        \caption{Encoder: 6-4-3}
    \end{subfigure}
    \begin{subfigure}{0.24\textwidth}
        \centering
        \includegraphics[width=1\linewidth]{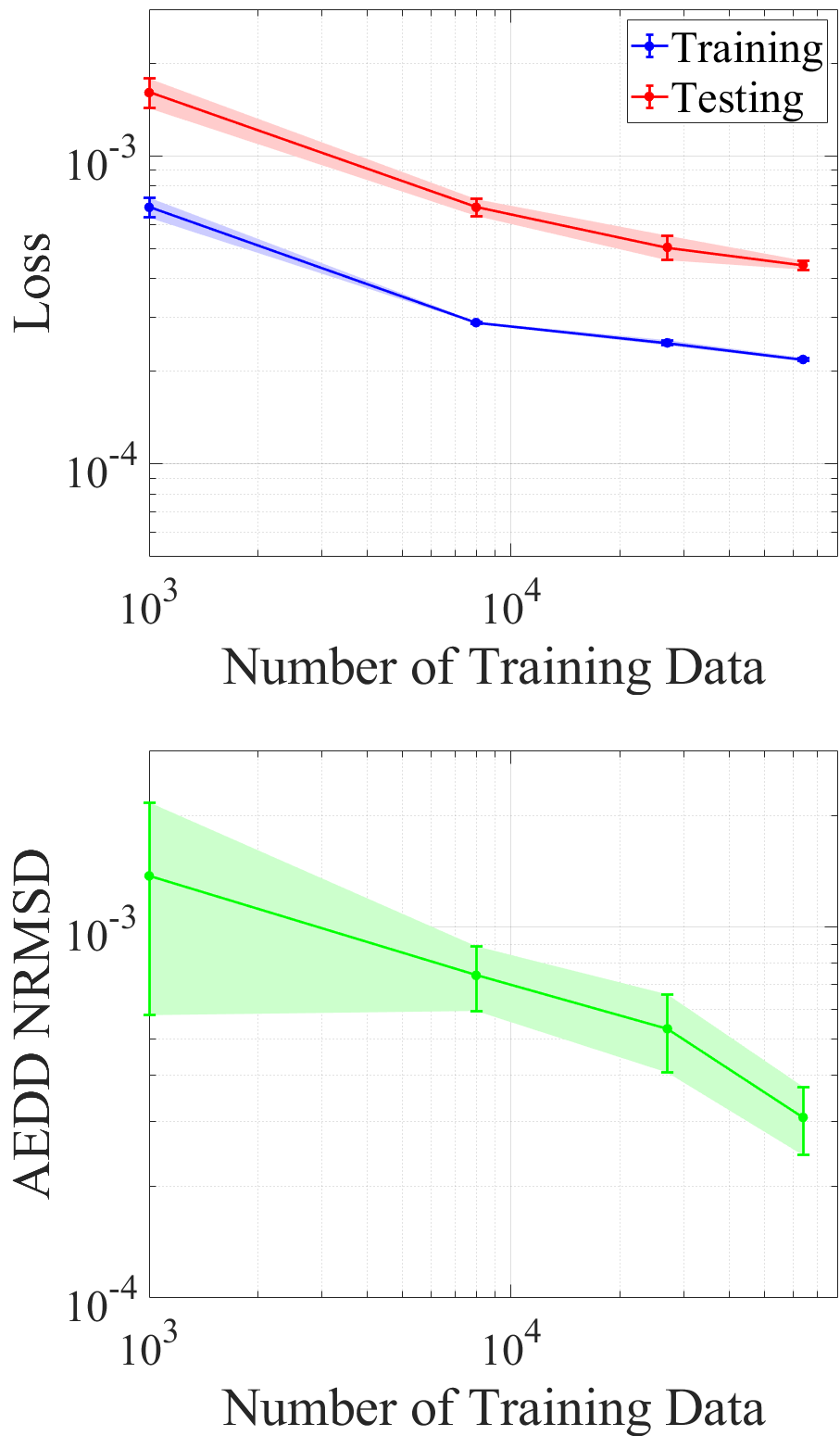}
        \caption{Encoder: 6-5-4}
    \end{subfigure}
    \begin{subfigure}{0.24\textwidth}
        \centering
        \includegraphics[width=1\linewidth]{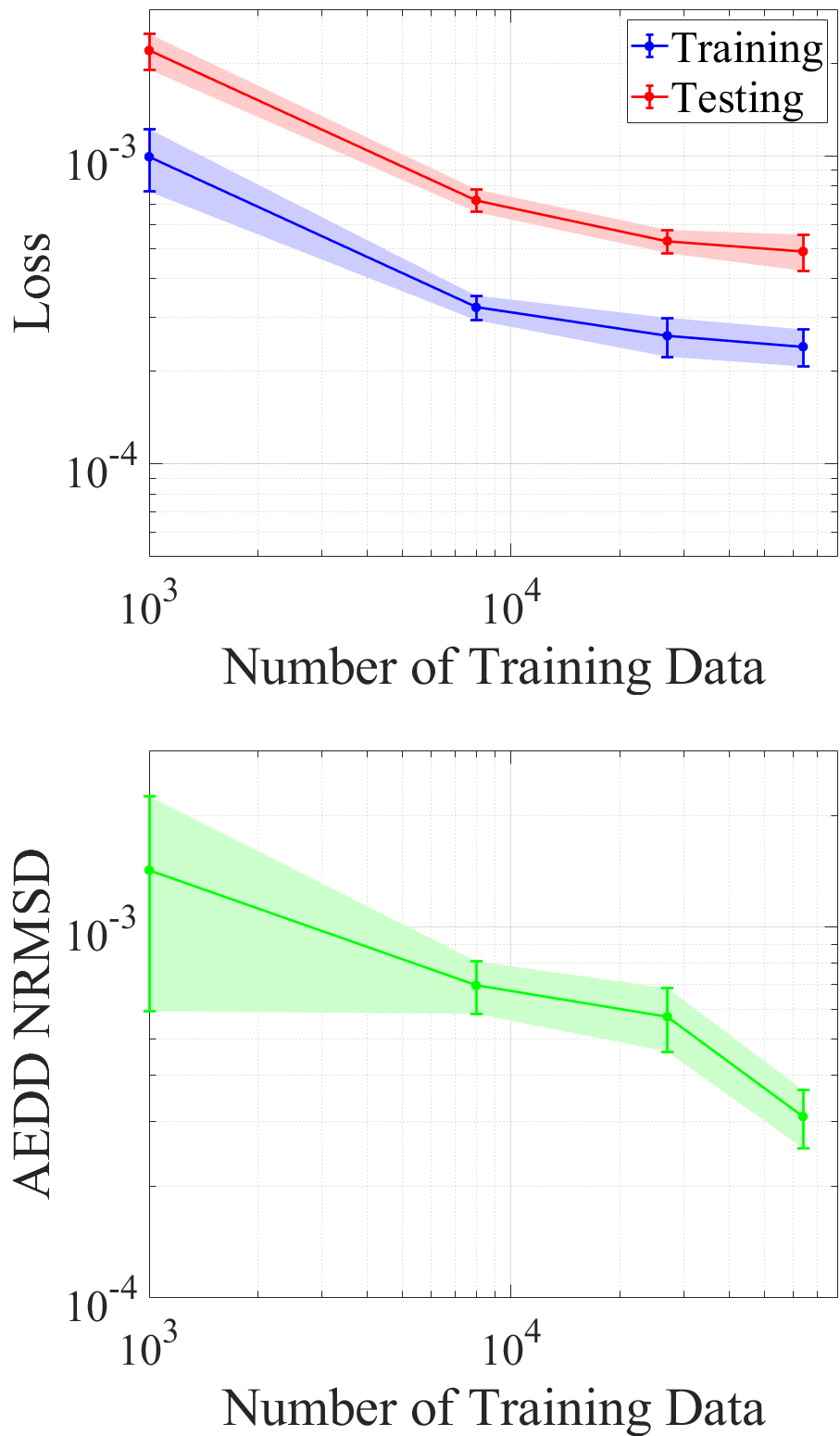}
        \caption{Encoder: 6-5-4-3}
    \end{subfigure}
    \begin{subfigure}{0.24\textwidth}
        \centering
        \includegraphics[width=1\linewidth]{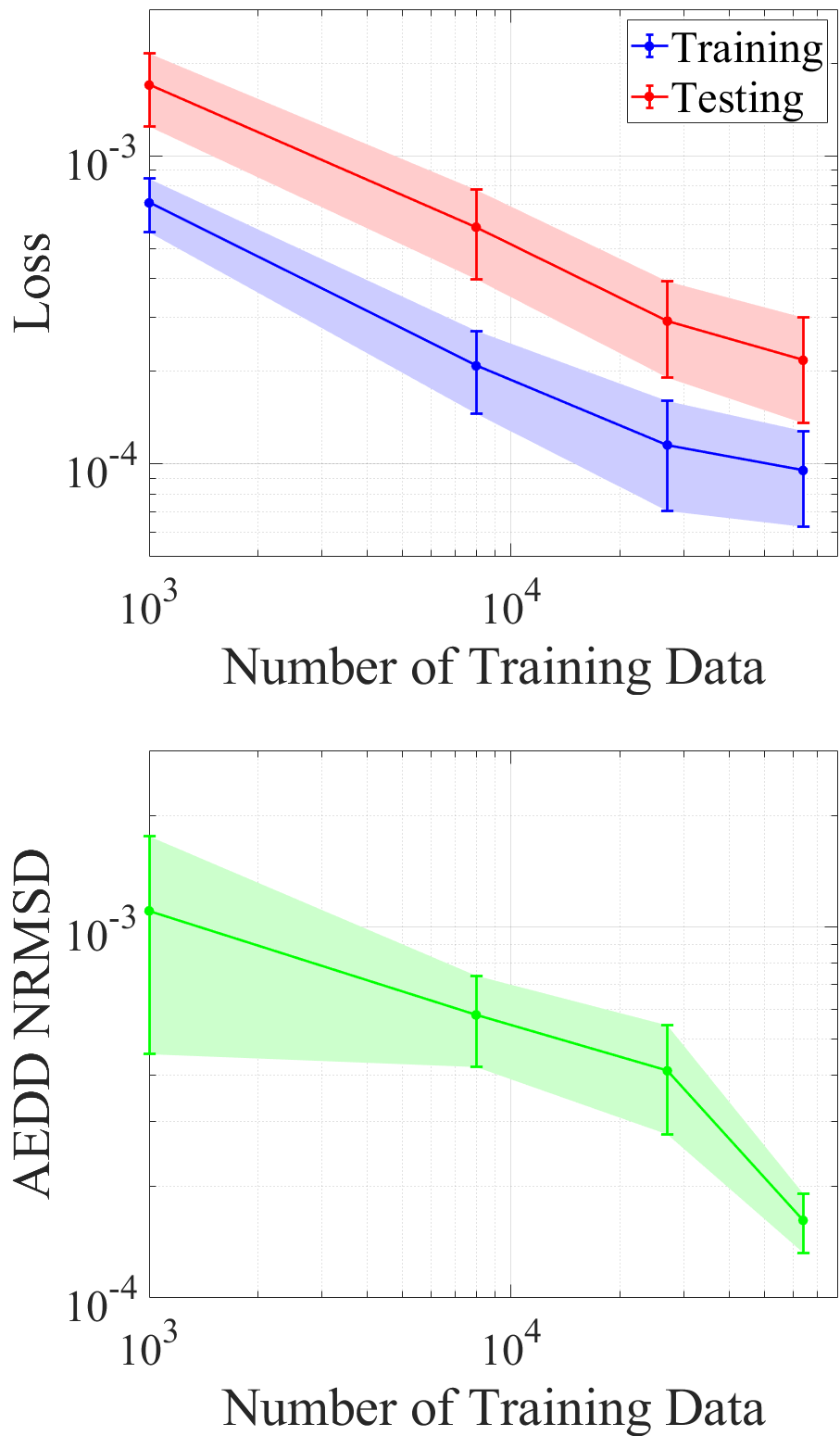}
        \caption{Encoder: 6-10-8-5}
    \end{subfigure}
\caption{Error curves (mean with standard deviations shaded) of four different encoder architectures: (a) $6-4-3$; (b) $6-5-4$; (c) $6-5-4-3$; (d) $6-10-8-5$ Top: final training and testing losses of autoencoders; Bottom: \textit{NRMSD} between AEDD and constitutive model-based solutions}
\label{fig.beam.errorCurve}
\end{figure}

The first row in Fig. \ref{fig.beam.errorCurve} shows the error curves (mean with standard deviations shaded) of the final training and testing losses against the size of training dataset for different autoencoder architectures.
Here, the noisy material datasets of different sizes, $M=10^3$, $20^3$, $30^3$, and $40^3$, that defined in Section \ref{sec3.1.1} are used for training the autoencoders, where autoencoders are trained with five random initialization for each case.
Besides, to fairly compare the testing errors between the autoencoders trained with various sizes of training data, we use the same test dataset consisting of 729 material data points that are generated from the same procedure in Section \ref{sec3.1.1} but not included in the given material datasets.

As we can see, all the selected autoencoders converge well, yielding smaller training and testing errors as the size of material dataset increases. Moreover, it is observed that the autoencoder with a larger architecture could lead to greater variation due to training randomness, indicated by the standard deviations. This is because the training algorithms used to minimize the loss function in Eq. (\ref{eq.loss_auto}) do not guarantee global minimization, and a larger DNN with more trainable parameters may cause higher randomness. However, the trained results shown here are satisfactory due to the employment of regularization. It also shows that the training and testing losses decrease as the dimension of the embedding layer increases.

\subsubsection{Data-driven modeling results}\label{sec3.1.3}
The data-driven solution is compared with the constitutive model-based reference solution using Eq. (\ref{eq.3.1.1}).
To better assess the accuracy of AEDD solutions, a normalized root-mean-square deviation (\textit{NRMSD})
is introduced
\begin{equation}\label{eq.3.1.2}
    \textit{NRMSD} = \sqrt{\sum_i^{N_{eval}}\frac{(\bar{w}_i^{AEDD} - \bar{w}_i^{ref})^2}{N_{eval}}} / (PL^2/EI),
\end{equation}
where $N_{eval}=200$ is the number of evaluation points, $\bar{w}_i^{AEDD}$ and $\bar{w}_i^{ref}$ are the normalized tip deflection obtained by AEDD and model-based reference solutions, respectively. In this cantilever beam case, the normalized tip deflection $\bar{w}_i=w_i/L$ is obtained at the maxiumn loading, i.e. $PL^2/EI = 10$, where $L$ and $H$ are the length and the width of the beam, respectively, and $I=H^3/12$, see Fig. \ref{fig.bean.schematic}.

\begin{figure}[!ht]
\centering
    \begin{subfigure}{0.49\textwidth}
        \centering
        \includegraphics[width=1\linewidth]{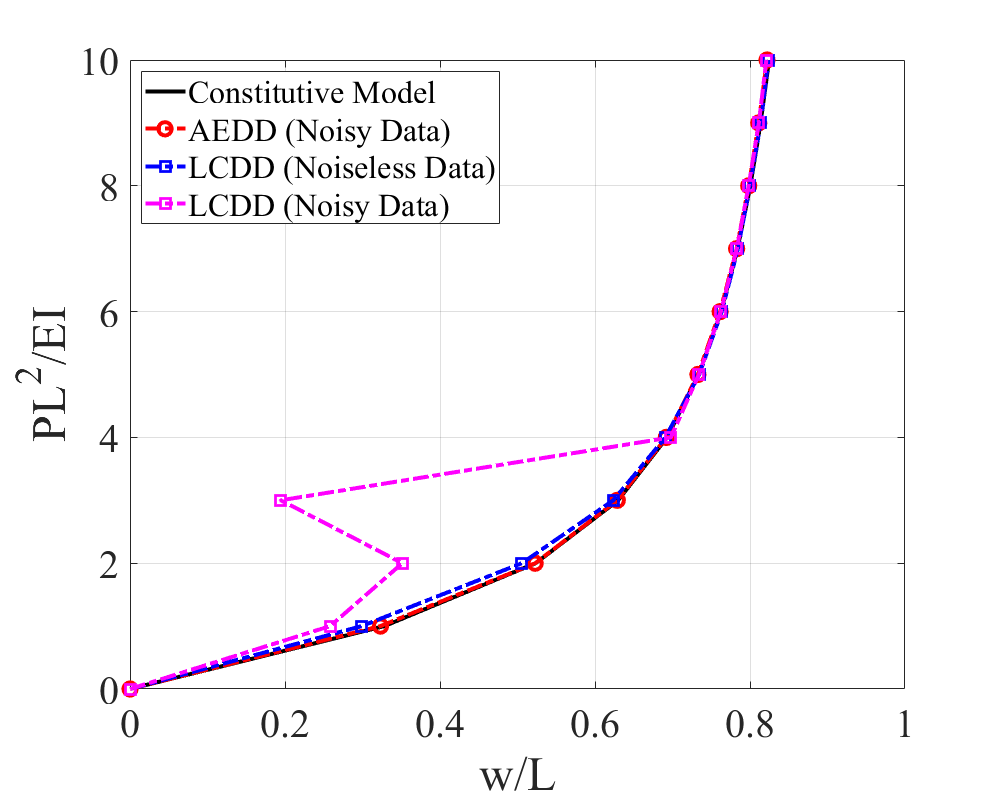}
        \caption{}
    \end{subfigure}
    \begin{subfigure}{0.49\textwidth}
        \centering
        \includegraphics[width=1\linewidth]{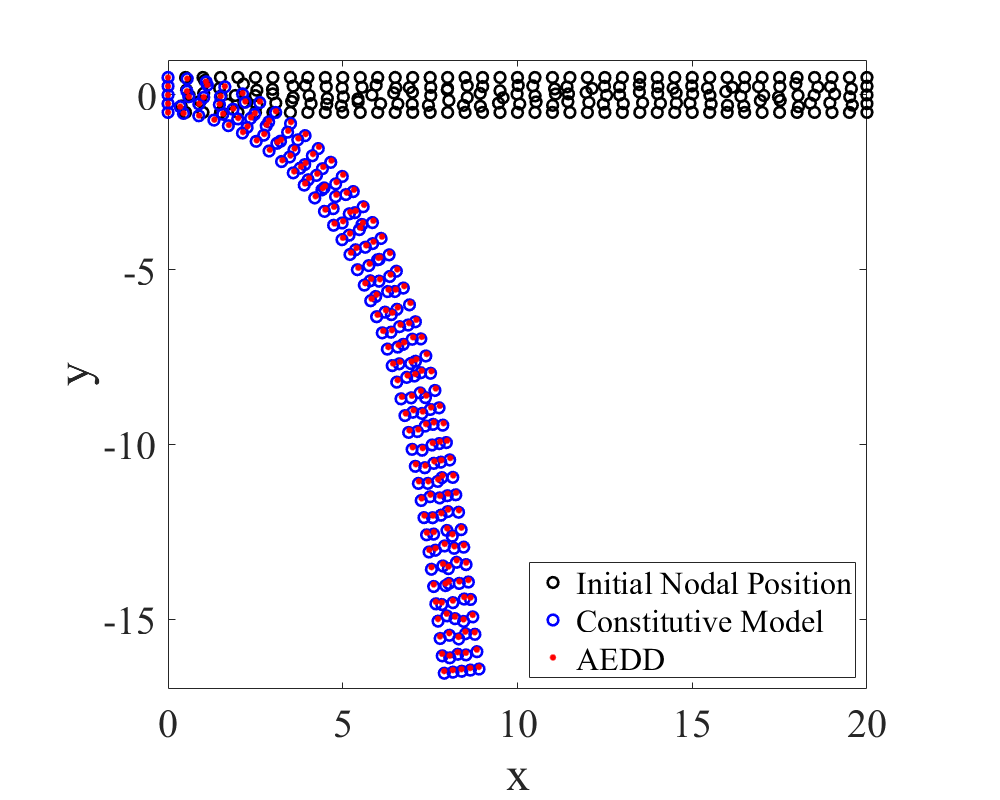}
        \caption{}
    \end{subfigure}
\caption{Comparison of constitutive model-based, LCDD, and AEDD solutions: (a) normalized tip deflection-loading; (b) initial and final nodal positions; The AEDD solution is obtained from using autoencoders trained with a material dataset of size $M=40^3$.}
\label{fig.beam.deflection}
\end{figure}

\begin{figure}[!ht]
\centering
    \begin{subfigure}{0.49\textwidth}
        \centering
        \includegraphics[width=1\linewidth]{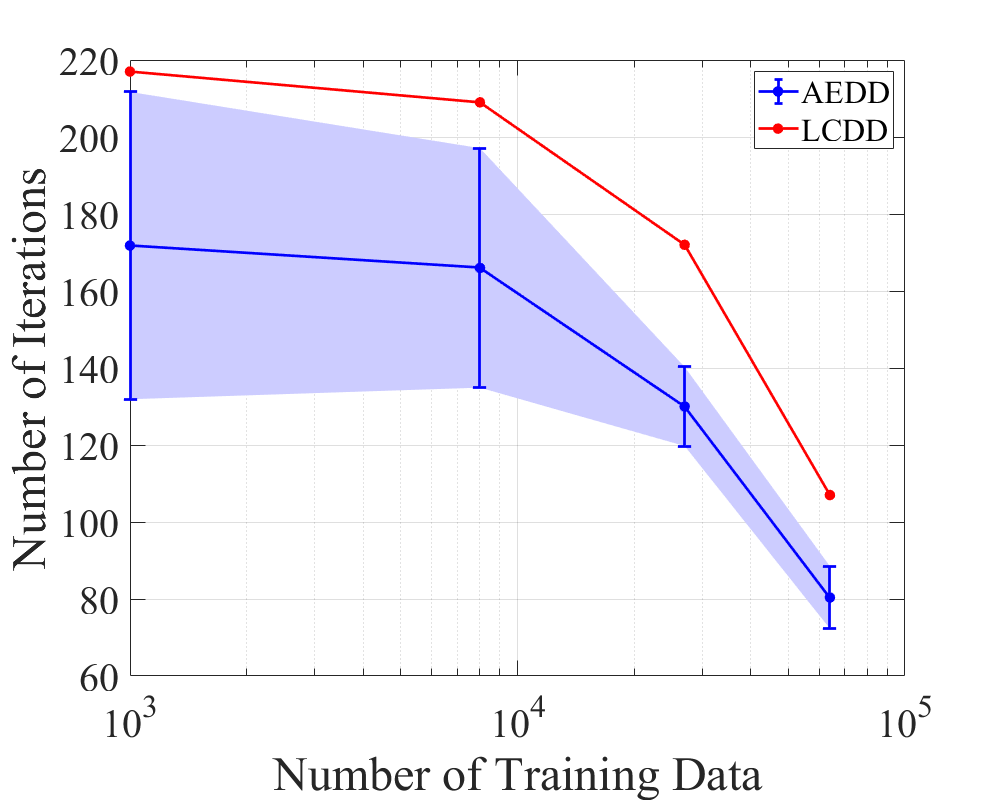}
        \caption{}
    \end{subfigure}
    \begin{subfigure}{0.49\textwidth}
        \centering
        \includegraphics[width=1\linewidth]{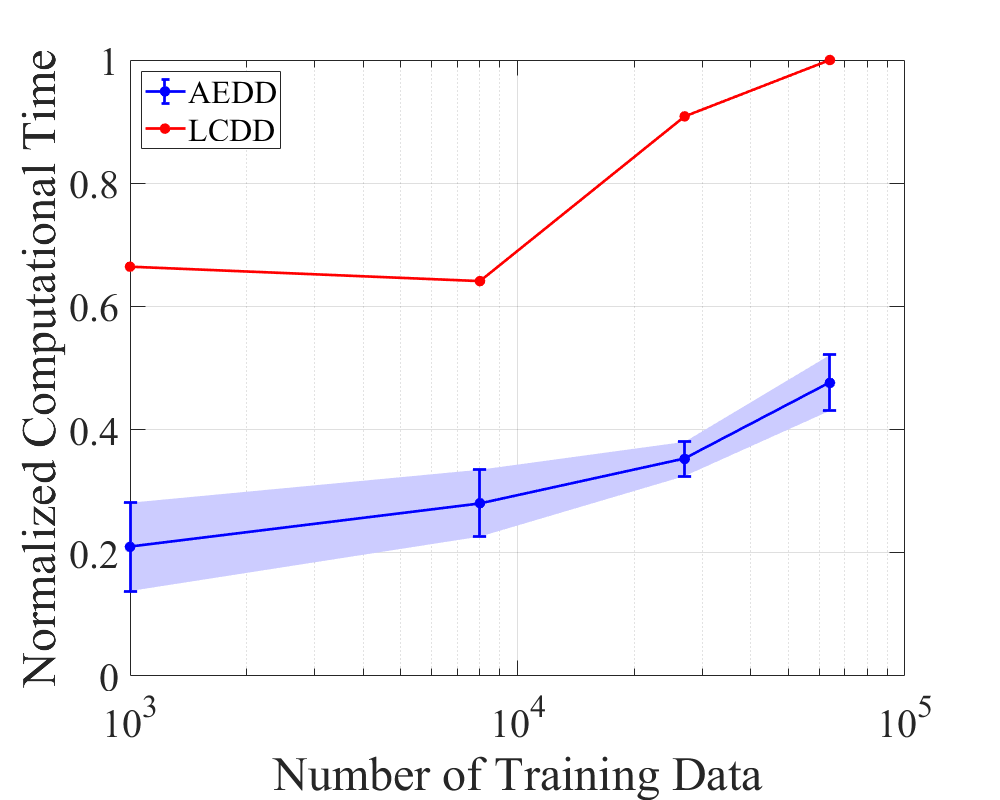}
        \caption{}
    \end{subfigure}
\caption{Comparison of LCDD and AEDD: (a) Number of iterations against number of training data; (b) normalized computational time against number of training data. The noisy datasets and the encoder architecture of 6-4-3 are employed in this test.}
\label{fig.beam.iter/time}
\end{figure}
The trained autoencoders corresponding to different material datasets and architetures are then applied to data-driven simulations, where the number of nearest neighbors used in locally convex reconstruction of the data-driven solver is set as 6.
\textit{NRMSD} of data-driven solutions with respect to the model-based reference solution is given in the bottom row of Fig. \ref{fig.beam.errorCurve}.
For all architectures, it can be observed that the AEDD solutions (both mean values and variation) improve as the number of training data increases, which suggests a good convergence property.
Although using an embedding dimension of 5 (encoder: $6-10-8-5$) yields the highest accuracy, the AEDD solutions obtained from using an embedding dimension of 3 and 4 are satisfactory. It also shows that the overall patterns of error convergence in \textit{NRMSDs} are similar across different encoder architectures using the same size of training dataset, indicating that the AEDD solutions are not sensitive to the width and depth of the encoder architecture as long as
autoencoders of a sufficiently large size are used.
Considering that using a more complex encoder architecture with a larger embedding dimension would increase computational cost in data-driven computing, an encoder architecture $6-4-3$ is used in the numerical examples.

Fig. \ref{fig.beam.deflection} shows that the normalized tip deflection-loading curve predicted by the proposed AEDD method agrees well with the model-based reference. The noisy data set of size $M=40^3$ is used in this case. The results obtained by LCDD are also provided for comparison in Fig. \ref{fig.beam.deflection}(a), where a few loading steps yield divergent data-driven solutions when the noisy material data is employed.
On the other hand, the AEDD method stays robust even with noisy data employed.
It is also worth noting that when using Solver I (Eq. (\ref{eq.AEDD_1})) in AEDD, we also observe unconverged solutions (which are not reported in the Figures). We attribute this to the information loss caused by the decoder functions. 
On the other hand, Solver II with the convex interpolation functions defined in the embedding space and the material data points in data space yields stable solutions.

The comparison of AEDD and LCDD with respect to the iteration number and the computational cost are given in Fig. \ref{fig.beam.iter/time}. In this case, the architecture of $6-4-3$ is used. While the number of iterations for convergence varies in AEDD due to the non-uniqueness of autoencoder training, it generally requires less data-driven iterations than LCDD to achieve converged solutions, as shown in Fig. \ref{fig.beam.iter/time}(a). This is because a more generalized embedding space is used in AEDD for computing the local material solution. We also observe that with less noisy material data, the required iteration number decreases regardless of the increase in data size, an attractive property for data-driven computing.
Moreover, because the data search and the convexity-preserving interpolation in AEDD local solver are performed in the low-dimensional embedding space instead of the high-dimensional data space, the computational cost of AEDD is substantially reduced compared to LCDD, as shown in Fig. \ref{fig.beam.iter/time}(b).

\subsubsection{Data-driven modeling with sparse noisy datasets}\label{sec3.1.4}
\begin{figure}[!ht]
\centering
    \begin{subfigure}{0.32\textwidth}
        \centering
        \includegraphics[width=1\linewidth]{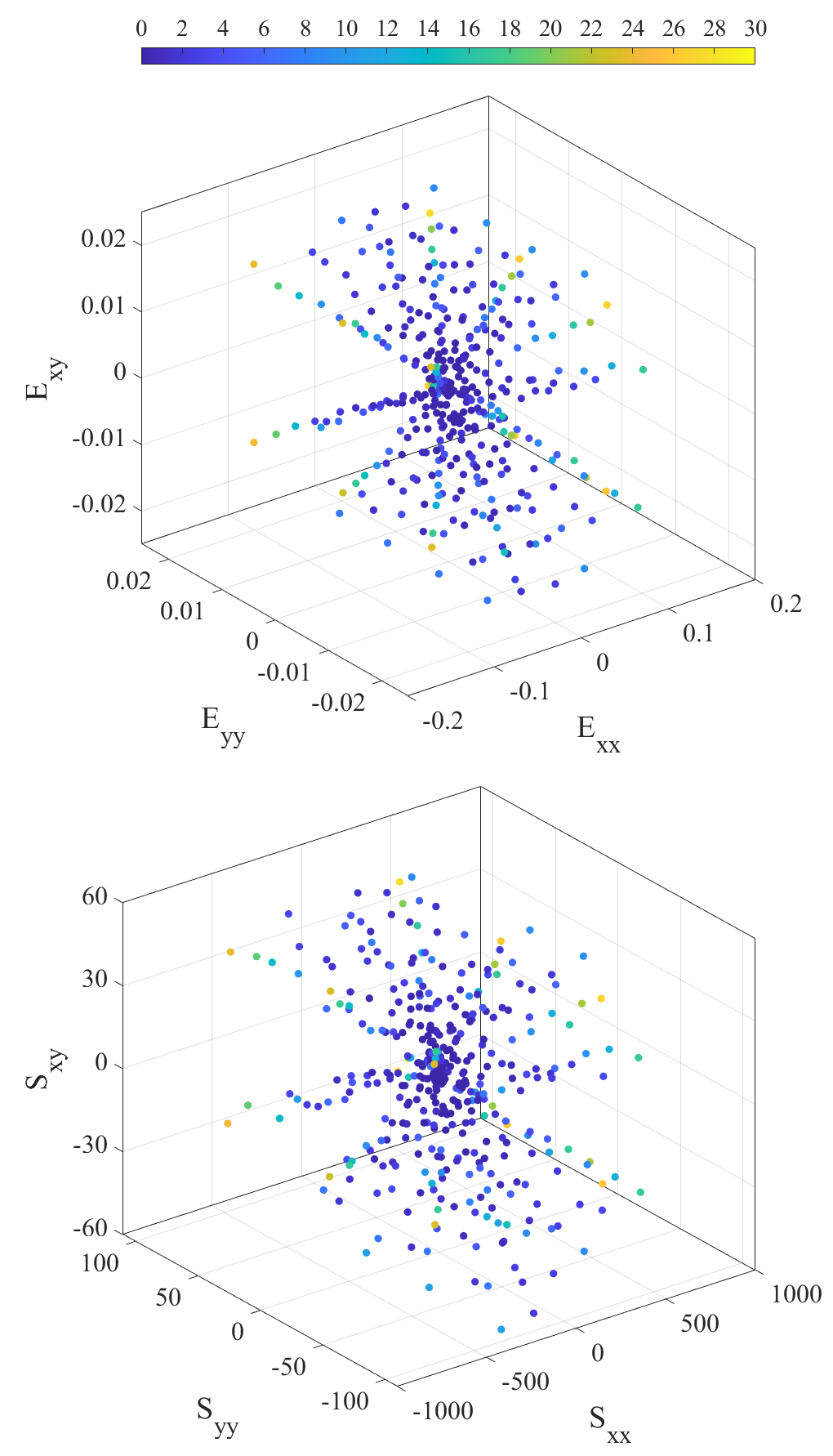}
        \caption{}
    \end{subfigure}
    \begin{subfigure}{0.32\textwidth}
        \centering
        \includegraphics[width=1\linewidth]{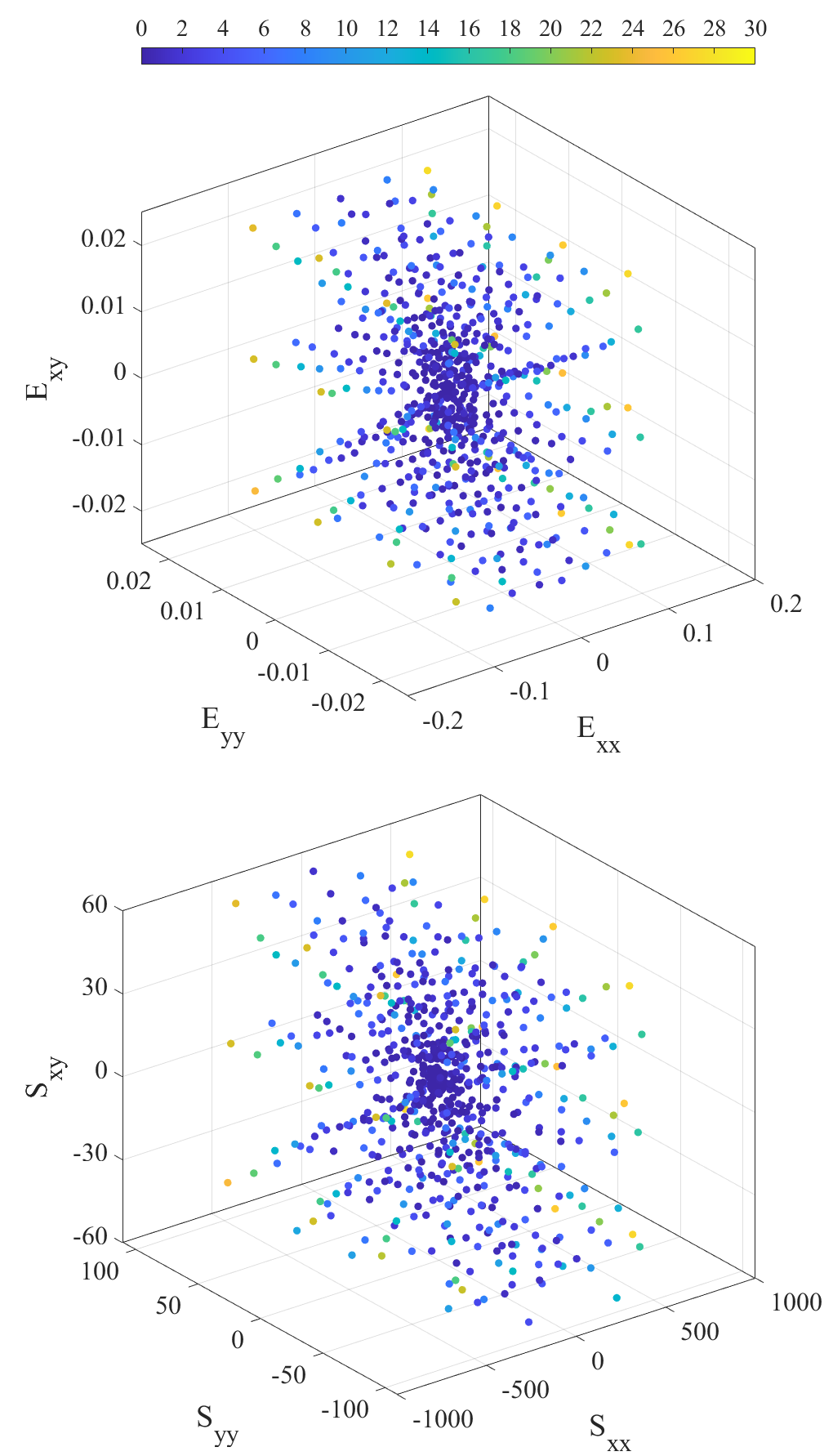}
        \caption{}
    \end{subfigure}
    \begin{subfigure}{0.32\textwidth}
        \centering
        \includegraphics[width=1\linewidth]{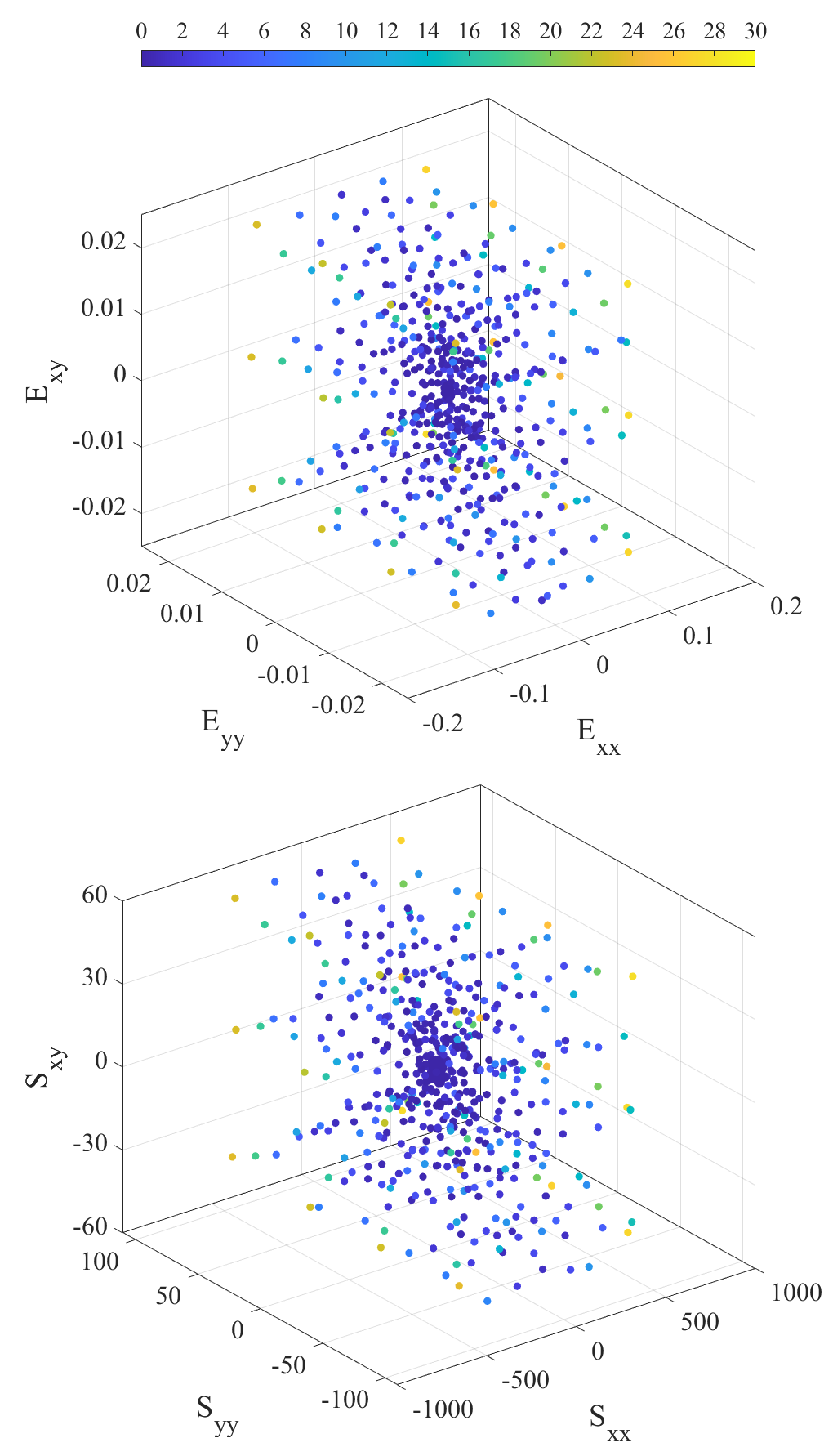}
        \caption{}
    \end{subfigure}
\caption{Sparse noisy material datasets: (a) sparse dataset 1; (b) sparse dataset 2; (c) sparse dataset 3; Top: strain components; Bottom: stress components}
\label{fig.beam.sparse_dataset}
\end{figure}

To evaluate the performance of the proposed AEDD approach when datasets are sparse, three noisy material datasets (Table \ref{tab.sparse_dataset}) are generated in a similar manner as described in Section \ref{sec3.1.1} but with fewer data points compared to Fig. \ref{fig.beam.noisyData2}.
First, several loading paths are selected with uniformly distributed Green-Lagrangian strains for each of the loading paths. The corresponding 2nd PK stresses are generated using the elastic tensor given in Eq. (\ref{eq.3.1.1}). Consequently, the sparse noisy material datasets are given in Fig. \ref{fig.beam.sparse_dataset}, where Gaussian perturbations scaled by $0.4\bar{\mathbf{z}}_{max}/\sqrt[3]{M}$ are added independently pointwise to both the strain and the stress data.


\begin{table}[h!]
    \centering
    \caption{Sparse material datasets}
    \begin{tabular}{c c c c}
    \hline
    \hline
    Sparse & Number of & Number of Data Points & Total Number\\
    Dataset & Loading Paths & per Loading Path & of Data Points ($M$)\\
    \hline
    1 & 56 & 10 & 560 \\
    2 & 98 & 10 & 980 \\
    3 & 98 & 8 & 784 \\
    \hline
    \hline
    \end{tabular}
    \label{tab.sparse_dataset}
\end{table}

\begin{figure}[!ht]
\centering
    \begin{subfigure}{0.32\textwidth}
        \centering
        \includegraphics[width=1\linewidth]{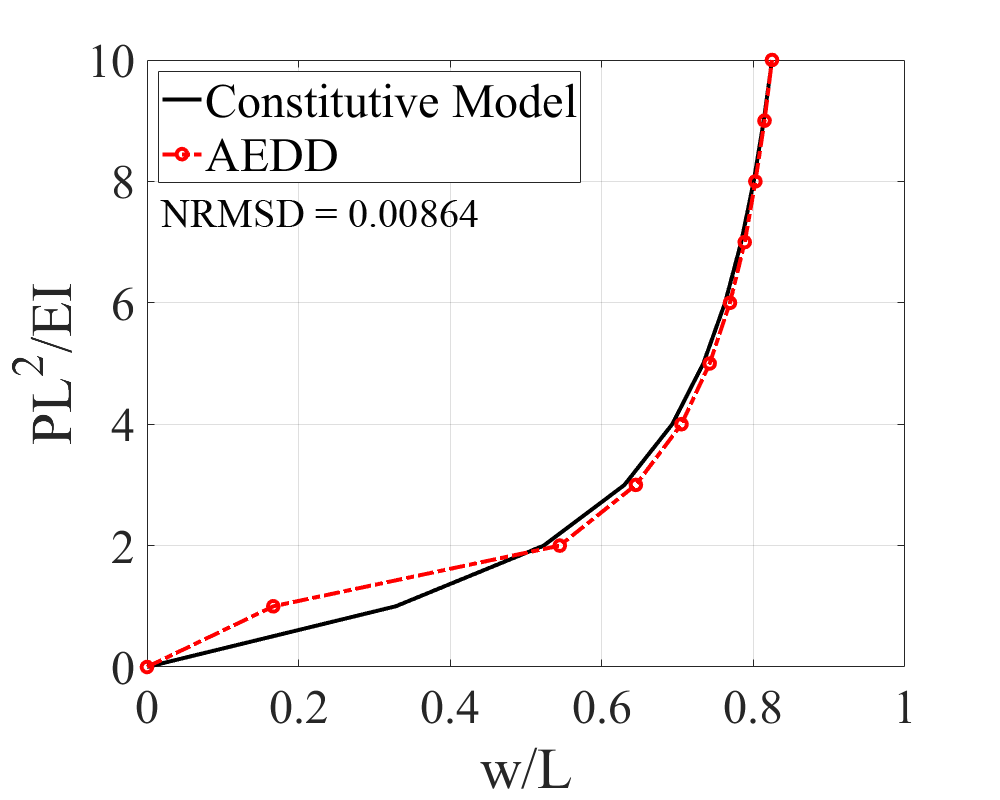}
        \caption{}
    \end{subfigure}
    \begin{subfigure}{0.32\textwidth}
        \centering
        \includegraphics[width=1\linewidth]{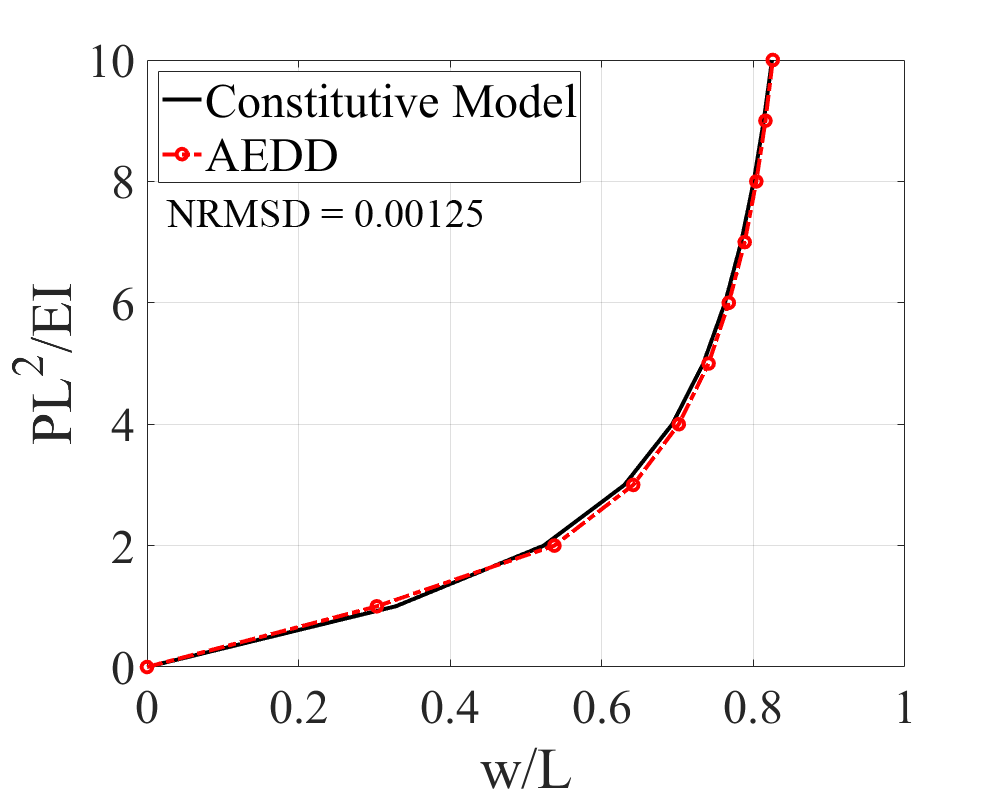}
        \caption{}
    \end{subfigure}
    \begin{subfigure}{0.32\textwidth}
        \centering
        \includegraphics[width=1\linewidth]{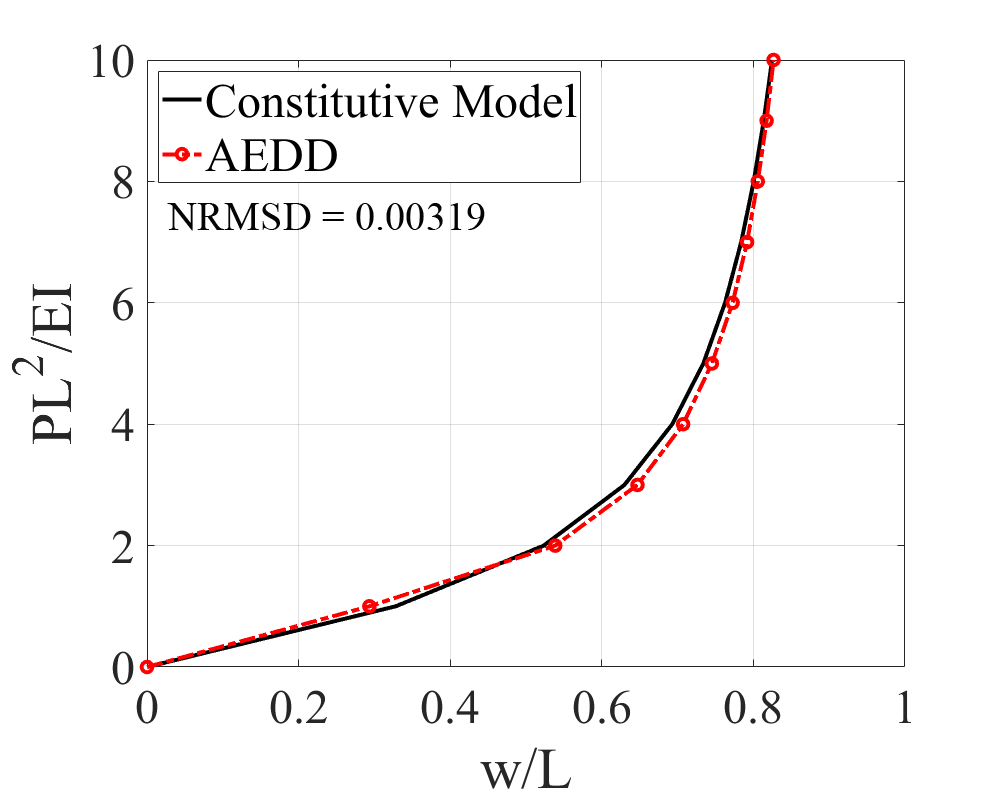}
        \caption{}
    \end{subfigure}
\caption{Comparison of constitutive model-based and AEDD solutions: (a) sparse dataset 1; (b) sparse dataset 2; (c) sparse dataset 3}
\label{fig.beam.sparse_result}
\end{figure}

An autoencoder (6-4-3) is trained using the sparse noisy datasets and used in AEDD modeling of the cantilever beam problem. The normalized tip deflection-loading responses predicted by the proposed AEDD method are compared with the constitutive model-based solutions, as shown in Fig. \ref{fig.beam.sparse_result}. The results demonstrate that the proposed AEDD method remains robust and accurate when dealing with noisy material datasets at different levels of sparsity and that the data-driven prediction accuracy improves as the data density increases.

\subsection{Biological tissue data-driven modeling}\label{sec3.2}
\begin{figure}[!ht]
\centering
    \begin{subfigure}{0.32\textwidth}
        \centering
        \includegraphics[width=1\linewidth]{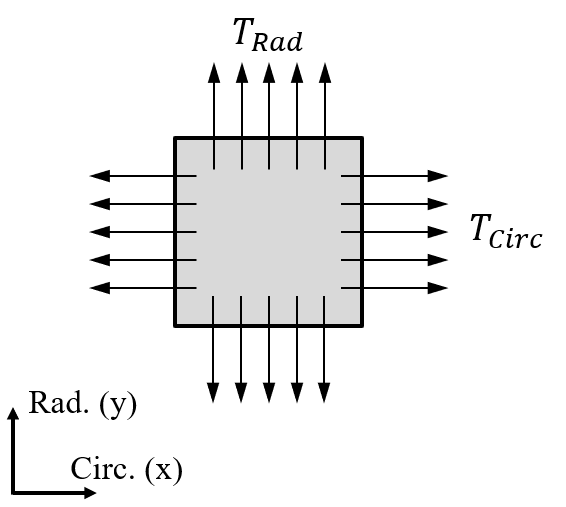}
        \caption{}
    \end{subfigure}
    \begin{subfigure}{0.25\textwidth}
        \centering
        \includegraphics[width=1\linewidth]{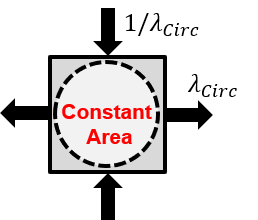}
        \caption{}
    \end{subfigure}
    \begin{subfigure}{0.41\textwidth}
        \centering
        \includegraphics[width=1\linewidth]{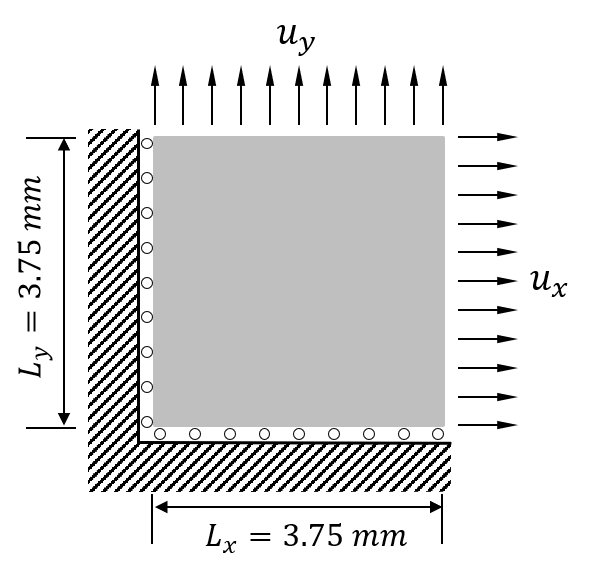}
        \caption{}
    \end{subfigure}
    \begin{subfigure}{0.49\textwidth}
        \centering
        \includegraphics[width=1\linewidth]{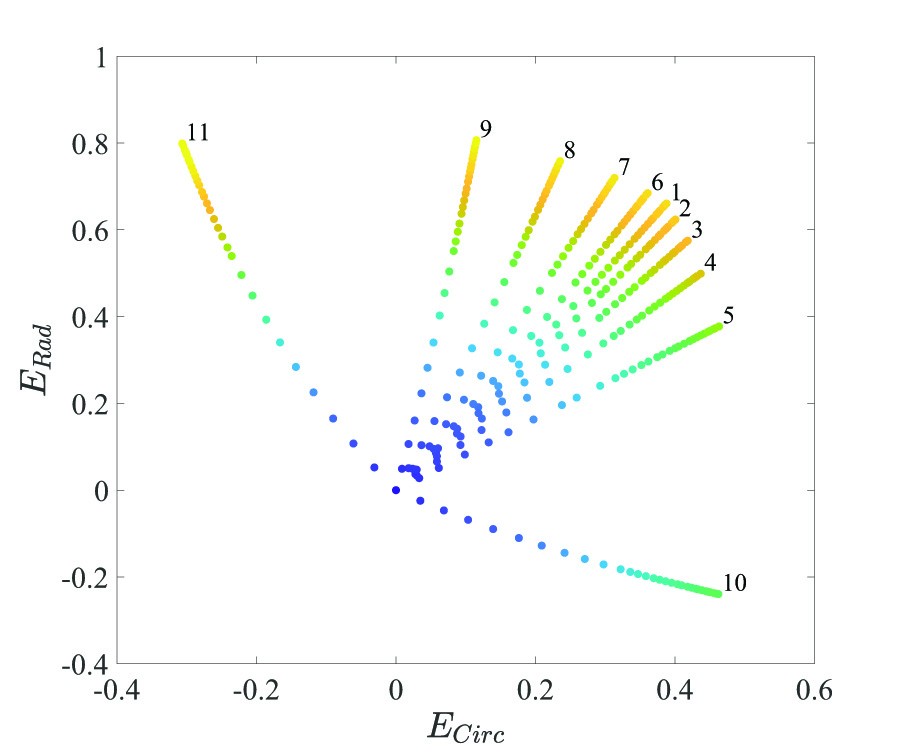}
        \caption{}
    \end{subfigure}
    \begin{subfigure}{0.49\textwidth}
        \centering
        \includegraphics[width=1\linewidth]{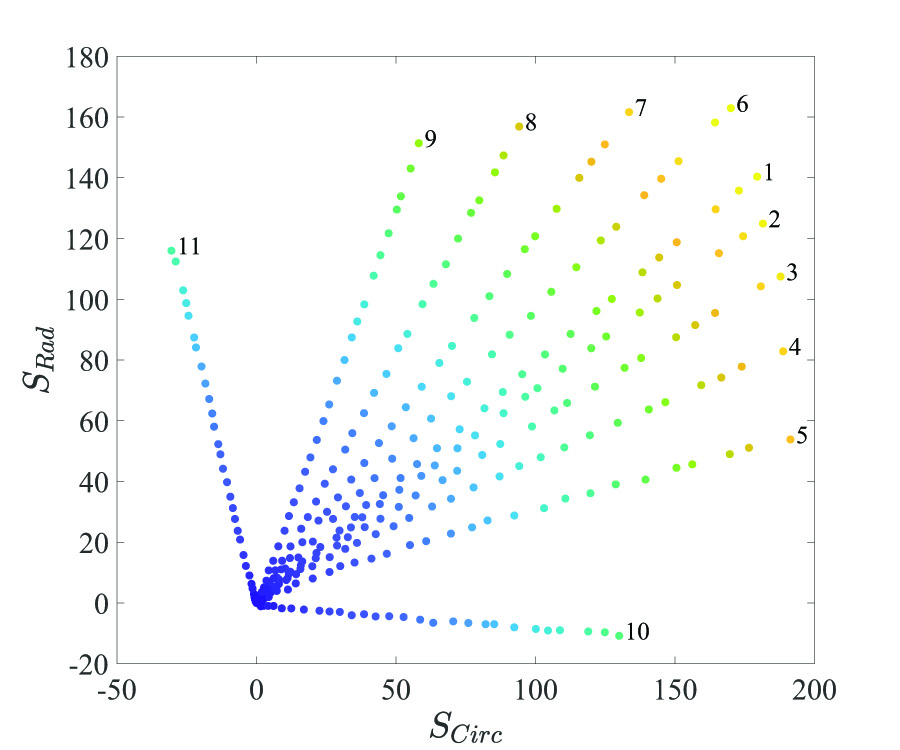}
        \caption{}
    \end{subfigure}
\caption{(a) Schematic of a mitral valve posterior leaflet (MVPL) specimen mounted on a biaxial testing system; (b) pure shear protocol 10 (\textit{x}: tension, \textit{y}: compression); (c) schematic of the model of biaxial testing in data-driven computation; (d) Green strain of all protocols; (e) 2nd-PK stress of all protocols}
\label{fig.bio.exp}
\end{figure}

The effectiveness of the proposed AEDD computational framework is examined by using the biological data from biaxial mechanical experiments of a porcine mitral valve posterior leaflet (MVPL) \cite{jett2018investigation}.
Fig. \ref{fig.bio.exp}(a) shows the schematic of a MVPL specimen with a dimension $7.5 mm \times 7.5 mm$ subjected to prescribed displacements, where the tissue's circumferential and radial directions are denoted as \textit{x} and \textit{y} axes, respectively, and the stretch ratios along these two directions are defined as $\lambda_{Circ}$ and $\lambda_{Rad}$.

A total of eleven protocols (Table \ref{tab.protocol}) includes \textit{nine biaxial tension protocols} with various tension ratios and \textit{two pure shear protocols}, as illustrated in Fig. \ref{fig.bio.exp}(b).
The normal components of the Green strain and the associated 2nd-PK stress tenors generated from the 11 biaxial mechanical testing are plotted in Fig. \ref{fig.bio.exp}(d) and Fig. \ref{fig.bio.exp}(e), respectively. It shows that the measured data points are sparse in the stress-strain phase space.
It is noted that in the mechanical testing the direct measurements are the applied membrane tensions, $T_{Rad}$ and $T_{Circ}$, and the displacements are estimated by digital image correlation techniques. Thus, the measured Green strain and 2nd-PK stress data are based on homogeneous deformation assumption in the test specimen. More details about the tissue strain and stress calculations as well as the experimental setting can be found in \cite{jett2018investigation,HeJBM2020}.


\begin{table}[h!]
    \centering
    \caption{Eleven biaxial mechanical testing protocols of a representative MVPL specimen and the corresponding measured displacements used in data-driven computations}
    \begin{tabular}{c c c c c c}
    \hline
    \hline
    Protocol ID & Protocol & $\lambda_{Circ}$ & $\lambda_{Rad}$ & $u_{Circ}$ (mm) & $u_{Rad}$ (mm)\\
    \hline
    1 & Biaxial Tension $T_{Circ}:T_{Rad}$=1:1 & 1.333 & 1.525 & 2.498 & 3.938 \\
    2 & Biaxial Tension $T_{Circ}:T_{Rad}$=1:0.8 & 1.342 & 1.499 & 2.564 & 3.744 \\
    3 & Biaxial Tension $T_{Circ}:T_{Rad}$=1:0.6 & 1.355 & 1.466 & 2.662 & 3.498 \\
    4 & Biaxial Tension $T_{Circ}:T_{Rad}$=1:0.4 & 1.369 & 1.415 & 2.770 & 3.110 \\
    5 & Biaxial Tension $T_{Circ}:T_{Rad}$=1:0.2 & 1.388 & 1.326 & 2.913 & 2.442 \\
    6 & Biaxial Tension $T_{Circ}:T_{Rad}$=0.8:1 & 1.313 & 1.541 & 2.344 & 4.055 \\
    7 & Biaxial Tension $T_{Circ}:T_{Rad}$=0.6:1 & 1.275 & 1.562 & 2.064 & 4.215 \\
    8 & Biaxial Tension $T_{Circ}:T_{Rad}$=0.4:1 & 1.213 & 1.588 & 1.596 & 4.409 \\
    9 & Biaxial Tension $T_{Circ}:T_{Rad}$=0.2:1 & 1.109 & 1.618 & 0.820 & 4.635 \\
    10 & Pure Shear in x & 1.387 & 0.721 & 2.903 & -2.093 \\
    11 & Pure Shear in y & 0.620 & 1.612 & -2.847 & 4.590 \\
    \hline
    \hline
    \end{tabular}
    \label{tab.protocol}
\end{table}

Five study cases are considered to evaluate the performance of the proposed AEDD framework, which is compared with that of the LCDD method \cite{he2019physics,HeJBM2020}. In these tests (Case 1--5), the experimental data (see Fig. \ref{fig.bio.exp}(d-e)) associated with the selected biaxial testing protocols, called \textit{training protocols}, are used for constructing material dataset $\mathbb{E}$, and different data-driven modeling approaches with the constructed material dataset are tested on other protocols (called \textit{testing protocols}) to assess their performance against the experimental results. The training and testing protocols of the Five study cases are described as below:
\begin{itemize}
    \item \textbf{Case 1}: \textit{Training Protocols}: 1, 3, 4, 7, and 8; \textit{Testing Protocols}: 2 and 5, used to investigate AEDD's performance in interpolative and extrapolative predictions.
    \item \textbf{Case 2}: \textit{Training Protocols}: 1, 3, 4, 7, 8, 10, and 11; \textit{Testing Protocols}: 2 and 5, used to investigate AEDD's performance in interpolative and extrapolative predictions.
    \item \textbf{Case 3}: \textit{Training Protocols}: 1, 2, 6, 10, and 11; \textit{Testing Protocols}: 3 and 4, used to investigate AEDD's performance in extrapolative prediction.
    \item \textbf{Case 4}: \textit{Training Protocols}: 2, 5, 7, and 8, which are asymmetrically distributed; \textit{Testing Protocols}: 1, 3 and 4, used to investigate AEDD's performance in intrapolative prediction.
    \item \textbf{Case 5}: \textit{Training Protocols}: 1 -- 9; \textit{Testing Protocols}: 10 and 11.
\end{itemize}
For the first three cases, the protocols used for training are symmetrically distributed, while the training protocols are asymmetrically distributed for the last case.

For AEDD, autoencoders are first trained offline using the training protocols and then employed in the local step of the data-driven solvers (Section \ref{AEDDsolvers}) of AEDD during the online computation. In the following study, Solver II (Section \ref{sec.AEDD_locII}) is employed and the number of nearest neighbors in locally convex reconstruction of the data-driven solver is set as 6.
A diagonal matrix is used as the weight matrix $\hat{\mathbb{C}}$ with each diagonal component being the ratio of the standard deviation of the associated component of the stress data to that of the strain data. This is similar to the normalization technique used in deep learning that applies the standard deviation of each input unit to inversely scale the input data \cite{goodfellow2016deep}.

The prediction of data-driven methods on testing protocols that are not included in the training dataset are compared with the corresponding experimental data.
The \textit{NRMSD} (Eq. (\ref{eq.3.1.2})) normalized with respect to the maximum stress of the experimental data is employed to assess the prediction performance of the methods. In the data-driven modeling, considering the symmetric geometry of the tissue specimen and the symmetric loading conditions, the upper right quarter of the sample is modelled with symmetric boundary conditions, as shown in Fig. \ref{fig.bio.exp}(c), and the prescribed displacements are applied to the top and the right boundaries.

\subsubsection{Case 1}\label{sec3.2.1}
\begin{figure}[!ht]
\centering
    \begin{subfigure}{0.32\textwidth}
        \centering
        \includegraphics[width=1\linewidth]{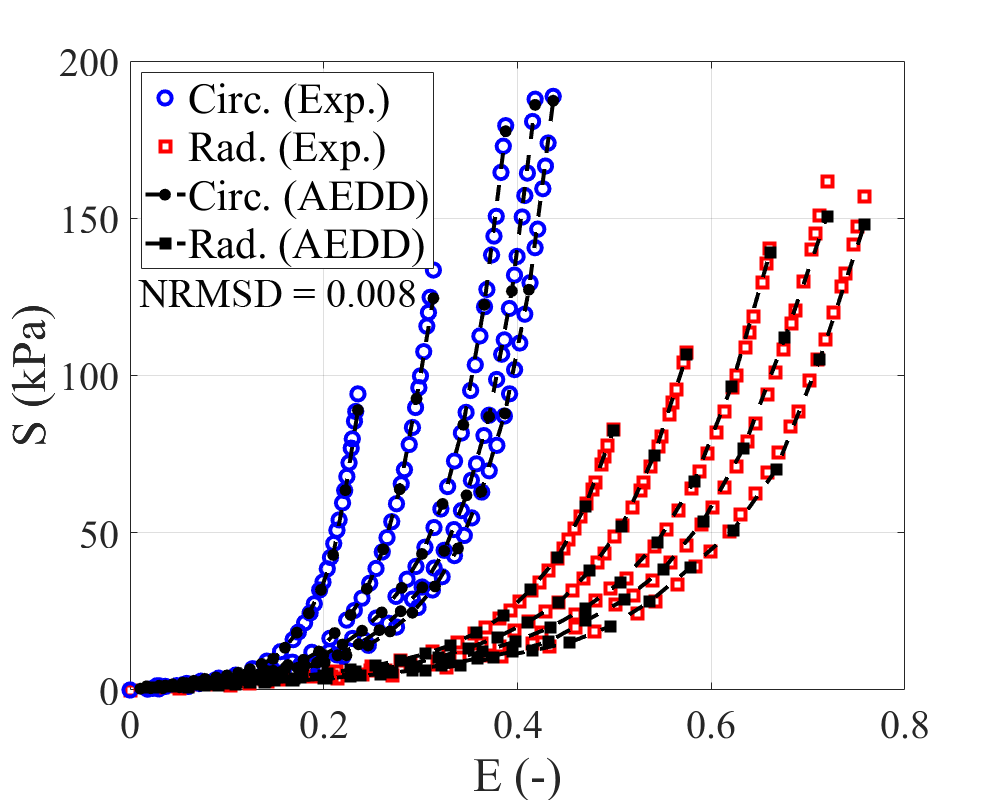}
        \caption{Protocols 1, 3, 4, 7, 8}
    \end{subfigure}
    \begin{subfigure}{0.32\textwidth}
        \centering
        \includegraphics[width=1\linewidth]{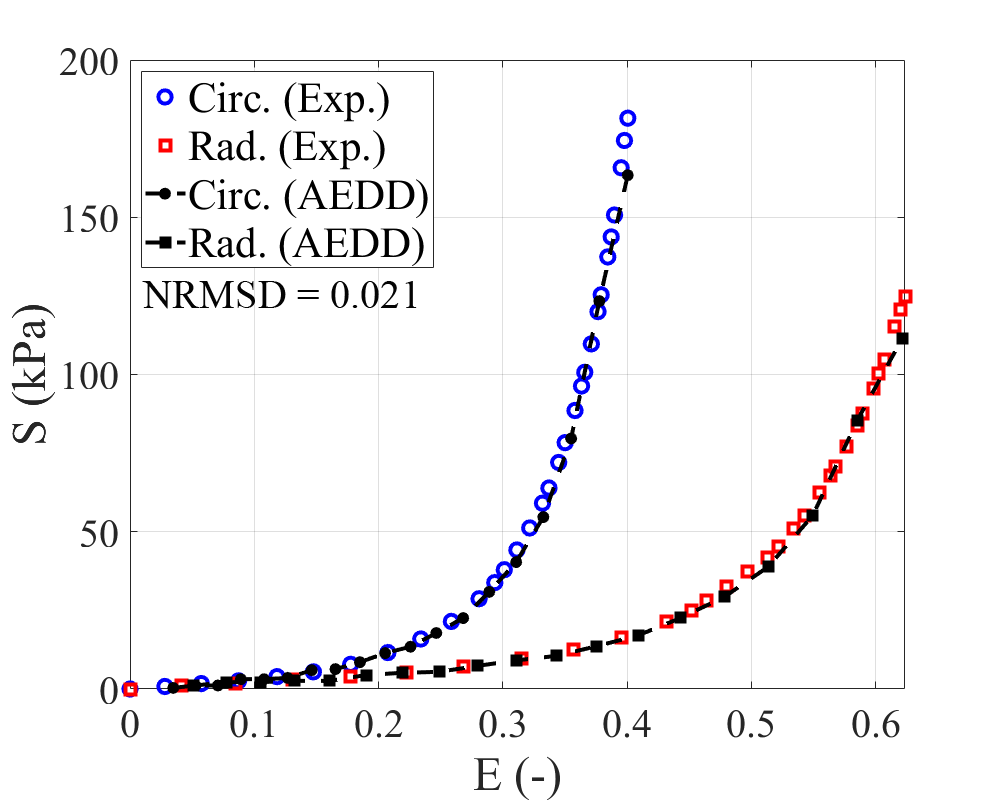}
        \caption{Protocol 2}
    \end{subfigure}
    \begin{subfigure}{0.32\textwidth}
        \centering
        \includegraphics[width=1\linewidth]{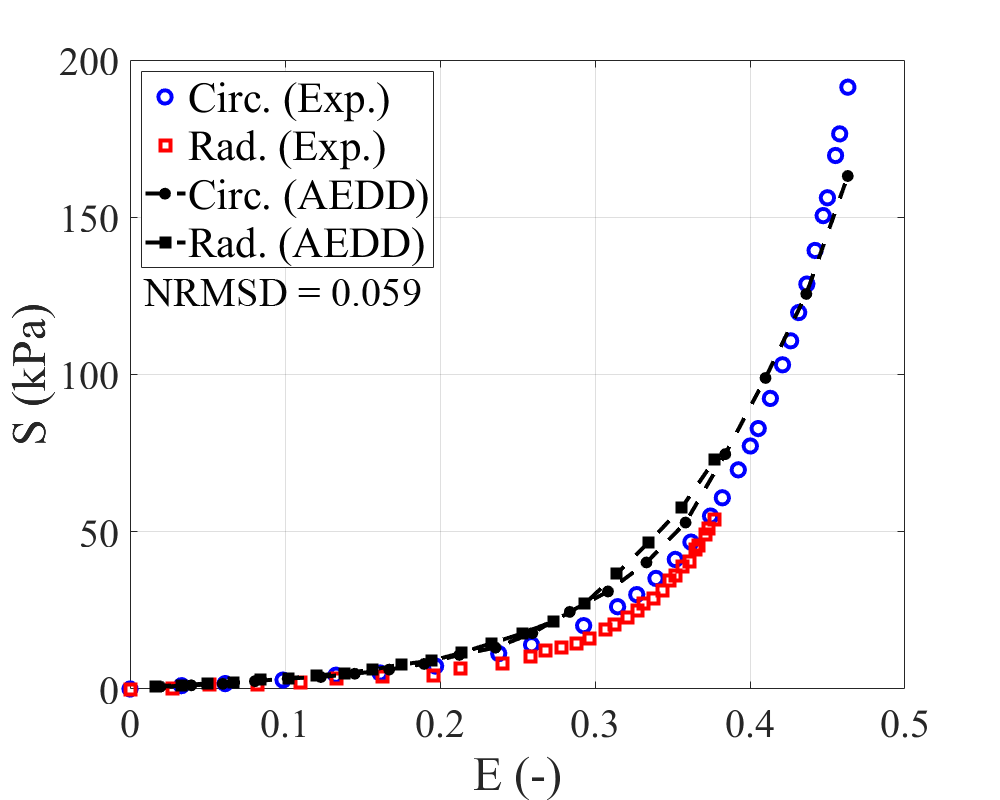}
        \caption{Protocol 5}
    \end{subfigure}
    \begin{subfigure}{0.32\textwidth}
        \centering
        \includegraphics[width=1\linewidth]{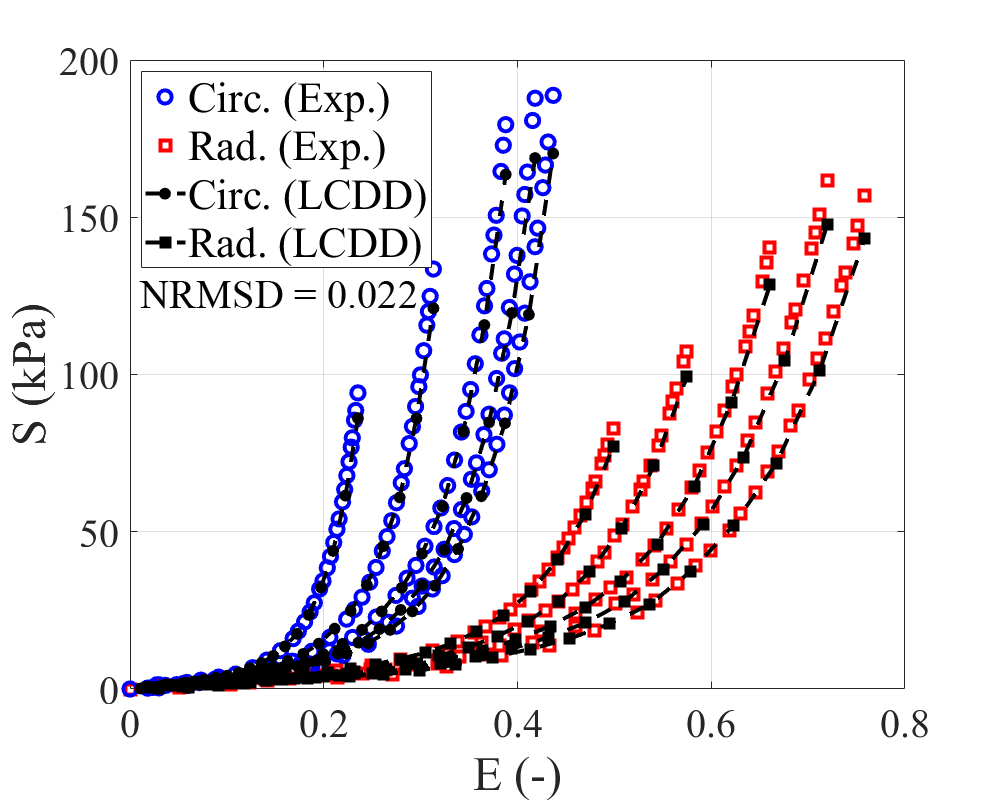}
        \caption{Protocols 1, 3, 4, 7, 8}
    \end{subfigure}
    \begin{subfigure}{0.32\textwidth}
        \centering
        \includegraphics[width=1\linewidth]{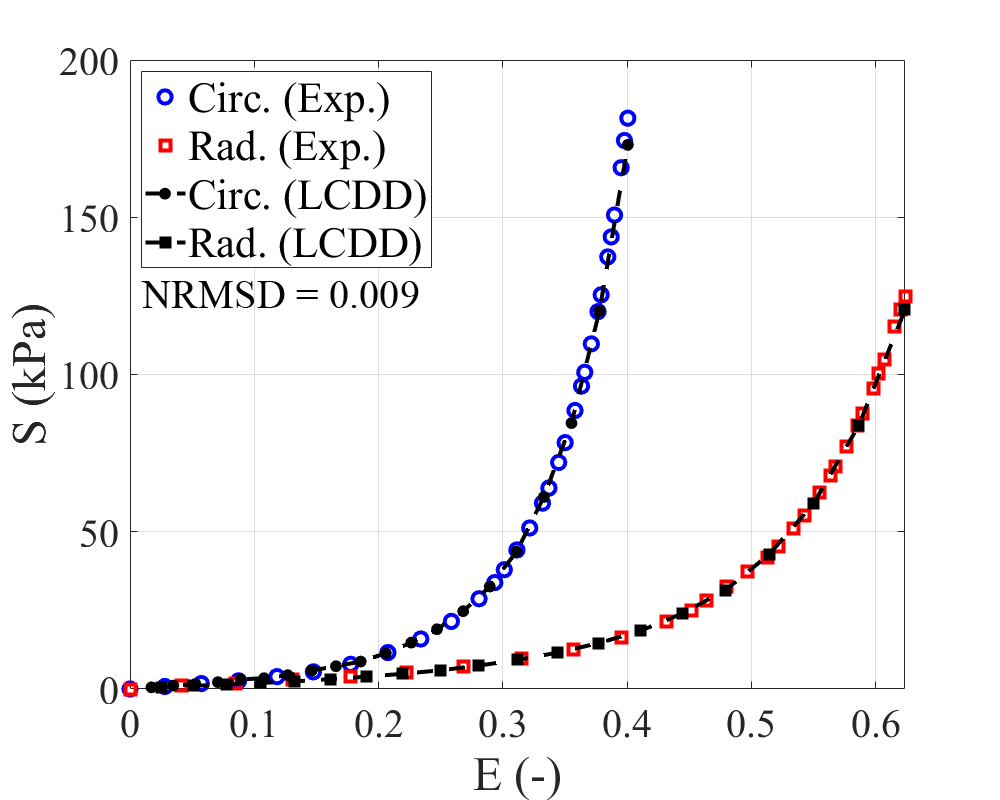}
        \caption{Protocol 2}
    \end{subfigure}
    \begin{subfigure}{0.32\textwidth}
        \centering
        \includegraphics[width=1\linewidth]{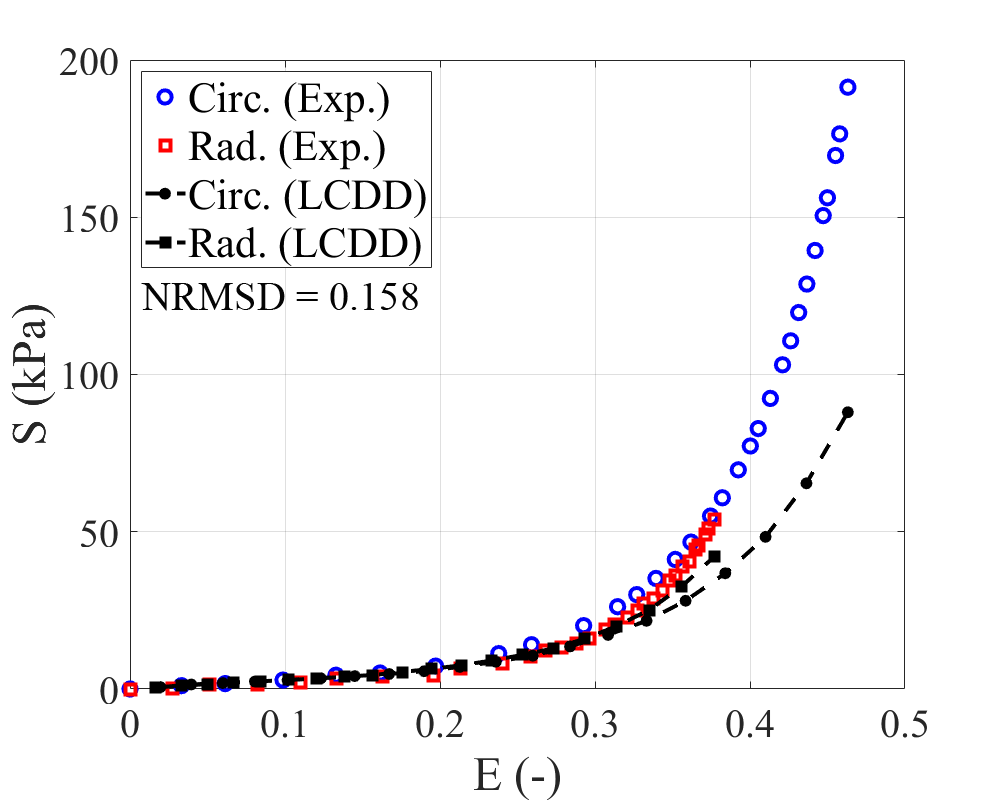}
        \caption{Protocol 5}
    \end{subfigure}
\caption{Comparison of interpolative (Protocol 2) and extrapolative (Protocol 5) predictability: (a) AEDD prediction on training Protocols 1, 3, 4, 7, 8; (b) AEDD prediction on Protocol 2; (c) AEDD prediction on Protocol 5; (d) LCDD prediction on training Protocols 1, 3, 4, 7, 8; (e) LCDD prediction on Protocol 2;  (f) LCDD prediction on Protocol 5. Protocols 1, 3, 4, 7, and 8 are used to train the autoencoder applied in AEDD}
\label{fig.bio.case1}
\end{figure}

We first examine the data fitting capability whereby the data-driven methods are tested on the training protocols 1, 3, 4, 7 and 8, as shown in Fig. \ref{fig.bio.case1}(a) and Fig. \ref{fig.bio.case1}(d). It shows that both AEDD (\textit{NRMSD}$_{AEDD}$=0.008) and LCDD (\textit{NRMSD}$_{LCDD}$=0.022) provide satisfactory fitting results, but AEDD yields a slightly higher accuracy.
Since the strains and stresses of testing protocols 2 and 5 lie inside and outside the domain covered by the data of the training protocols, respectively, as shown in Fig. \ref{fig.bio.exp}(d-e), the AEDD predictions on the testing protocols 2 and 5 are interpolative and extrapolative predictions, respectively. For the interpolative prediction test on Protocol 2, the results of these two approaches also agree well with the experimental data, as shown in Fig. \ref{fig.bio.case1}(b) and (e). The \textit{NRMSD} errors indicate that LCDD achieves a higher accuracy, i.e. $\textit{NRMSD}_{LCDD}=0.009 < \textit{NRMSD}_{AEDD}=0.021$.
However, its extrapolative prediction on Protocol 5 is worse than that from AEDD ($\textit{NRMSD}_{LCDD}=0.158 > \textit{NRMSD}_{AEDD}=0.059$), see Fig. \ref{fig.bio.case1}(c) and (f).
The results demonstrate better extrapolative generalization ability of AEDD. It could be attributed to the underlying low-dimensional global material manifold learned by the autoencoders. Specifically, AEDD performs local neighbor searching and locally convex reconstruction of optimal material state based on geometric distance information in the low-dimensional global embedding space, which contains the underlying manifold structure of the material data and contributes to a higher solution accuracy and better generalization performance. In contrast, LCDD performs local neighbor searching and locally convex reconstruction purely from the existing material data points without any generalization, leading to lower extrapolative generalization ability.  

\begin{figure}[!ht]
\centering
    \begin{subfigure}{0.49\textwidth}
        \centering
        \includegraphics[width=1\linewidth]{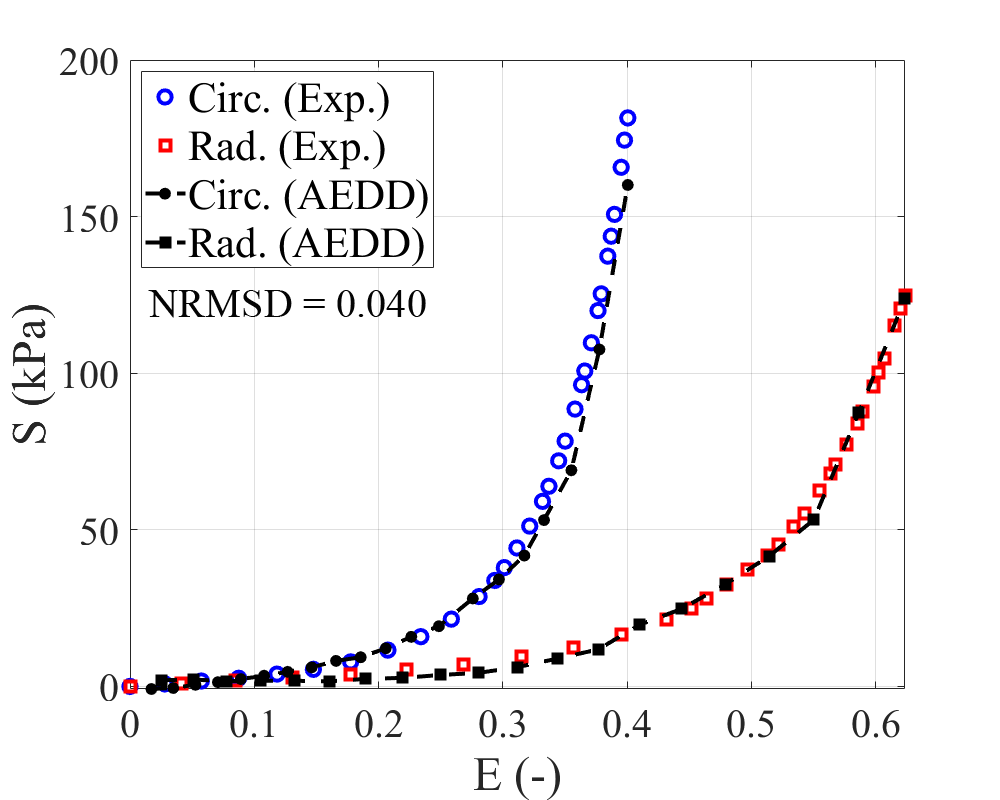}
        \caption{Protocol 2}
    \end{subfigure}
    \begin{subfigure}{0.49\textwidth}
        \centering
        \includegraphics[width=1\linewidth]{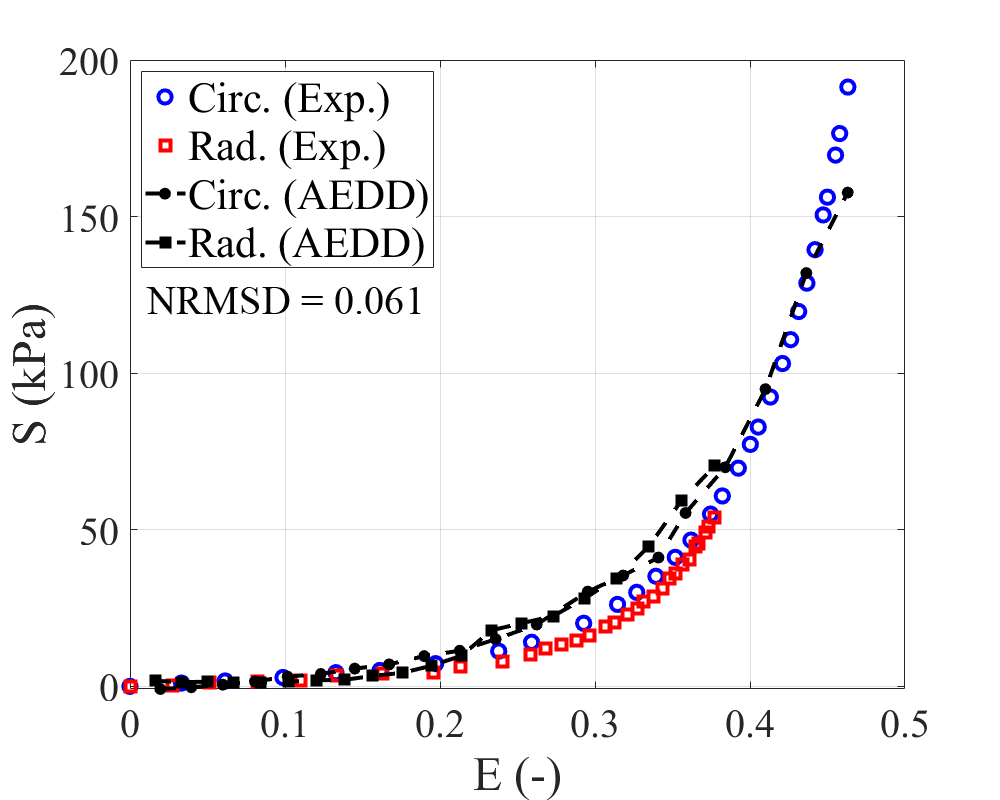}
        \caption{Protocol 5}
    \end{subfigure}
\caption{Data-driven prediction by AEDD with Solver I on (a) Protocol 2 and (b) Protocol 5. Protocols 1, 3, 4, 7, and 8 are used to train the autoencoder}
\label{fig.bio.case1_I}
\end{figure}

Another proposed AEDD method with Solver I (Section \ref{sec.AEDD_locI}) using the same training protocols as material dataset are also investigated, as shown in Fig. \ref{fig.bio.case1_I}. As expected, compared to the results obtained by using Solver II, see Fig. \ref{fig.bio.case1}(b) and (c), the prediction capability by Solver I decreases on both testing protocols. Especially in the interpolative prediction test Protocol 2, the \textit{NRMSD} error increases to $0.04$ from $0.021$. We attribute the larger errors with Solver I to the employment of decoders in constructing the optimal material state. Since we have demonstrated that the AEDD approach with Solver II provides better data-driven prediction results, we only consider this approach in the following study.

\subsubsection{Case 2}\label{sec3.2.2}
In this case study, the objective is to verify how the incorporation of material data of different deformation modes affects the interpolative and extrapolative predictability in the proposed data-driven modeling. Two pure shear protocols are introduced in the training material dataset in addition to the biaxial tension protocols used in Case 1. The two pure shear protocols (10 and 11) in the training dataset exhibit different material behaviors from the remaining biaxial tension protocols (1, 3, 4, 7, and 8). The AEDD predictions on the testing protocols 2 and 5 are interpolative and extrapolative predictions, respectively.

\begin{figure}[!ht]
\centering
    \begin{subfigure}{0.32\textwidth}
        \centering
        \includegraphics[width=1\linewidth]{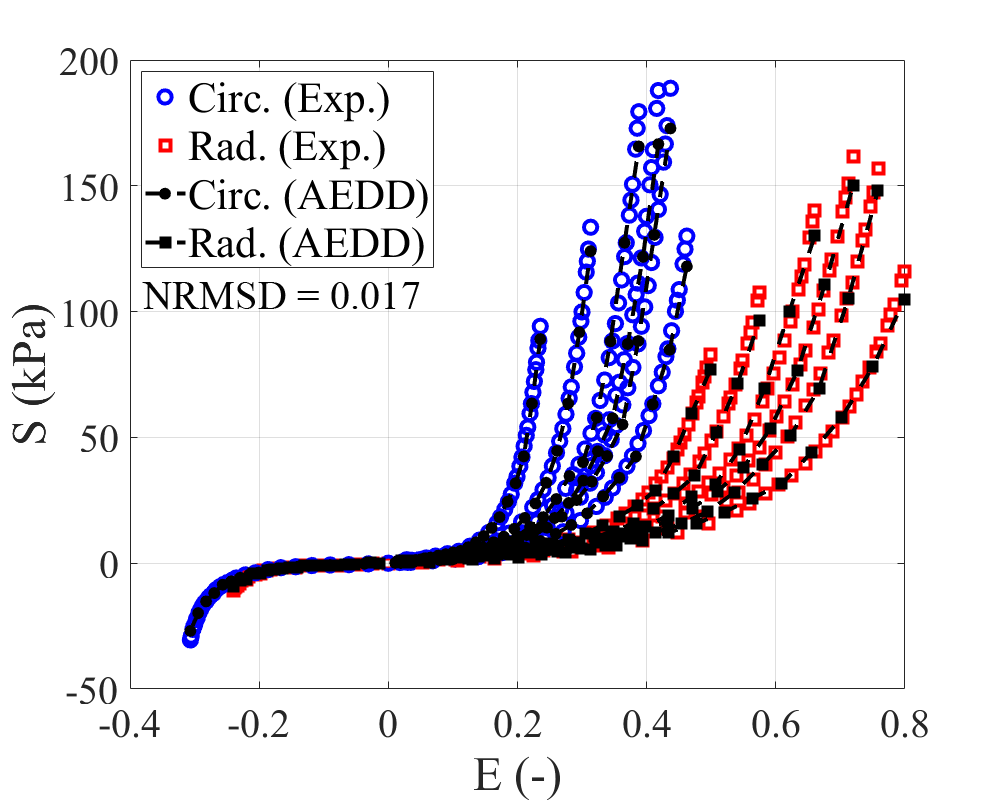}
        \caption{Protocols 1, 3, 4, 7, 8, 10, 11}
    \end{subfigure}
    \begin{subfigure}{0.32\textwidth}
        \centering
        \includegraphics[width=1\linewidth]{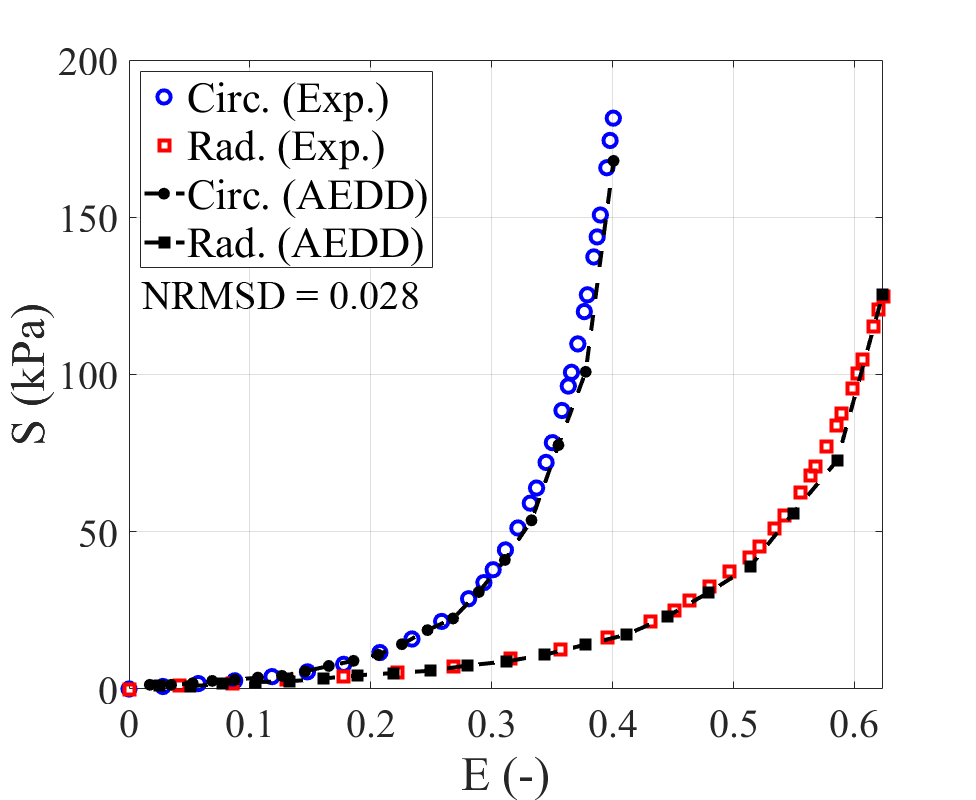}
        \caption{Protocol 2}
    \end{subfigure}
    \begin{subfigure}{0.32\textwidth}
        \centering
        \includegraphics[width=1\linewidth]{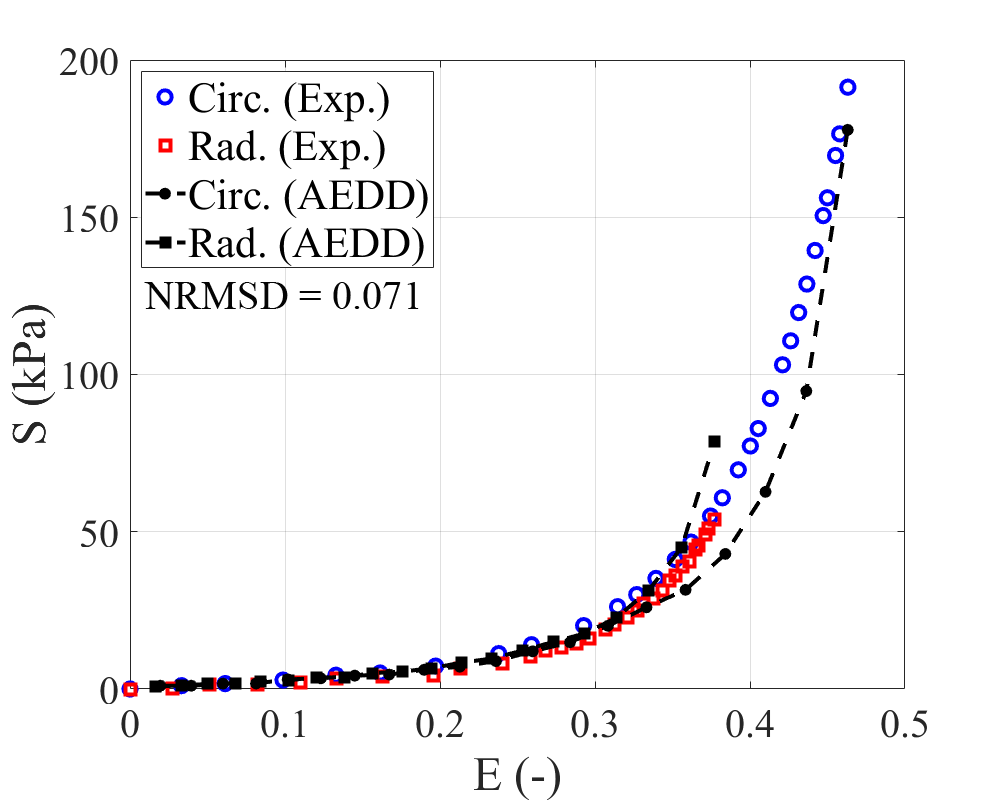}
        \caption{Protocol 5}
    \end{subfigure}
    \begin{subfigure}{0.32\textwidth}
        \centering
        \includegraphics[width=1\linewidth]{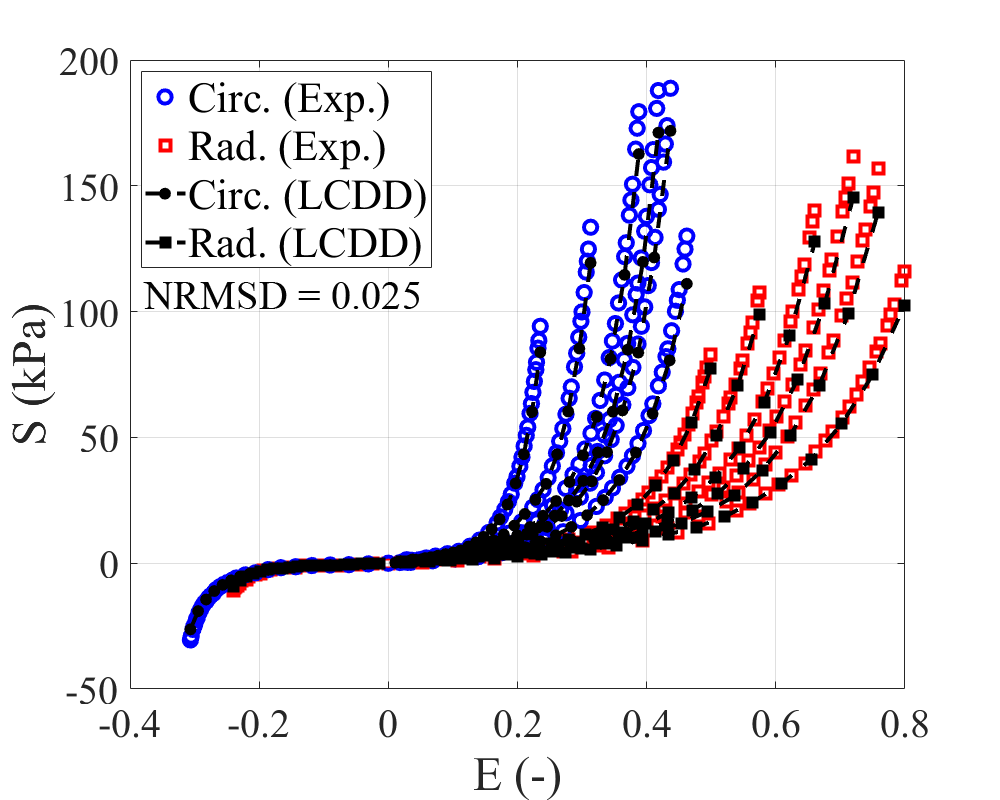}
        \caption{Protocols 1, 3, 4, 7, 8, 10, 11}
    \end{subfigure}
    \begin{subfigure}{0.32\textwidth}
        \centering
        \includegraphics[width=1\linewidth]{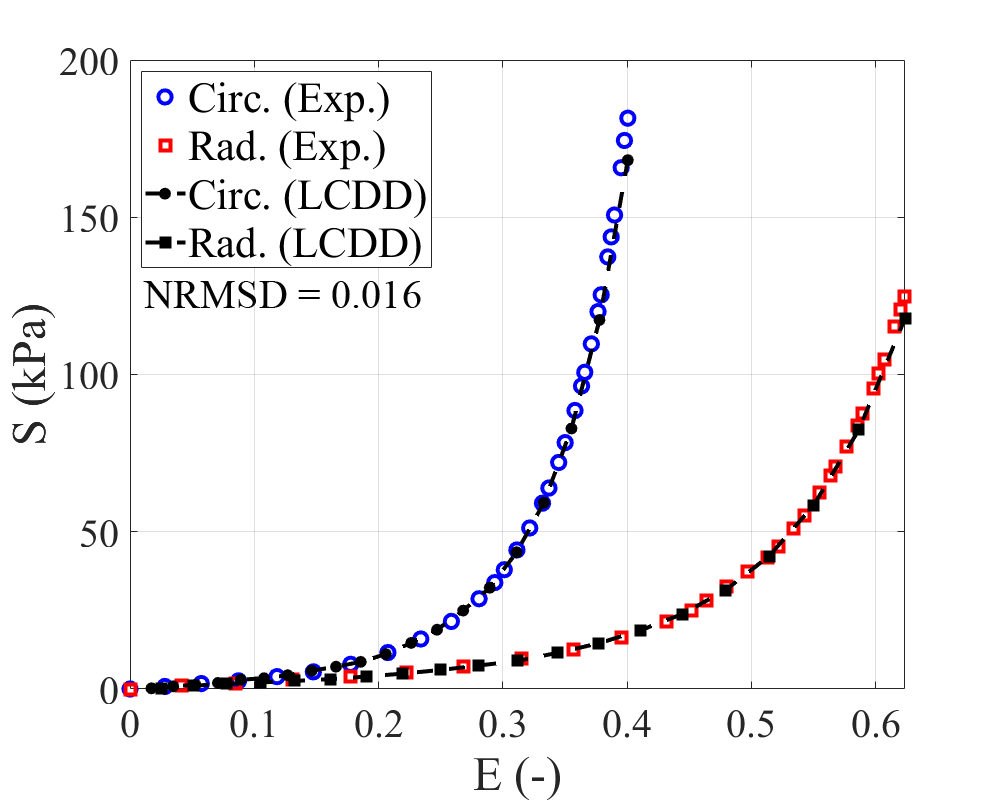}
        \caption{Protocol 2}
    \end{subfigure}
    \begin{subfigure}{0.32\textwidth}
        \centering
        \includegraphics[width=1\linewidth]{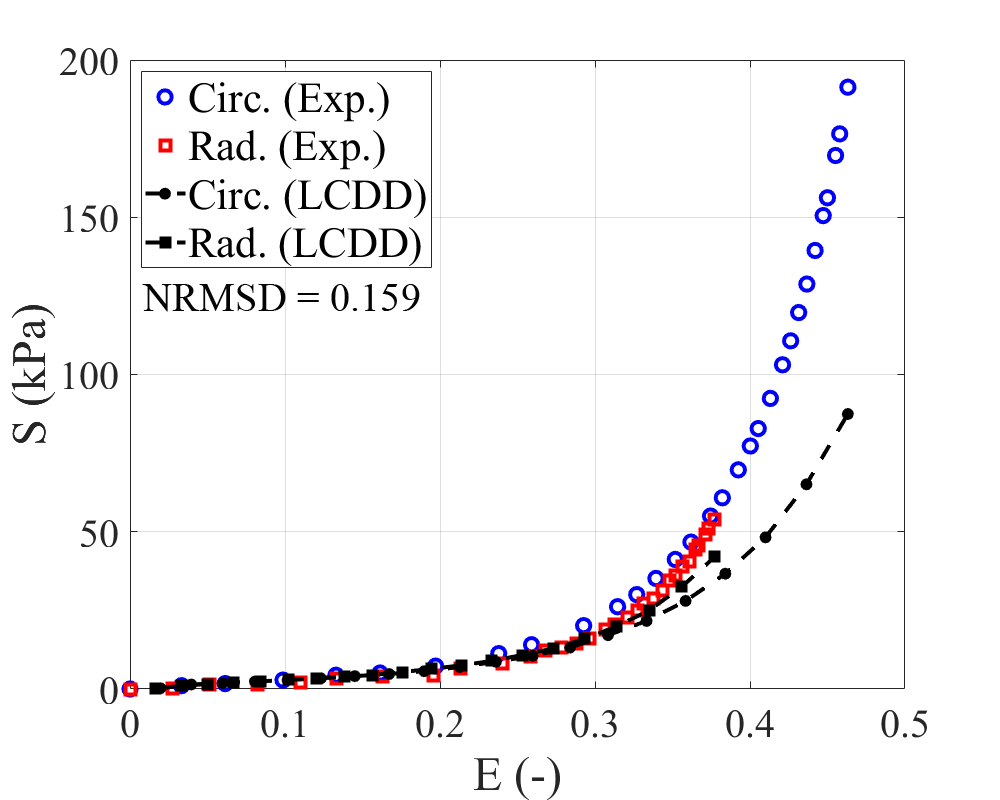}
        \caption{Protocol 5}
    \end{subfigure}
\caption{Comparison of interpolative (Protocol 2) and extrapolative (Protocol 5) predictability: (a) AEDD prediction on training Protocols 1, 3, 4, 7, 8, 10, 11; (b) AEDD prediction on Protocol 2;  (c) AEDD prediction on Protocol 5; (d) LCDD prediction on training Protocols 1, 3, 4, 7, 8, 10, 11; (e) LCDD prediction on Protocol 2;  (f) LCDD prediction on Protocol 5. Protocols 1, 3, 4, 7, 8, 10, and 11 are used to train the autoencoder applied in AEDD}
\label{fig.bio.case2}
\end{figure}

As can be seen from Fig \ref{fig.bio.case2}(a) and (d), both LCDD and AEDD maintain good fitting performance for all the biaxial tension and pure shear training protocols. 
They also perform well for the testing Protocol 2 (Fig. \ref{fig.bio.case2}(b) and (e)) with almost the same accuracy in Case 1.
This is a desirable property in data-driven methods.
AEDD again yields higher accuracy than LCDD in the extrapolative test (Protocol 5), as evidenced by the smaller NRMSD value 0.071 in the AEDD prediction over 0.159 in the LCDD prediction. This further demonstrates the enhanced extrapolative generalization in the proposed autoencoder-based approach.
Moreover, compared with Case 1, Fig. \ref{fig.bio.case2}(c) shows that AEDD with the material data from the pure shear protocols improves the prediction for strain $E < 0.35$ but results in slightly more discrepancies in the high strain range.


\subsubsection{Case 3}\label{sec3.2.3}
\begin{figure}[!ht]
\centering
    \begin{subfigure}{0.32\textwidth}
        \centering
        \includegraphics[width=1\linewidth]{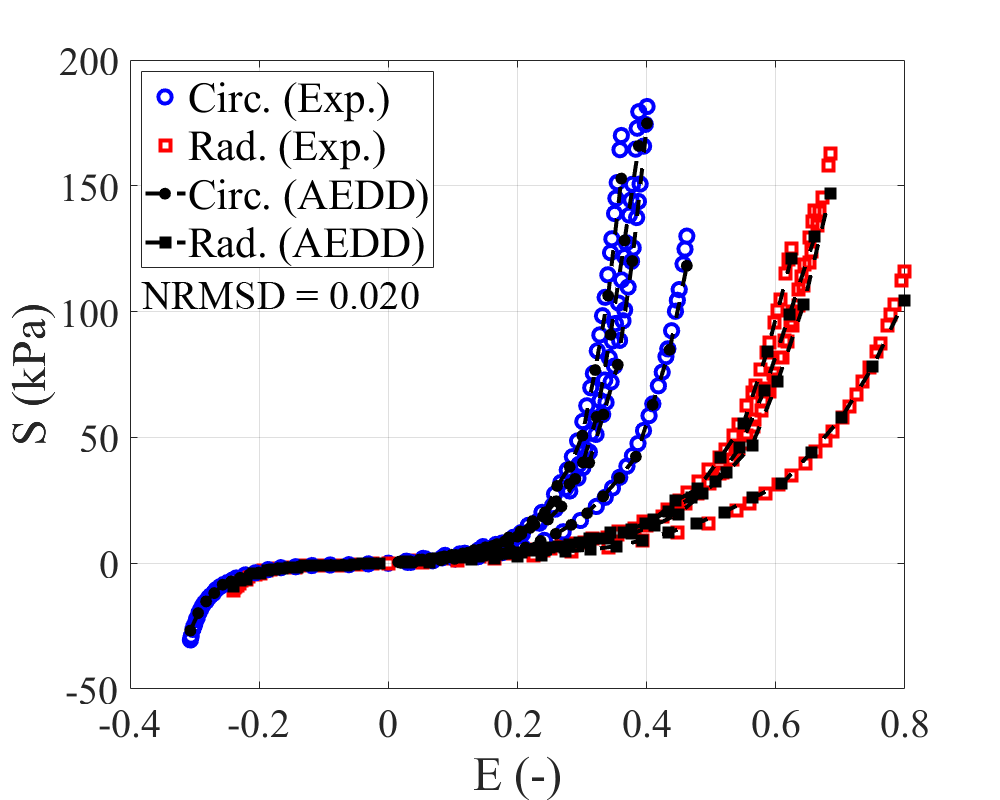}
        \caption{Protocols 1, 2, 6, 10, 11}
    \end{subfigure}
    \begin{subfigure}{0.32\textwidth}
        \centering
        \includegraphics[width=1\linewidth]{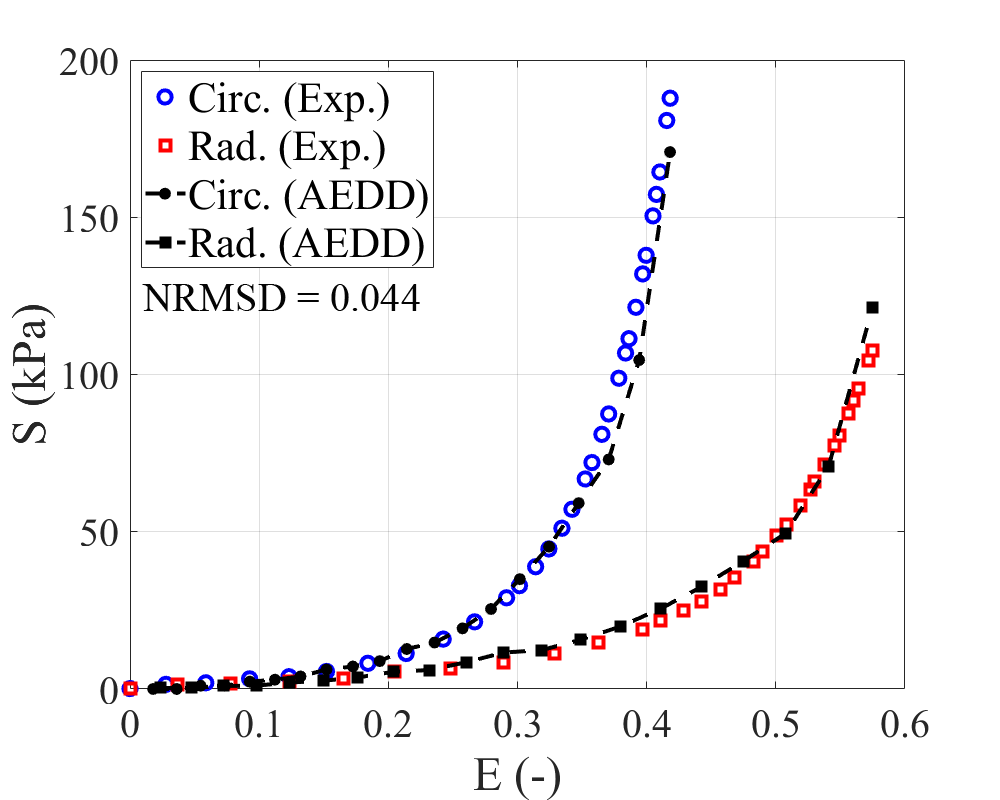}
        \caption{Protocol 3}
    \end{subfigure}
    \begin{subfigure}{0.32\textwidth}
        \centering
        \includegraphics[width=1\linewidth]{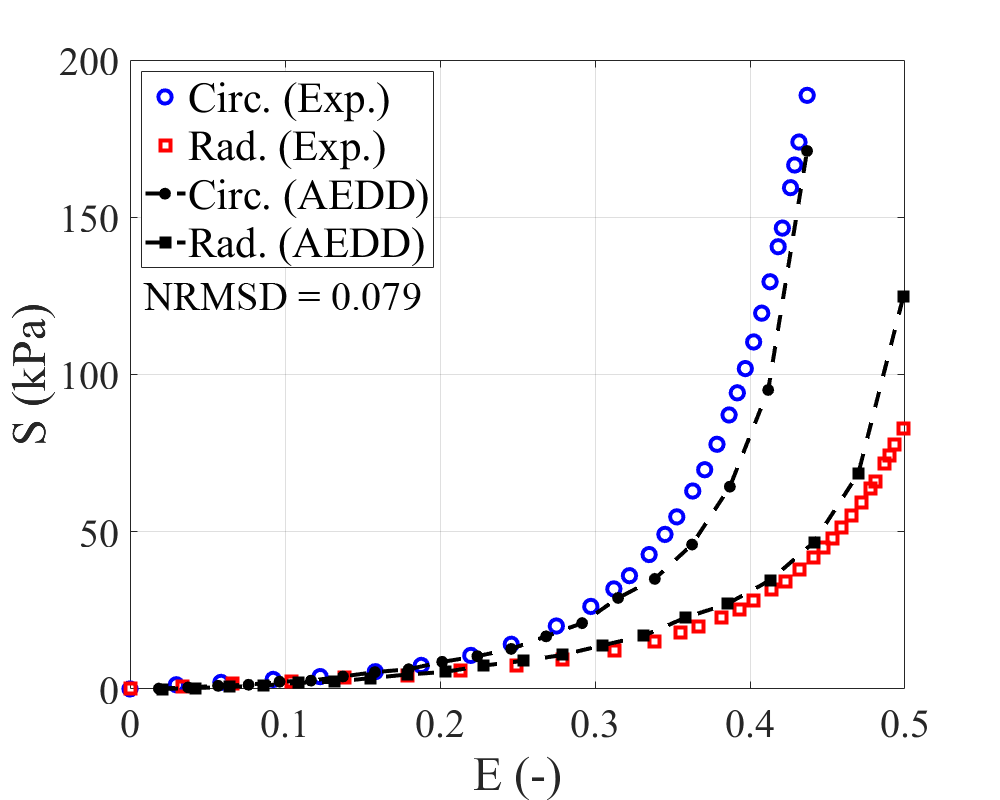}
        \caption{Protocol 4}
    \end{subfigure}
    \begin{subfigure}{0.32\textwidth}
        \centering
        \includegraphics[width=1\linewidth]{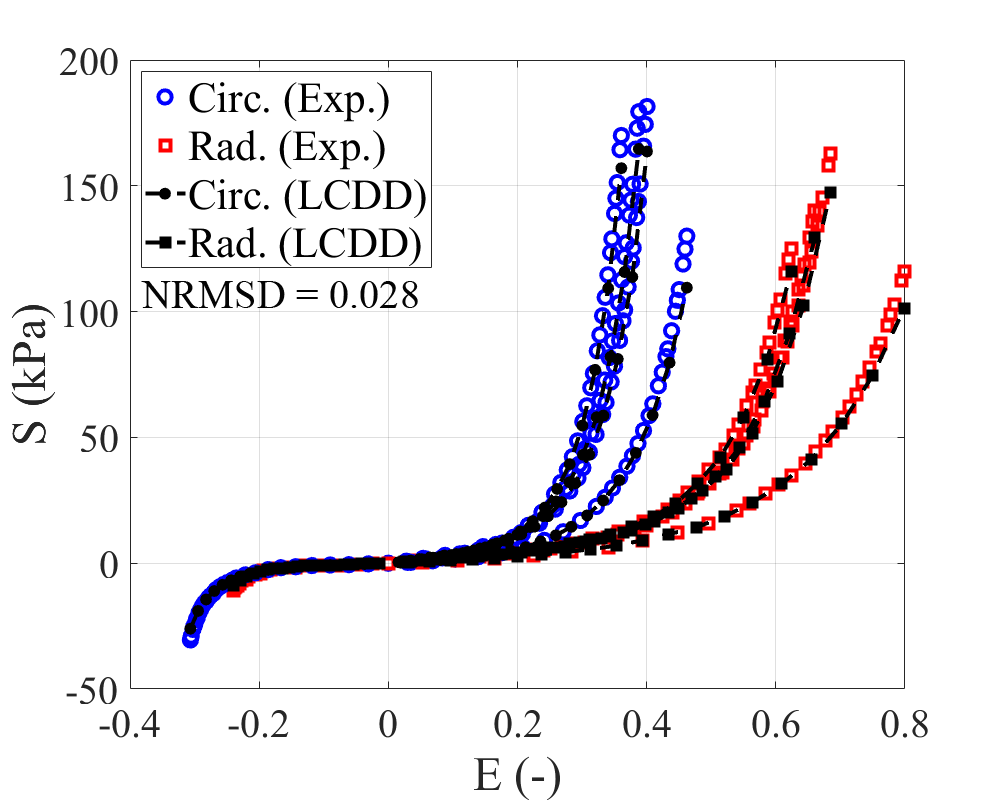}
        \caption{Protocols 1, 2, 6, 10, 11}
    \end{subfigure}
    \begin{subfigure}{0.32\textwidth}
        \centering
        \includegraphics[width=1\linewidth]{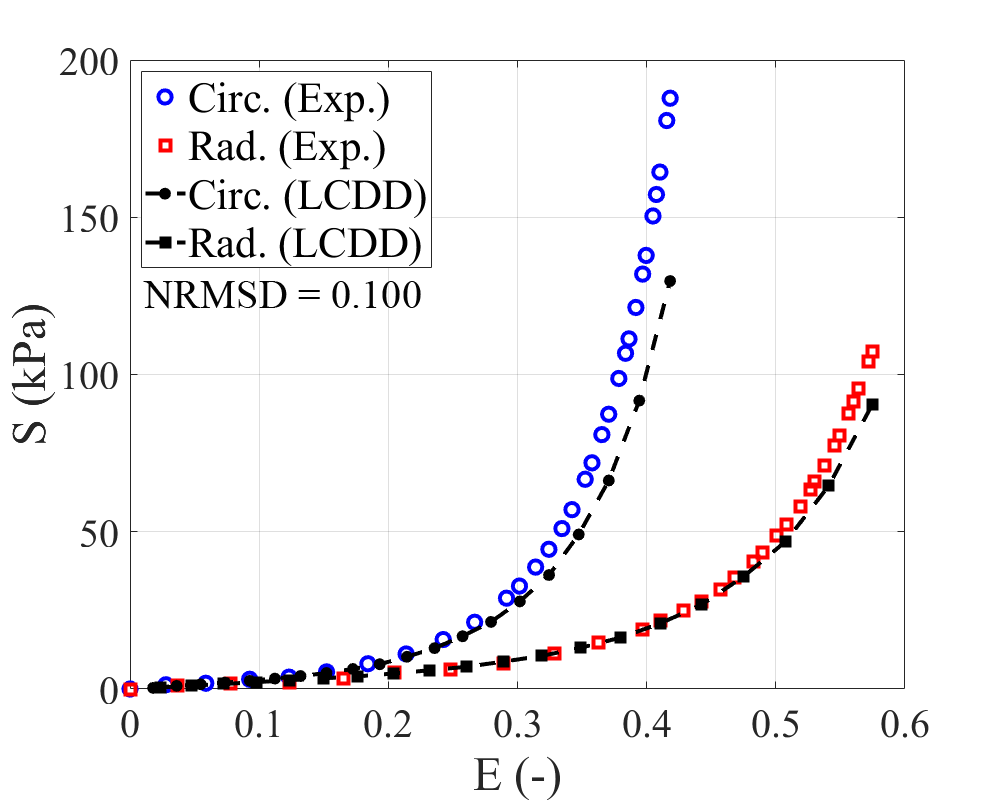}
        \caption{Protocol 3}
    \end{subfigure}
    \begin{subfigure}{0.32\textwidth}
        \centering
        \includegraphics[width=1\linewidth]{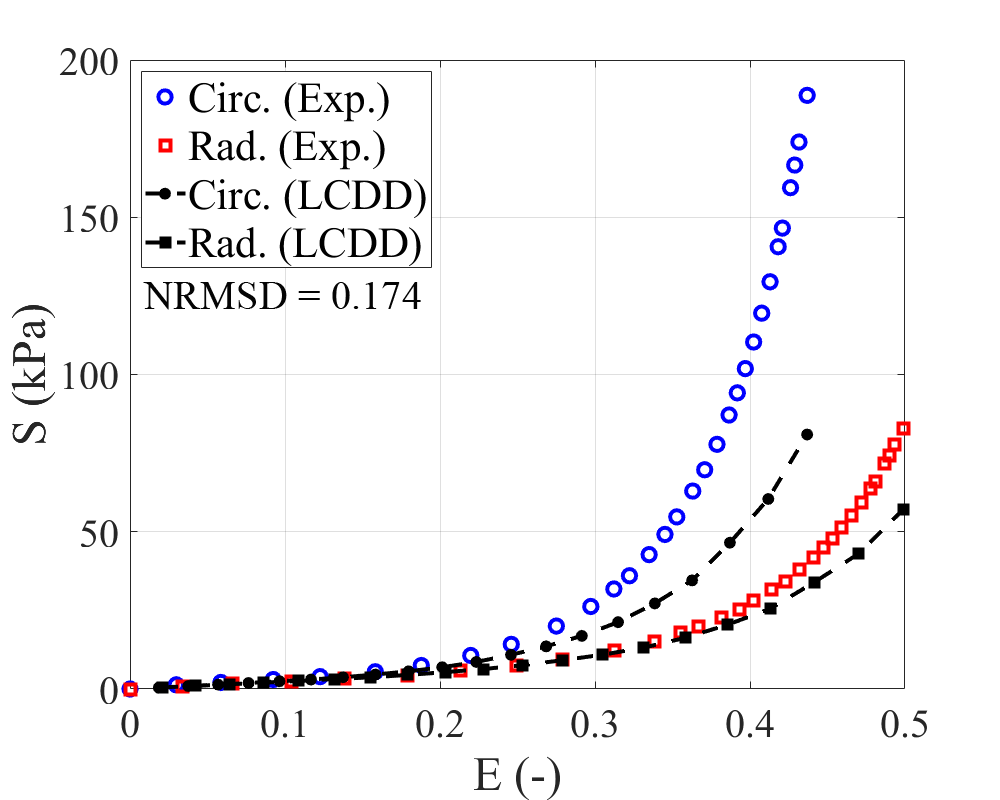}
        \caption{Protocol 4}
    \end{subfigure}
\caption{Comparison of extrapolative predictability: (a) AEDD prediction on training Protocols 1, 2, 6, 10, 11; (b) AEDD prediction on Protocol 3;  (c) AEDD prediction on Protocol 4; (d) LCDD prediction on training Protocols 1, 2, 6, 10, 11; (e) LCDD prediction on Protocol 3;  (f) LCDD prediction on Protocol 4. Protocols 1, 2, 6, 10, and 11 are used to train the autoencoder applied in AEDD}
\label{fig.bio.case3}
\end{figure}

The extrapolative prediction performance of AEDD is further explored in this case study.
Here, three biaxial tension protocols (Protocols 1, 2, and 6) with similar loading patterns, as illustrated by the experimental data in Fig. \ref{fig.bio.exp}(d) and (e), and two pure shear protocols (Protocols 10 and 11) are used for the material training dataset. The AEDD and LCDD approaches are tested on two testing protocols (Protocols 3 and 4) subjected to larger loading ratio differences between tissue’s circumferential and radial directions.
Again, Fig. \ref{fig.bio.case3} shows that AEDD outperforms LCDD in both training and testing protocols. In the testing cases (Protocols 3 and 4), while the LCDD results show clear discrepancies from the experimental data, AEDD provides a better accuracy, as evidenced by reducing the NRMSD with more than $50\%$ from the LCDD prediction. This example further verifies better extrapolative generalization capability of AEDD.

\subsubsection{Case 4}\label{sec3.2.4}
\begin{figure}[!ht]
\centering
    \begin{subfigure}{0.49\textwidth}
        \centering
        \includegraphics[width=1\linewidth]{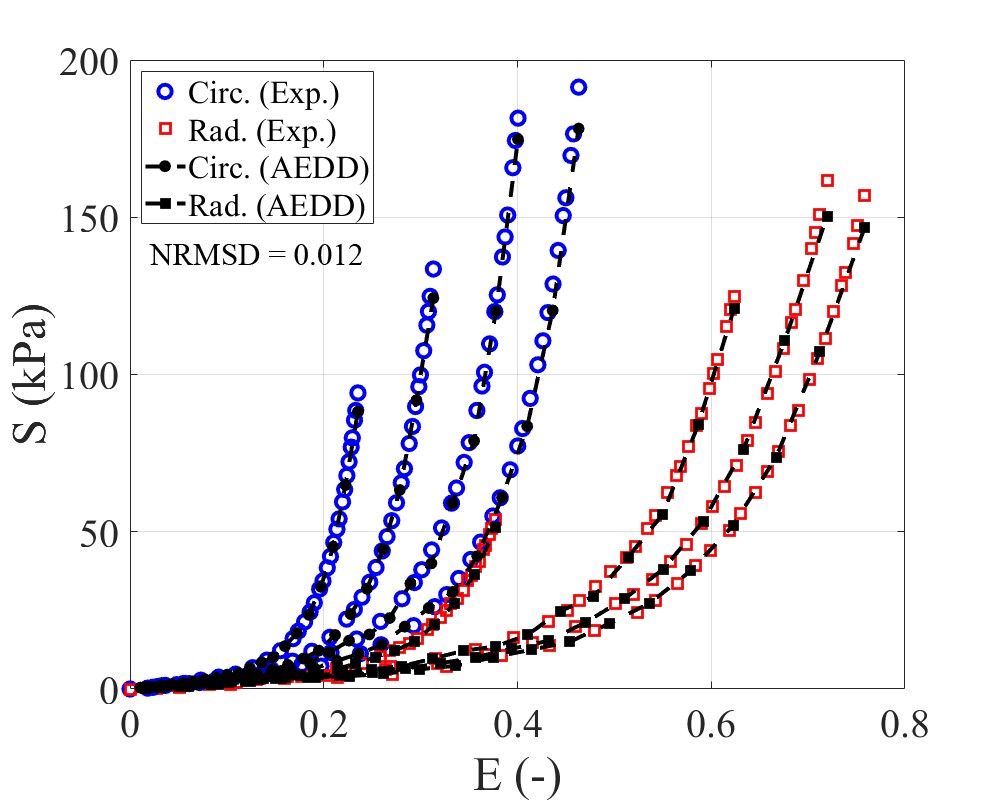}
        \caption{Protocols 2, 5, 7, 8}
    \end{subfigure}
    \begin{subfigure}{0.49\textwidth}
        \centering
        \includegraphics[width=1\linewidth]{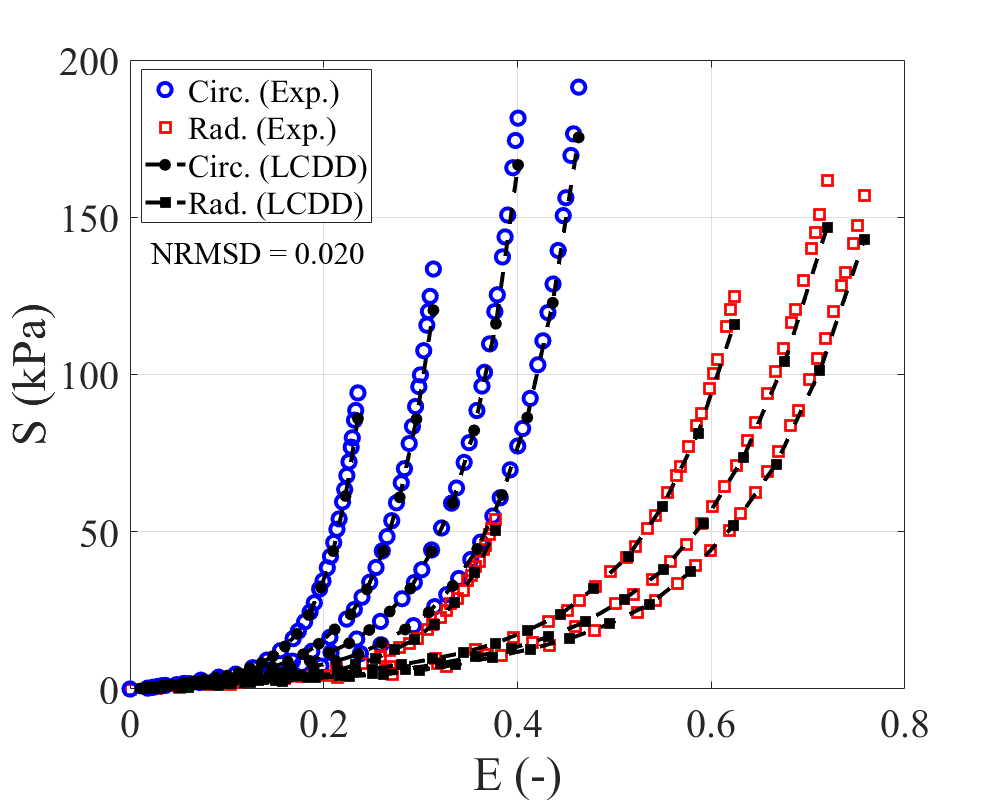}
        \caption{Protocols 2, 5, 7, 8}
    \end{subfigure}
    \begin{subfigure}{0.32\textwidth}
        \centering
        \includegraphics[width=1\linewidth]{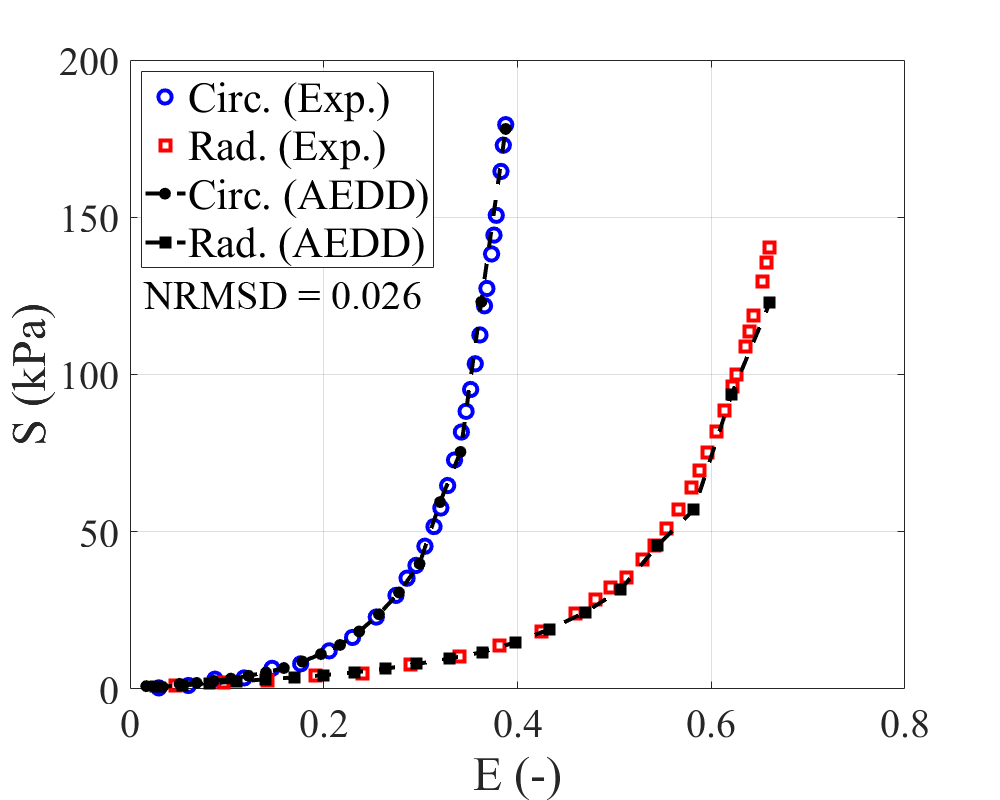}
        \caption{Protocol 1}
    \end{subfigure}
    \begin{subfigure}{0.32\textwidth}
        \centering
        \includegraphics[width=1\linewidth]{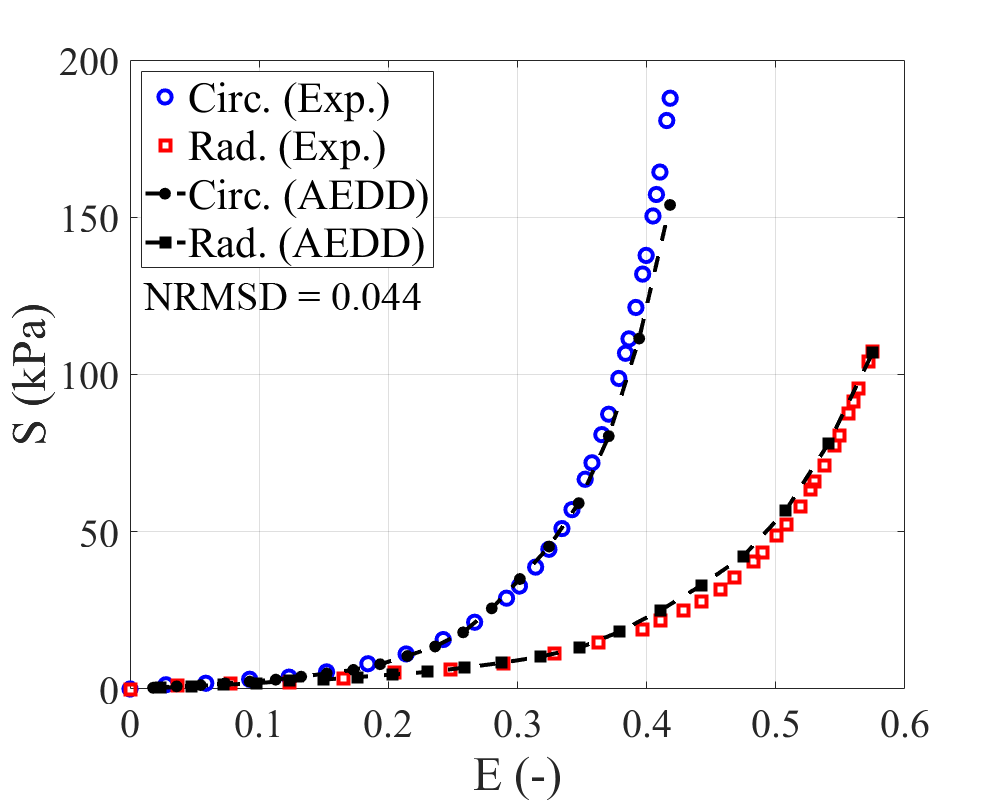}
        \caption{Protocol 3}
    \end{subfigure}
    \begin{subfigure}{0.32\textwidth}
        \centering
        \includegraphics[width=1\linewidth]{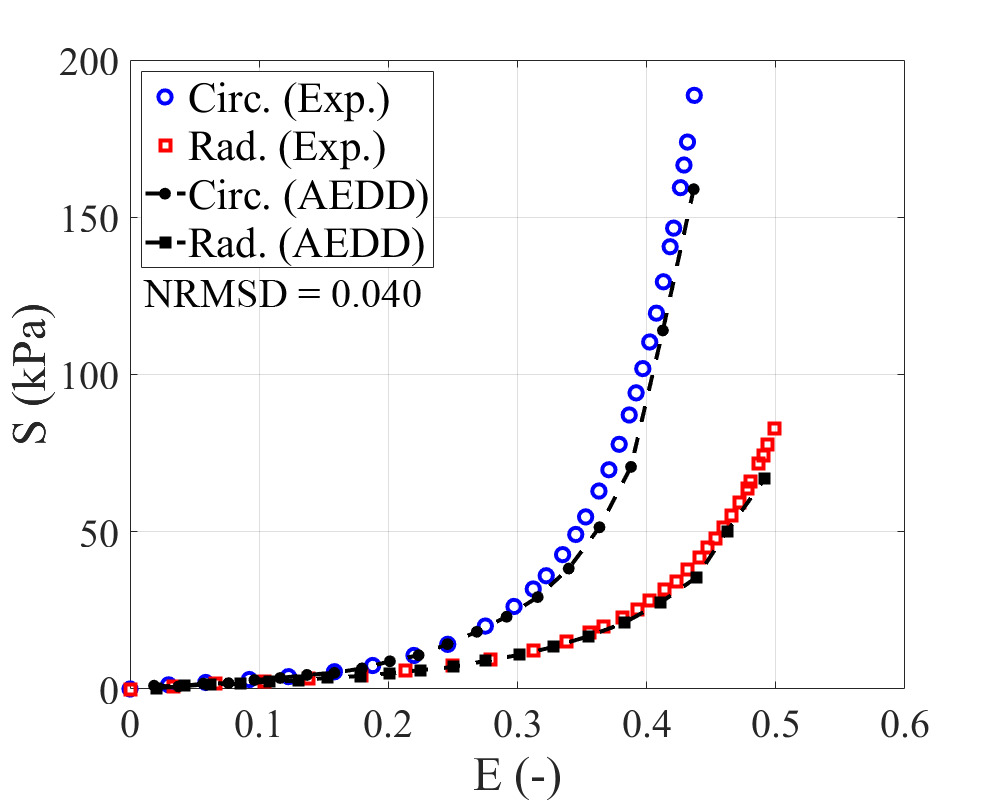}
        \caption{Protocol 4}
    \end{subfigure}
    \begin{subfigure}{0.32\textwidth}
        \centering
        \includegraphics[width=1\linewidth]{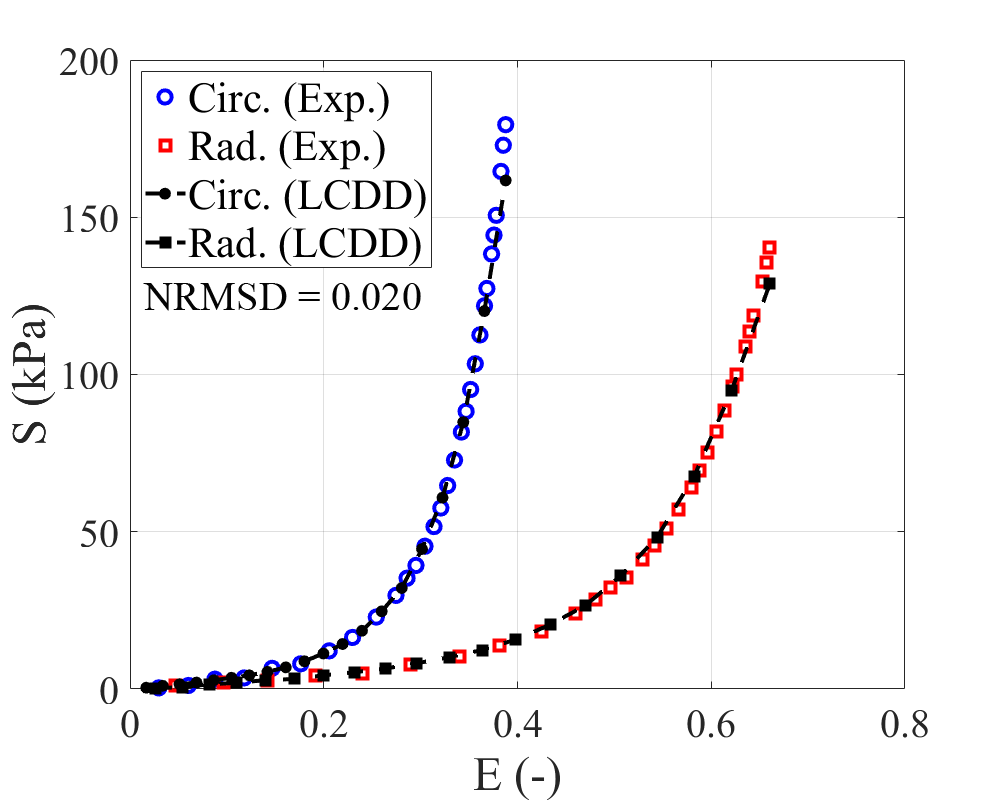}
        \caption{Protocol 1}
    \end{subfigure}
    \begin{subfigure}{0.32\textwidth}
        \centering
        \includegraphics[width=1\linewidth]{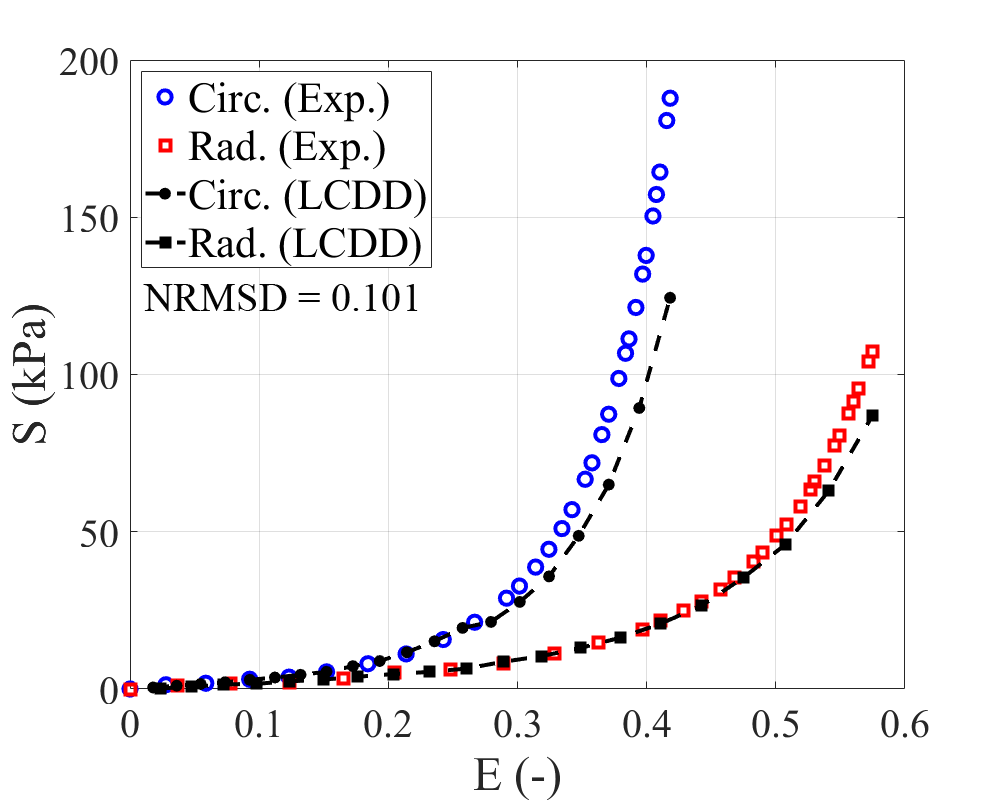}
        \caption{Protocol 3}
    \end{subfigure}
    \begin{subfigure}{0.32\textwidth}
        \centering
        \includegraphics[width=1\linewidth]{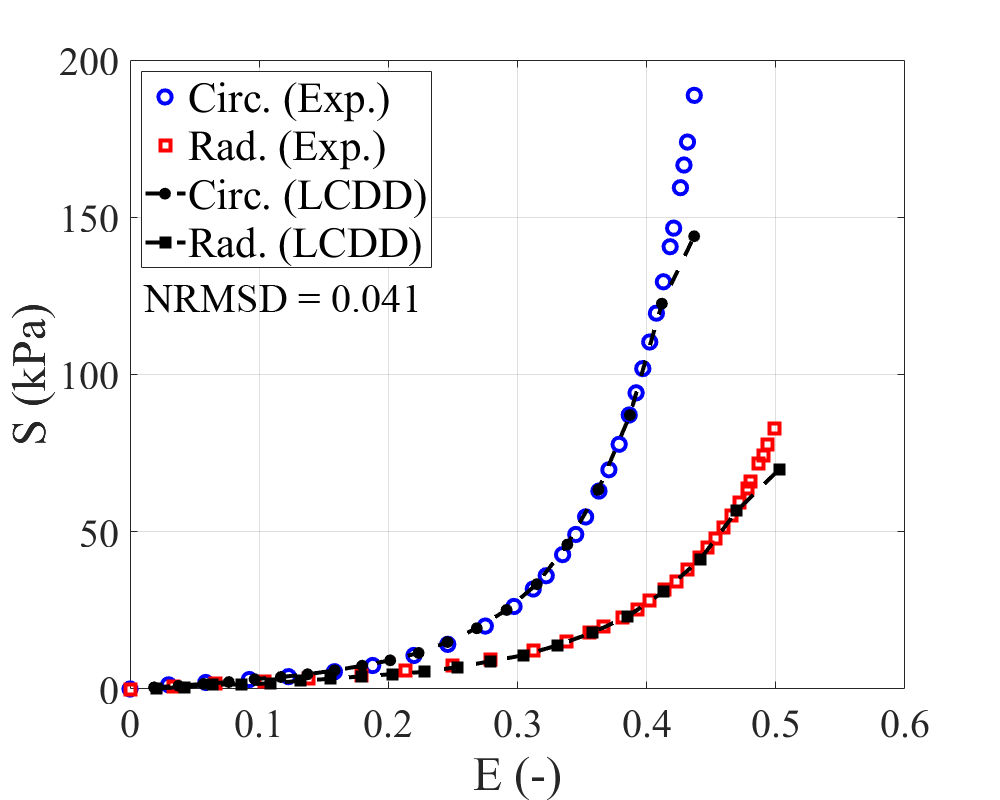}
        \caption{Protocol 4}
    \end{subfigure}
\caption{Comparison of interpolative predictability: (a) AEDD prediction on training Protocols 2, 5, 7, 8; (b) LCDD prediction on training Protocols 2, 5, 7, 8; (c) AEDD prediction on Protocol 1;  (d) AEDD prediction on Protocol 3; (e) AEDD prediction on Protocol 4; (f) LCDD prediction on Protocol 1; (g) LCDD prediction on Protocol 3;  (h) LCDD prediction on Protocol 4. Protocols 2, 5, 7, and 8 are used to train the autoencoder applied in AEDD}\label{fig.bio.case4}
\end{figure}

As can be seen from Fig. \ref{fig.bio.case1} and \ref{fig.bio.case2}, both LCDD and AEDD work well for interpolative testing cases when using training protocols with symmetrical loading conditions.
In Case 4, we investigate how a material training dataset from asymmetrically distributed protocols (biaxial tension Protocols 2, 5, 7, and 8) affects the interpolative prediction performance.
Although the simulation results on the training protocols from both AEDD and LCDD agree well with experimental data, as shown in Fig. \ref{fig.bio.case4}(a) and (b), the accuracy of LCDD deteriorates substantially on the testing protocols compared to AEDD, as shown in Fig. \ref{fig.bio.case4}(g). The results demonstrate that AEDD's performance is more robust when dealing with irregular training datasets, which could be attributed to the underlying material manifold learned by the autoencoders.


\subsubsection{Case 5}\label{sec3.2.5}
\begin{figure}[!ht]
\centering
    \begin{subfigure}{0.49\textwidth}
        \centering
        \includegraphics[width=1\linewidth]{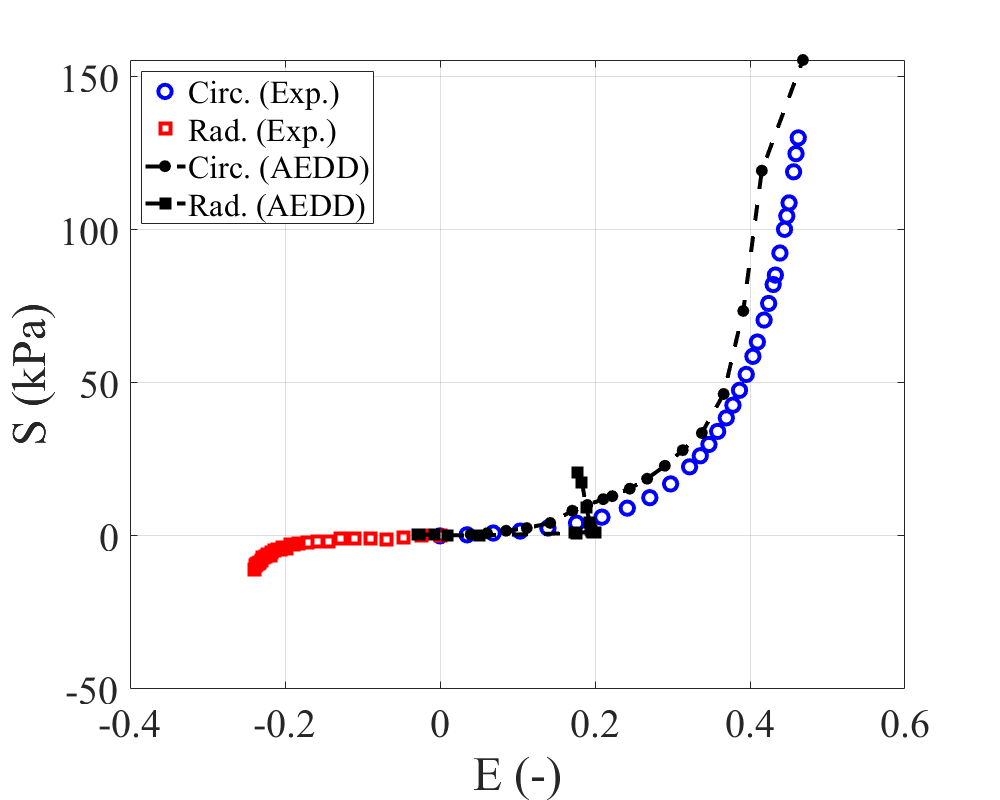}
        \caption{Protocol 10}
    \end{subfigure}
    \begin{subfigure}{0.49\textwidth}
        \centering
        \includegraphics[width=1\linewidth]{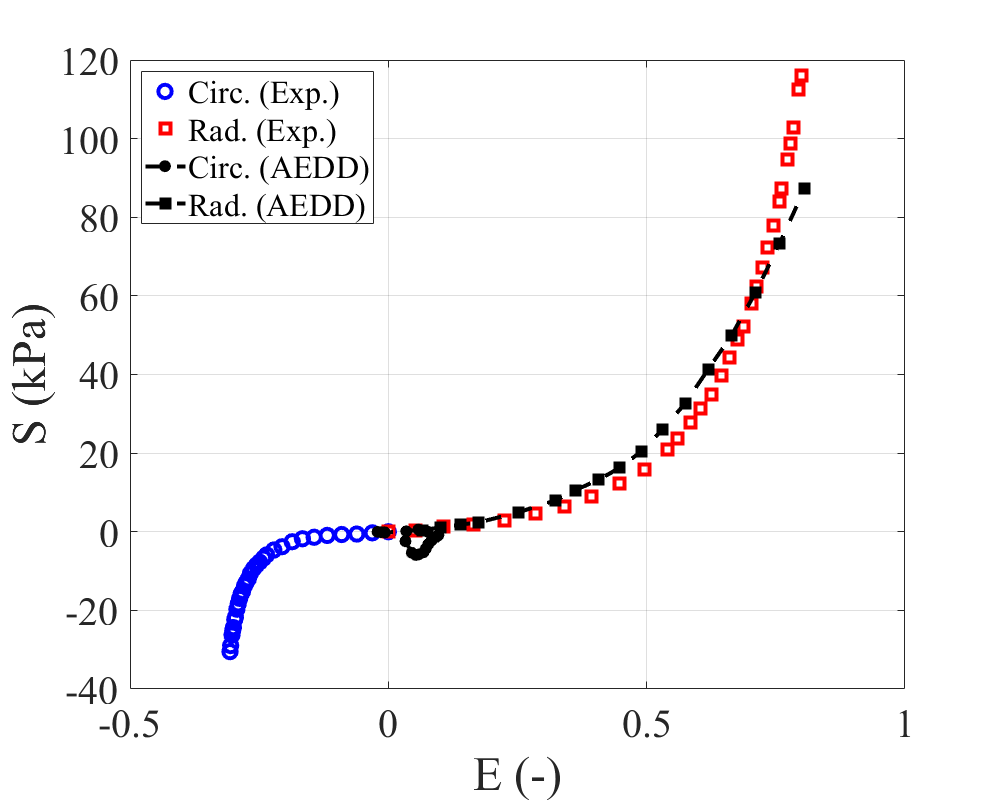}
        \caption{Protocol 11}
    \end{subfigure}
\caption{Data-driven prediction by AEDD on (a) Protocol 10 and (b) Protocol 11. Protocols 1--9 are used to train the autoencoder}
\label{fig.bio.case5}
\end{figure}

The results in Cases 1--4 have demonstrated that AEDD yields improved interpolative and extrapolative prediction compared to the LCDD approach by introducing autoencoders in the material data-driven local solver. In this last case, we investigate the performance of the AEDD method on the testing dataset that are fully unrelated to the training dataset.
Specifically, autoencoders were trained using the biaxial tension protocols 1--9 for AEDD predictions on the pure shear protocols 10 and 11. As displayed in Fig. \ref{fig.bio.case5}, AEDD predictions on the pure shear protocols (10 and 11) show some deviations from the experimental data. It is because the training protocols are all biaxial tension protocols that do not contain any information about the material behaviors in the pure shear protocols. These results demonstrate that the predictive capability of the machine learning techniques such as AEDD depends on the richness and quality of the given training data.

\section{Conclusion}\label{sec.conclusion}
In this study, we introduced the deep manifold learning approach via autoencoders to learn the underlying material data structure and incorporated it into the data-driven solver to enhance solution accuracy, generalization ability, efficiency, and robustness in data-driven computing. The proposed approach is thus named auto-embedding  data-driven  (AEDD)  computing. In this approach, autoencoders are trained in an offline stage and thus consume little computational overhead in solution procedures. 
The trained autoencoders are then applied in the proposed data-driven solver during online computation. The trained encoders and decoders define the explicit transformation between low- and high-dimensional spaces of material data, enabling efficient embedding extension to new data points. A simple Shepard convex interpolation scheme is employed in the proposed data-driven solver to preserve convexity in the local data reconstruction, enhancing the robustness of the data-driven solver.

A parametric study is conducted in the beam problem to investigate the effects of noise in material datasets, the size and sparsity of datasets, neural networks initialization during training, and autoencoder architectures on the performance of autoencdoers and data-driven solutions. Autoencoders with four different architectures are trained with synthetic noisy material datasets generated from a phenomenological model and different random initialization. The parametric study shows the performance of the offline trained autoencoders improves as the amount of training data increases regardless of the examined autoencoder architectures and neural network initialization. AEDD predictions are accurate and robust when dealing with sparse noisy datasets with the solutions converging to the constitutive model-based reference solutions as the number of material data and data density increase. In addition, with the offline trained autoencoders and efficient Shepard convex reconstruction for online computation, AEDD shows enhanced computational efficiency compared to the LCDD approach \cite{HeJBM2020}. The effectiveness of the proposed framework is further examined by modeling biological tissues using experimental data. The proposed AEDD framework shows a good performance in modeling complex materials. Through five study cases, the proposed approach show stronger generalization capability and robustness than the LCDD approach \cite{HeJBM2020}. This is attributed to the fact that the local neighbor searching and locally convex reconstruction in the proposed data-driven solver is based on geometric distance information in the filtered global embedding space learned by autoencoders, which contains the underlying manifold structure of the material data. The results of the last case also demonstrates the effects of richness and quality of the training data on the predictive capability of the AEDD method.

Although using 6 nearest neighbors in the locally convex reconstruction of the data-driven solver and an empirical weight matrix $\hat{\mathbb{C}}$ based on statistical information of data is sufficient for AEDD to produce accurate and robust solutions in the problems of this study, choosing the optimal number of nearest neighbors and the optimal weight matrix requires further investigations. The results of the proposed data-driven approach demonstrate the promising performance by integrating the autoencoder enhanced deep manifold learning into data-driven computing of systems with complex material behaviors.

\section*{Acknowledgements}
The support of this work by the National Science Foundation under Award Number CCF-1564302 to University of California, San Diego, is greatly appreciated. 
Q. H. acknowledges partial support from Pacific Northwest National Laboratory (PNNL) under the Collaboratory on Mathematics and Physics-Informed Learning Machines for Multiscale and Multiphysics Problems (PhILMs) project. 
PNNL is operated by Battelle for the DOE under Contract DE-AC05-76RL01830. 
\appendix
\section{Reproducing kernel approximation}\label{appendix:RKPM}
The displacement field $\mathbf{u(x)}$ and the Lagrange multiplier $\boldsymbol{\lambda}(\mathbf{x})$ in weak-form equations Eq. (\ref{eq.global}) are approximated by
\begin{subequations}\label{eq.2.2.1}
\begin{align}
  \begin{split}
    \mathbf{u}^h(\mathbf{x}) = \sum_{I=1}^{NP} \Psi_I(\mathbf{x}) \mathbf{d}_I,
  \end{split}
  \\
  \begin{split}
      \boldsymbol{\lambda}^h(\mathbf{x}) = \sum_{I=1}^{NP} \Psi_I(\mathbf{x}) \mathbf{\Lambda}_I,
  \end{split}
\end{align}
\end{subequations}
where $\mathbf{d}_I$ and $\mathbf{\Lambda}_I$ are the nodal coefficients associated with the fields $\mathbf{u(x)}$ and  $\boldsymbol{\lambda}(\mathbf{x})$, respectively, and $\Psi_I(\mathbf{x})$ is the reproducing kernel (RK) approximation function expressed as
\begin{equation}\label{eq.2.2.1b}
    \Psi_I(\mathbf{x}) = \mathbf{H}^T(\mathbf{x}-\mathbf{x}_I) \mathbf{b(x)} \phi_a(\mathbf{x}-\mathbf{x}_I),
\end{equation}
where $\mathbf{H}^T(\mathbf{x}-\mathbf{x}_I) = [1, x_1-x_{1I},x_2-x_{2I},x_3-x_{3I}, ..., (x_3-x_{3I})^n]$ is a vector of monomial basis functions up to the $n$-th order, and $\phi_a(\mathbf{x}-\mathbf{x}_I)$ is a kernel function with a local support size "$a$", controlling the smoothness of the RK approximation function,
for example, the cubic B-spline kernel function:
\begin{equation}\label{eq.2.2.1c}
    \phi_a(y) = \begin{cases}
        \frac{2}{3} - 4y^2 + 4y^3, & 0 \leq y < \frac{1}{2} \\
        \frac{4}{3} - 4y + 4y^2 - \frac{4}{3}y^3, & \frac{1}{2} \leq y < 1\\
        0, & y \geq 1
    \end{cases} \hspace{0.5cm} \text{with} \hspace{0.1cm} y = \frac{||\mathbf{x}-\mathbf{x}_I||}{a}.
\end{equation}
In Eq. (\ref{eq.2.2.1b}), $\mathbf{b(x)}$ is a parameter vector determined by imposing the $n$-th order reproducing conditions \cite{liu1995reproducing,chen1996reproducing},
\begin{equation}\label{eq.2.2.1d}
    \sum_{I=1}^{NP} \Psi_I(\mathbf{x}) x_{1I}^i x_{2I}^j x_{3I}^k = x_{1}^i x_{2}^j x_{3}^k, \hspace{0.5cm} |i+j+k|=0,1,...,n.
\end{equation}
Substituting Eq. (\ref{eq.2.2.1b}) into Eq. (\ref{eq.2.2.1d}) yields $\mathbf{b(x)}=\mathbf{M}^{-1}(\mathbf{x}) \mathbf{H}(\mathbf{0})$, where  $\mathbf{M}(\mathbf{x})$ is a moment matrix given by
\begin{equation}\label{eq.2.2.3}
    \mathbf{M}(\mathbf{x}) = \sum_{I=1}^{NP} \mathbf{H}(\mathbf{x}-\mathbf{x}_I) \mathbf{H}^T(\mathbf{x}-\mathbf{x}_I) \phi_a(\mathbf{x}-\mathbf{x}_I).
\end{equation}
The RK approximation function is then obtained as,
\begin{equation}\label{eq.2.2.2}
    \Psi_I(\mathbf{x}) = \mathbf{H}^T(\mathbf{0}) \mathbf{M}^{-1}(\mathbf{x}) \mathbf{H}(\mathbf{x}-\mathbf{x}_I) \phi_a(\mathbf{x}-\mathbf{x}_I).
\end{equation}

\section{Nodal integration scheme}\label{appendix:SCNI}
The stabilized conforming nodal integration (SCNI) approach is employed for the domain integration of the weak form (Eq. (\ref{eq.global})) to achieve computational efficiency and accuracy when using RK shape functions with nodal integration quadrature schemes.

\begin{figure}[!h]
    \centering
    \includegraphics[width=0.4\linewidth]{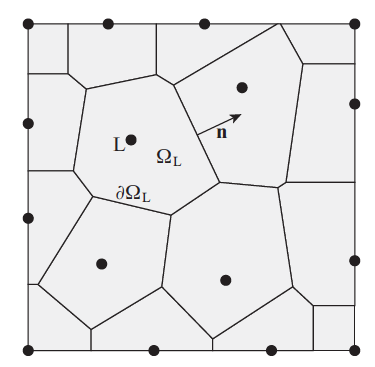}
    \caption{Illustration of Voronoi diagram for SCNI.}\label{fig.SCNI}
\end{figure}

The key idea behind SCNI is to satisfy the linear patch test (thus, ensure the linear consistency) by leveraging a condition, i.e. the divergence constraint on the test function space and numerical integration \cite{chen2002non}, expressed as:
\begin{equation}\label{eq.integration_constraint}
    \hat{\int_{\Omega}}\nabla\Psi_I d\Omega = \hat{\int_{\partial\Omega}} \Psi_I \mathbf{n} d\Gamma, 
\end{equation}
where '\^{}' over the integral symbol denotes numerical integration. In SCNI, an effective way to achieve Eq. (\ref{eq.integration_constraint}) is based on nodal integration with gradients smoothed over conforming representative nodal domains, as shown in Fig. \ref{fig.SCNI}, converted to boundary integration using the divergence theorem
\begin{equation}\label{eq.smooth_gradient}
    \Tilde{\nabla}\Psi_I(\mathbf{x}_L) = \frac{1}{V_L} \int_{\Omega_L} \nabla\Psi_I d\Omega = \frac{1}{V_L} \int_{\partial\Omega_L} \Delta\Psi_I \mathbf{n}d\Gamma,
\end{equation}
where $V_L = \int_{\Omega_L} d\Omega$ is the volume of a conforming smoothing domain associated with the node $\mathbf{x}_L$, and $\Tilde{\nabla}$ denotes the smoothed gradient operator. In this method, smoothed gradients are employed for both test and trial functions, as the approximation in Eq. (\ref{eq.smooth_gradient}) enjoys first order completeness and leads to a quadratic rate of convergence for solving linear solid problems by meshfree Galerkin methods. As shown in Fig. \ref{fig.SCNI}, the continuum domain $\Omega$ is partitioned into $N$ conforming cells by Voronoi diagram, and both the nodal displacement vectors and the state variables (e.g., stress, strain) are defined at the set of nodes $\{\mathbf{x}_L\}_{L=1}^N$.

Therefore, if we consider two-dimensional elasticity problem under the SCNI framework, the smoothed strain-displacement matrix $\Tilde{\mathbf{B}}_I(\mathbf{x}_L)$ used in (16) is expressed as:
\begin{equation}\label{eq.smooth_gradient_matrix}
    \Tilde{\mathbf{B}}_I(\mathbf{x}_L) =
    \begin{bmatrix}
        \Tilde{b}_{I1}(\mathbf{x}_L) & 0 \\
        0 & \Tilde{b}_{I2}(\mathbf{x}_L) \\
        \Tilde{b}_{I2}(\mathbf{x}_L) & \Tilde{b}_{I1}(\mathbf{x}_L)
    \end{bmatrix},
\end{equation}
with
\begin{equation}\label{eq.smooth_gradient_component}
    \Tilde{b}_{Ii}(\mathbf{x}_L) = \frac{1}{V_L} \int_{\partial\Omega_L} \Psi_I(\mathbf{x})n_i(\mathbf{x}) d\Gamma.
\end{equation}
Since the employment of the smoothed gradient operator in Eq. (\ref{eq.smooth_gradient}) and Eq. (\ref{eq.smooth_gradient_component}) satisfies the divergence constraint regardless of the numerical boundary integration, a trapezoidal rule for each segment of $\partial\Omega_L$ is used in this study.

\bibliography{Reference}

\end{document}